\documentclass[11pt,reqno]{amsart}
\usepackage{color}
\usepackage{amsfonts}
\usepackage{mathrsfs}
\usepackage{amssymb, amsmath}
\usepackage{graphicx}
\usepackage{pstricks}
\usepackage{epsfig,subfigure}
\usepackage{color}

\newcommand{\newcom}{\newcommand}

\newcom{\al}{\alpha}
\newcom{\be}{\beta}
\newcom{\eps}{\epsilon}
\newcom{\del}{\delta}
\newcom{\Del}{\Delta}
\newcom{\ga}{\gamma}
\newcom{\Ga}{\Gamma}
\newcom{\ka}{\kappa}
\newcom{\Lam}{\Lambda}
\newcom{\lam}{\lambda}
\newcom{\Om}{\Omega}
\newcom{\om}{\omega}
\newcom{\Si}{\Sigma}
\newcom{\si}{\sigma}
\newcom{\tht}{\theta}
\newcom{\dtri}{\nabla}
\newcom{\td}{\tilde}
\newcom{\tri}{\triangle}
\newcom{\oo}{\infty}
\newcom{\vphi}{\varphi}

\newcom{\cA}{{\mathcal A}}
\newcom{\cB}{{\mathcal B}}
\newcom{\cC}{{\mathcal C}}
\newcom{\cD}{{\mathcal D}}
\newcom{\cF}{{\mathcal F}}
\newcom{\cI}{{\mathcal I}}
\newcom{\cL}{{\mathcal L}}
\newcom{\cM}{{\mathcal M}}
\newcom{\cN}{{\mathcal N}}
\newcom{\cP}{{\mathcal P}}
\newcom{\cR}{{\mathcal R}}
\newcom{\cS}{{\mathcal S}}
\newcom{\cQ}{{\mathcal Q}}
\newcom{\caly}{{\mathcal Y}}
\newcom{\calz}{{\mathcal Z}}
\newcom{\bfz}{{\bf Z}}
\newcom{\R}{\Bbb R}
\newcom{\N}{\Bbb N}
\newcom{\Z}{\Bbb Z}
\newcom{\C}{\Bbb C}
\newcom{\E}{\Bbb E}

\newcom{\bx}{\bar x}
\newcom{\bz}{\bar z}
\newcom{\f}{\frac}
\newcom{\di}{\displaystyle\int}
\newcom{\ds}{\displaystyle\sum}
\newcom{\dl}{\displaystyle\lim}
\newcom{\ov}{\overline}
\newcom{\sset}{\subset}
\newcom{\wt}{\widetilde}
\newcom{\p}{\partial}
\newcom\na{\nabla}
\newcom{\co}{\cdot}
\newcom{\suml}{\sum\limits}
\newcom{\supl}{\sup\limits}
\newcom{\intl}{\int\limits}
\newcom{\infl}{\inf\limits}
\newcom{\disp}{\displaystyle}
\newcom{\non}{\nonumber}
\newcom{\no}{\noindent}
\newcom{\QED}{$\square$}
\def\ef{\hphantom{MM}\hfill\llap{$\square$}\goodbreak}

\newtheorem{athm}{\bf \t}[section]
\newenvironment{thm} [1] {\def\t{#1}\begin{athm} \bf \rm} {\end {athm}}
\newcom{\bthm}{\begin{thm}}\newcom{\ethm}{\end{thm}}

\newcom{\beq}{\begin{equation}}
\newcom{\eeq}{\end{equation}}
\newcom{\ben}{\begin{eqnarray}}
\newcom{\een}{\end{eqnarray}}
\newcom{\beno}{\begin{eqnarray*}}
\newcom{\eeno}{\end{eqnarray*}}

\numberwithin{equation}{section}

\begin{document}

\title{Elliptic estimates for Dirichlet-Neumann operator on a corner domain}

\author{Mei MING}
\address{School of Mathematics\\ Sun Yat-sen Univeristy\\ Guangzhou 510275,China }
\email{mingm@mail.sysu.edu.cn}

\author{Chao WANG}
\address{School of Mathematical Sciences\\ Peking University\\ Beijing 100871,China}
\email{wangchao@math.pku.edu.cn}

\begin{abstract} We consider the elliptic estimates for Dirichlet-Neumann operator related to the water-wave problem on a two-dimensional corner domain in this paper. Due to  the singularity of the boundary, there will be singular parts in the solution of the elliptic problem  for D-N operator.  To begin with,  we  study  elliptic problems with mixed boundary condition to derive singularity decompositions and estimates. Based on the analysis,   we present the estimates for both  D-N operator and its shape derivative with the existence of singular parts. 
\end{abstract}
\maketitle

\tableofcontents

\section{Introduction}
The water-wave problem investigates the movements of the ideal inviscid incompressible fluid in some domain with a free surface. Arising from the water-wave problem, Dirichlet-Neumann (D-N) operator  plays an important role in the research. 

We want to study D-N operator related to the 2D water-wave problem on a corner domain, where the free surface intersects with the oblique flat beach at the contact point. Only when the properties of D-N operator  are well understood,  the water-wave problem on a corner domain can be solved. In fact, this problem corresponds to the scene of sea waves moving near a beach,  therefore a detailed study on this problem is interesting and important both in mathematical theory and physical applications.  

\begin{figure}
\begin{pspicture}(0,-3)(5,2.5)
\pscurve[showpoints=false](-0.45,-0.17)(0,-0.03)(0.8,0.1)(1.2,0.15)(2,0.03)(3,0.1)(4,0)(5,0.15)(6,0)(7,0.02)
\psline[linewidth=1.5pt](-2,0.5)(4,-2.15)
\pscurve[linewidth=1.5pt](4,-2.15)(4.5,-2.23)(5,-2.2)(5.7,-2.12)(6.5,-2.2)(7.2,-2.13)
\rput(3.1,1){$\vec n_t$}
\psline{->}(3,0.1)(2.97,0.7)
\rput(2.3,0.5){$\vec\tau_t$}
\psline{->}(3,0.1)(2.3,0.096)
\rput(1.8,-2){$\vec n_b$}
\psline{->}(2.5,-1.5)(2.3,-2)
\rput(2.8,-2){$\vec \tau_b$}
\psline{->}(2.5,-1.5)(3,-1.75)
\rput(6.5,-1.2){$\vec g$}
\psline{->}(6.2,-.7)(6.2,-1.3)
\psline(1,-0.5)(2,-0.5)
\psline(3,-0.5)(3.5,-0.5)
\psline(5,-0.5)(6,-0.5)
\psline(2.5,-1)(3,-1)
\psline(4.5,-1)(5.5,-1)
\psline(4,-1.5)(4.5,-1.5)
\psline(5.5,-2)(6,-2)
\rput(-0.7,-0.6){$O(0,0)$}
\rput(5.5,-1.5){$\Omega_t$}
\rput(5.3,0.7){$\Ga_t:\,z=\eta(t,x)$}
\rput(5,-2.8){$\Ga_b:\,z=l(x)$}
\rput(0.7, -1.4){$\Ga_b:\,z=-\ga x$}
\end{pspicture}
\caption{The corner domain}
\label{figdomain}
\end{figure}
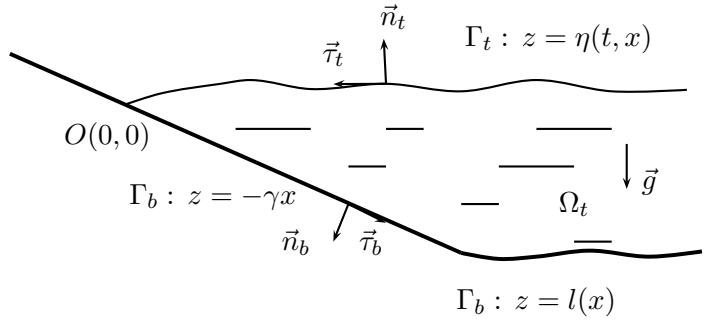

We consider in this paper the (horizontally) unbounded corner domain $\Omega$ containning a free surface $\Ga_t=\{(x,z)\,|\,z=\eta(t,x)\}$ and a fixed smooth and bounded bottom $\Ga_b=\{(x,z)\,|\,z=l(x)\}$, that is
\[\Omega=\{(x,z)\,|\,l(x)<z<\eta(t,x),\,x\ge 0\}\] where $\eta(t,x)$ is the deviation of the moving free surface near its equilibrium state $z=0$ at time $t$ and 
the bottom 
\[l(x)=-\ga x,\quad \hbox{when}\quad x\le x_0\] for some fixed constant $x_0$ and slope
$-\ga$. We write $\eta(x)$ for $\eta(t,x)$ in this paper since {\it we only consider about D-N operator for some fixed time $t$ here}. The dependence on time $t$ for D-N operator and the elliptic system will be discussed in a subsequent paper. 

Through this paper we consider that $\eta(x)$ is defined on $\R^+$, and the  surface and the bottom intersects only at the origin $O=(0,0)$ that is,  $\eta(x)$ satisfies
\[
\eta(0)=0,\quad\hbox{and}\quad \eta(x)-l(x)>0,\quad\forall x>0.
\]
Based on the corner domain $\Om$, the (non-scaled) D-N operator $G(\eta)$ related to the water-wave problem is defined by
\[
G(\eta)\psi=\p_{n_t}\phi|_{\Ga_t}
\]
with the associated elliptic problem for the velocity potential $\phi$
\[
\left\{\begin{array}{ll}
\Delta\phi=0,\qquad \hbox{in}\quad \Omega\\
\phi|_{\Ga_t}=\psi,\qquad \p_{n_b}\phi|_{\Ga_b}=0.
\end{array}\right.
\] Sometimes we also use the definition of the scaled D-N operator
\[
G(\eta)\psi=\sqrt{1+|\eta'|^2}\p_{n_t}\phi|_{\Ga_t}.
\]The study for D-N operator is essentially the study of the elliptic system above. Since there are  mixed-boundary conditions  in this system, we will focus only on the mixed problem in this paper.   

Before stating our main result, we should recall some previous results about general linear elliptic systems on non-smooth domains first. To begin with, when we say non-smooth, we always refer to Lipschitz. {\it As a matter of fact, when the boundary $\Ga$ of the domain $\Om$ is Lipschitz, the classical elliptic theory for a smoother boundary doesn't work any more because it's impossible to straighten the non-smooth boundary by using a smooth transformation. } A theory for non-smooth boundaries is therefore required. In fact,  the non-smooth elliptic theory has been fully developed already, starting from the fundamental works by Kondrat'ev \cite{Kon63, Kon67}. Some more first works are done by Birman and Skvortsov \cite{BS}, Eskin \cite{E}, Lopatinskiy \cite{Lop},  Maz'ya \cite{Maz64, Maz67}, Kondrat'ev and Oleinik \cite{KO}, Maz'ya and Plamenevskiy \cite{MazP77}, Maz'ya and Rossmann \cite{MazR92},  Grisvard \cite{PG1}, Dauge \cite{D},  Costabel and Dauge \cite{CD93I, CD93II} etc.. These works provide regularity results in Sobolev space or  some weighted Sobolev space for general linear elliptic problems on Lipschitz domains  under assumptions of sufficient smoothness for both the boundary away from the corner point and the system coefficients. In these works the authors considered the Dirichlet-boundary conditions or lower-order regularities most of the time.

In  a classical elliptic problem,   the solution can be as smooth as possible when the boundary and the right side of the problem are smooth enough. As a contrast, when the boundary is Lipschitz (piece-wisely smooth), there would be corner points on the boundary (singular points). Consequently, the violation of the boundary smoothness will lead to  singularities in the solutions near  corner points, no matter how smooth the right sides of the problems may become, see \cite{Fich}.

On the other hand, there always exists a solution in $H^1$ for an elliptic problem on a Lipschitz domain, see for example Grisvard \cite{PG1}.  Singularities arise when  the higher-order regularity is referred to.  Kondrat'ev \cite{Kon67} proved the decomposition of the solution into a regular part and a singular part, where the singular part is indeed  a finite sum of some singular functions with constant coefficients. Maz'ja and Plamenevskii \cite{MazP84}, Dauge, Nicaise, Bourlard and Lubuma \cite{DNBL}, Grisvard \cite{PG2} etc. gave formulas for the singular coefficients in different cases respectively. In fact, the singular functions have the general ansatz 
\[r^\lam\log^q r\varphi(\tht)\]  near the corner point, where $\lam$, $q$ are constants and $\varphi(\tht)$ is a bounded trigonometric function. These singular functions depend on the shape of the corner as well as  the left-side operators of the elliptic system including the boundary operators instead of the right side of the system.  Moreover, the number of singular functions gets larger when a higher regularity of the solution is considered.   Therefore the singular part is also named as {\it asymptotics}  sometimes. Another fact is,  as mentioned for example in \cite{PG2}, the smaller the contact angle becomes, the smoother the solution would become (i.e. the number of  singular functions will decrease), and a more detailed discussion will be in Rmk \ref{critical angle}.  In this paper, the singular functions would be in form of $r^\lam \varphi(\tht)$ due to the range of the contact angle $\om$ concerned.

The estimates for the regular part and the singular coefficients in Sobolev space or $L^p$ space are proved in Maz'ya and Plamenevskiy \cite{MazP80}, Grisvard \cite{PG1}, Costabel and Dauge \cite{CD93I, CD93II} etc. for 2 or 3 dimensional domains. In Grisvard \cite{PG1,PG2}, the $H^2$ estimate  and the formula for the singular coefficients are given in a simple case when the domain is a 2 dimensional rectilinear polygon and the elliptic operator is Laplacian. Costabel and Dauge \cite{CD93I, CD93II} considered 3 dimensional bounded piecewise analytic Lipschitz domain with a general 2nd-order elliptic operator and some general boundary conditions. They proved higher-order estimates for the regular parts of the solutions in Sobolev space, as well as estimates for the singular coefficients. One can tell from  \cite{CD93I, CD93II} that, when the boundaries are curved, the situation is much more complicated.

Since only mixed-type boundary conditions are considered in this paper, we recall some results on  mixed problems. Grisvard \cite{PG1,PG2} proved the singular decompositions for general boundary conditions including the mixed type, in the case of the 2 dimensional rectilinear boundary and the Laplacian operator. Banasiak and Roach \cite{B, BR1} showed $H^2$ decompositions for the curvilinear boundary and general 2nd-order elliptic operators. Dauge \cite{D92} studied mixed problems with  general 2nd-order elliptic operators in n dimensional polyhedron, where the boundary is of piecewise $\cC^\rho$ type and the author proved the $L^p$ regularity result without estimates. Costabel and Dauge \cite{CD94} talked about the mixed boundary case but they omit it for the sake of clarity, while as mentioned before, in \cite{CD93II} they proved an estimate for general boundary condition with the boundary piecewise analytic. 

Compared to these works, we consider in this paper the mixed boundary problem with a limited smoothness on the boundary, and {\it we  trace the dependence of the solution on the regularity  of the upper boundary in Sobolev norms in order to study D-N operator}, which is not  clearly stated in the works mentioned above. Moreover, based on the previous results \cite{PG1, PG2, BR1}, we prove the higher-order decompositions  in a detailed way.  In fact, starting from the unique existence of the variational solution $u\in H^1(\Om)$ for the mixed problem, the solution can be decomposed as 
\[u=u_r+\sum_i \chi_i(\om)c_{si}S_i\] with $u_r$  the regular part, $\chi_i(\om)$ some characteristic function for $\om$, $c_{si}$ the singular coefficient, and $S_i$ the singular function with some explicit formula. The decomposition depends on the order of the regularity we need for the solution. As stated above, when the order increases, the number of singular functions also increases.  On the other hand,  we  also consider higher-order estimates for the regular part as well as for singular coefficients, where we trace the dependence on the free boundary function $\eta$. The estimates for regular parts are similar as classical elliptic estimates. 

After studying the properties for mixed problems, we turn to  D-N operator, which is  an order-1 elliptic operator and plays an important role in the water-wave theory,  as already showed in  some previous works on water-wave problems essentially in smooth domains: Zakharov \cite{Z}, Craig, Sulem and Sulem \cite{CSS}, Wu \cite{Wu1}, Lannes \cite{Lannes}, Shatah and Zeng \cite{ShZeng}, Alazard, Burq and Zuily \cite{ABZ}. In \cite{Lannes,ShZeng}, the authors considered  the estimate for D-N operator under different formulations. Meanwhile, the expression for the shape derivative is given, and the estimate for the shape derivative is also proved. In fact, the shape derivative is important in the energy estimates for the linearized water-wave problem.

Similarly as in \cite{Lannes},  in this paper we will firstly give the estimate for D-N operator in presence of a singular part from the related elliptic system. Indeed, corresponding to the decomposition in the elliptic system, the estimates for D-N operator can be decomposed naturally into singular parts and regular parts.   Secondly, we plan to derive the expression of the shape derivative for the corner-domain case, which is motivated by \cite{ShZeng}. Consequently,  an estimate for the shape derivative will also be proved.  Moreover, a special case for the shape derivative is also discussed, which has a similar form as in \cite{Lannes}.

Compared to the classical works essentially on smooth domains,  water-wave problems with boundary singularities (corners) is a completely new topic in recent years. In fact, there are  works done quite recently by  Kinsey and Wu \cite{WuK}, Wu \cite{Wu2, Wu3} on 2 dimensional water-wave problems with angled crests. In these papers, the authors give a priori estimates and prove the well-posedness for the water-wave problem. For D-N operator part, they use a conformal mapping to convert the boundary singularity to the singularity for the mapping itself and they analyze the properties of this conformal mapping. Alazard \cite{A} considers three-dimensional fluid in a rectangular tank, bounded by a flat bottom, vertical walls and a free surface evolving under the influence of gravity. He proves that one can estimate its energy by looking only at the motion of contact points between the free surface and the vertical walls. 

Our results in this paper would be different from the classical works on D-N operators and water-wave problems, since we consider the problem on a corner domain. Meanwhile, our work is also in a  different formulation from the results  by \cite{WuK, Wu2, Wu3} for two reasons. Firstly,   our domain has an oblique bottom (which is different from the vertical walls in \cite{A} as well), so the free surface couldn't be naturally symmetrized with respect to the contact point. Secondly, we keep the boundary singularity and give complete estimates for D-N operators.  

Now, we state the main results of this paper about  estimates and decompositions for D-N operator and its shape derivative, and  one can find  more  details in Section 6. We consider the contact angle $\om\in (0,\f\pi2)$ through our paper, see Section 3. In the following results, we always try to split the singular part away from the regular part, which leads to some extra terms compared to \cite{Lannes, ShZeng}.

\bthm{Theorem} {\it If $\eta\in H^{K+2}(\R^+)$ with $K\in\N$ and  $\om\neq \f\pi{2n}$  for some $n\in\N$, then for all $f\in H^{K+\f32}(\R^+)$, one  can have the estimate for the  (scaled) D-N operator as
\beno
\big|G(\eta)f-\sqrt{1+|\eta'|^2}\p_{n_t} u_{sK}|_{\Gamma_t}\big|_{H^{K+\f12}}\leq C( |\eta|_{H^{K+2}})|f|_{H^{K+\f32}}.
\eeno where $f$ is the Dirichlet boundary condition and $u_{sK}$ is the singular part arising from the related elliptic problem \eqref{DN elliptic system} for $u$.
}\ethm
The shape derivative for D-N operator is derived following the idea of \cite{ShZeng}, and the estimate is given based on our elliptic-decomposition theory. The theorem below is proved in Section 6.2.1.
 \bthm{Theorem}\label{DN op geometric}
{\it Let $\om \in (0,\,\f\pi2)$ be the contact angle  satisfying $\om\neq \f\pi{2n}$ for any $n\in \N$. If $\eta\in H^{K+2}(\R^+)$ and  $f\in H^{K+\f32}(\Ga_t)$,  the shape derivative for  the (non-scaled) D-N operator related to system \eqref{DN elliptic system} of $u$ can be decomposed into
\[
D_sG(\eta_s)f_s|_{s=0}-G(\eta)\big(D_sf_s|_{s=0}\big)=\na_{n_t} u_2-\na_{n_t}{\bf w}\cdot\na u-\na_{(\na u)^T}{\bf w}\cdot{\bf n}_t\big|_{\Ga_t}=l_r+l_s,
\] 
where $\bf w$ is the variational vector field on $\Ga_t$, $D_s$ is the variational material derivative with the parameter $s$, $u_2$ satisfies system \eqref{u2 elliptic system} with the singular decomposition $u_2=u_{2rK}+u_{2sK}$. 

Moreover, the regular part $l_r\in H^{K+\f12}(\Ga_t)$ and the singular part $l_s\in H^{\f12}(\Ga_t)$ are defined by
\[\begin{split}
&l_r=\na_{n_t} u_{2rK}-\na_{n_t}{\bf w}\cdot\na u_{rK}-\na_{(\na u_{rK})^T}{\bf w}\cdot{\bf n}_t\big|_{\Ga_t},\\
&l_s=\na_{n_t} u_{2sK}-\na_{n_t}{\bf w}\cdot\na u_{sK}-\na_{(\na u_{sK})^T}{\bf w}\cdot{\bf n}_t\big|_{\Ga_t}
\end{split}\] and the following  estimate holds:
\[
|l_r|_{H^{K+\f12}}+|l_s|_{H^{\f12}}\le C\big(|\eta|_{H^{K+2}}\big)|f|_{H^{K+\f32}}|{\bf w}|_{H^{K+\f32}}
\] with the constant $C$ depending on $|\eta|_{H^{K+2}}$ and $\Ga_b$. 
}
\ethm
In the case of infinite depth (the bottom $\Ga_b$ is simply $z=-\ga x$) and assuming that $T_S$ near the corner is a global straightening transformation, we also prove a formal expression for the shape derivative in the following theorem, which is similar as the shape derivative in \cite{Lannes}.
 \bthm{Theorem}\label{thm:D-N-S} 
{\it For a given Dirichlet boundary condition $f$ on ${\Ga_t}$,  the (formal) shape derivative of the scaled D-N operator $G(\eta)$ related to system \eqref{DN elliptic system} is 
\[
d_\eta G(\eta)f\cdot h=G(\eta)\Big(\big(r_1+r_2\bar B\big)f'+r_3\bar B G(\eta)f\Big)+\Big(r_3\bar B f'+\big(r_1-r_2\bar B\big)G(\eta)f\Big)'
\] where $\bar B$ is expressed in \eqref{bar B} and 
\[
r_1=\f{h}{\ga+\eta'},\quad r_2=\f{\ga\eta'-1}{1+(\eta')^2},\quad r_3=\f{\ga+\eta'}{1+(\eta')^2}.
\]
}
\ethm
\medskip

Besides these results above, we need to stress a special and interesting case near the corner,  when the contact angle $\om=\f\pi{2n}$ for some $n\in\N$. In fact, {\it the solution to the related elliptic problem for D-N operator has no singular part in this case, which means that we can have the same estimates for D-N operator without singular parts as in classical cases}. The price here is that one would require additional compatible conditions on boundary conditions of the elliptic system, which will  be discussed in detail in Section 5 as well.  

In the end of this paper, we also consider about  regularizing transformations in order to reach the `natural' and proper regularity for the upper surface $z=\eta(x)$ through all the estimates we have: As long as one has $|\eta|_{H^{K+2}}$ in the coefficients, it can be always replaced by $|\eta|_{H^{K+\f32}}$.  

After dealing with the estimates for D-N operator, there will be our subsequent work on the well-posedness of the related 2 dimensional water-wave problem. In this paper,  we only consider some fixed contact angle $\om$ without the time variable, so the changing contact angle and the continuous dependence for the singularities with respect to the time will be discussed in the subsequent paper. 

\noindent{\bf Organization of this paper}. To start with, Section 3 discusses the range of the contact angle $\om$, which is important in the following sections. Section 4 gives the Trace Theorem for  corner domains. Section 5 provides us the decompositions for the variational solution of the mixed problem, as well as the higher-order estimates for both the regular part and singular coefficients. Section 6 shows the properties of D-N operator in the presence of singularities. Finally, Section 7 discusses about the regularizing transformations.

\section{Notations}
\noindent -  We denote by $\cS$ the corner domain $\cS=\{(x,z)|\,0\le z \le kx,\,x\ge 0\}$ for some constant $k>0$.\\
\noindent -  We denote by $\cS_0$ the corner domain $\cS_0=\{(x,z)|-\ga x\le z\le \eta'(0) x,\,x\ge 0\}$ straighten from $\Om$. Sometimes we take $\cS_0$ as a bounded triangle $\{(x,z)|-\ga x\le z\le \eta'(0) x, \,0\le x\le \del\}$ with some small constant $\del$,  since we focus only on a small neighborhood of the corner when we use $\cS_0$. The new right segment of the boundary is denoted by $\Ga_\del$.\\
\noindent -  We denote by $\om_1$ the angle of the upper boundary such that $\eta'(0)=\tan \om_1$, and $\om_2$  the angle of the bottom such that $\ga=\tan\om_2$.\\
\noindent -  We denote by $\om$ the contact angle for the corner domain $\Om$ and $\cS_0$. We have $\om=\om_1+\om_2$.\\ 
\noindent - When no confusion will be made, we  denote uniformly by $\Ga_t$ ($\Ga_b$) the upper (lower) boundary for the domain $\Om$, $\cS$ or $\cS_0$. \\ 
\noindent -  The unit tangential vector for the upper (lower) boundary of $\Om$, $\cS$ or $\cS_0$ is denoted by ${\bf\tau}_t$ ($\bf\tau_b$) when no confusion will be made. Similarly we denote by ${\bf n}_t$ ($\bf n_b$) the unit outward normal vector for the upper (lower) boundary. The orientation for the pairs $(\bf n_t,\bf \tau_t)$ and $(\bf n_b,\bf \tau _b)$ is the same as $(x,y)$ coordinates.\\
\noindent - The notation $P^t$ stands for the transpose of a vector or a matrix $P$. \\
\noindent -  We define the cut-off function $\be(r)$ compactly supported  in $\cS_0$ with $r$ the polar radius. We also use the same notation $\be$ as the corresponding cut-off function in $\Om$ or $\cS$ when no confusion will be made.\\
\noindent -  The constant  $b$ from the boundary condition for $\Ga_b$ can be written as $b=\tan\Phi$ for some constant angle $\Phi$.\\
\noindent -  We denote by $V$ or $V(\cS_0)$ the space  on the triangle $\cS_0$ such that
\[
V(\cS_0)=\{u\in H^1(\cS_0)\big| \,u|_{\Ga_t}=u|_{\Ga_\del}=0\}.
\]
\noindent -  We denote by $V^2_b(\cS_0)$ the space on  the triangle $\cS_0$ such that 
\[
V^2_b(\cS_0)=\big\{u\in H^2(\cS_0)\big|\, u|_{\Ga_t}=u|_{\Ga_\del}=0,\,\p_{n_b} u+b\p_{\tau_b} u|_{\Ga_b}=0\big\}.
\] 
\noindent -  $\cN_b$ is defined as the orthogonal space of $\Del(V^2_b)$ in $L^2(\cS_0)$.\\
\noindent - We denote by $\dot H^{s}(\Ga)$ the homogenous space of $H^s(\Ga)$, where $\Ga$ is some part of the boundary for $\Om$.\\
\noindent - The space $\tilde H^{s}(\Ga)$ ($s>0$) is defined as
\[\tilde H^{s}(\Ga)=\Big\{u\in \dot{H}^{s}(\Ga)\Big| \f{D^\al u}{\rho^{\si}}\in L^2(\Ga),\,|\al|=m\Big\}\] where $\rho=\rho(X)$ is the distance (arc length) between the point $X\in \Ga$ and the end points of $\Ga$, and $s=m+\si$ with $\si \in [0,1)$.  The norm is defined as 
\[
|u|^2_{\tilde H^{s}}=|u|^2_{H^s}+\int_\Ga \f{|D^\al u|^2}{\rho^{2\si}}dX.
\] 
Moreover,  we use $(\tilde H^s)^*$ to denote the dual space of $\tilde H^s$. We will use in this paper the cases $s=1/2$ and $s=3/2$.\\
\noindent - We denote by $|\cdot|_{H^s}$  the Sobolev norm  for $H^s(\Ga_t)$ or $H^s(\Ga_b)$ or $H^s(\R^+)$ on the boundary with $s>0$, and by $|\cdot|_{W^{m,\infty}}$ the $W^{m,\infty}$ norm on the boundary with $m\in\N$. \\
\noindent - We denote by $\|\cdot\|_{m,2}$ the Sobolev norm for $H^m$ in the domain $\Om$, $\cS$ or $\cS_0$, and by $\|\cdot\|_{m,\infty}$  $W^{m,\infty}$ norm in the domain with $m\in\N$. When $m=0$, we use $\|\cdot\|_2$ for $L^2$ norm and $\|\cdot\|_\infty$ for $L^\infty$ norm.\\
\noindent - When we use $\p^m$ on $\Om$, $\cS$ or $\cS_0$ for some $m\in \N$, we mean $\p^{\al_1}_x\p^{\al_2}_z$ with some $\al_1+\al_2=m$. \\
\noindent - We denote by $\bar\eta(x)$ the function $\bar\eta(x)=\eta(x)+\ga x$ on $\R^+$, and we denote by $\eta'(x)$ the derivative of $\eta$ with respect to $x$.\\
\noindent - The function $d(z)$ is defined as 
\[
d(z)=\f 1k-\f1{\ga+\eta'\big(\bar\eta^{-1}(z)\big)},
\]
where $\bar \eta^{-1}$ is the inverse of $\bar\eta$ if exists,  and $d_\infty$ means $d(z)|_{k=\infty}$. Moreover, we set $d_0=d(0)$.\\
\noindent - We define the eigenvalue $\lam=-\f{\pi}{2\om}$.

\section{The contact angle}
We focus first  on the contact point $O(0,0)$ to evaluate the range of the contact angle $\om$. One will see in the following that the range of the contact angle $\om$ plays an important role in the singular analysis.

First of all, from the boundary  condition on the bottom for the  water-wave problem  we have  \[V\cdot {\bf n}_b|_{\Ga_b}=0\] where $V$ is the velocity of the fluid in the domain. This condition still holds at the contact point $O$ as long as  $V$ exists point-wisely (the regularity in our paper will do). Taking the material derivative $D_t=\p_t+V\cdot\na_{x,z}$ on this condition (since the fluid particles flow along $\Ga_b$) and applying Euler equation on $\Ga_b$, one finds  near the contact point $O$ that
\[0=D_t V\cdot{\bf n}_b|_{\Ga_b}=(-\na_{x,z} p-g{\bf e}_z)\cdot{\bf n}_b|_{\Ga_b}\] with $p$  the pressure, $g$ the gravity coefficient, ${\bf e}_z=(0,1)^t$ while noticing that ${\bf n}_b=-(1+\ga^2)^{-\f12}(\ga,1)^t$ near $O$.  Naturally, this gives the condition for  the pressure near $O$:
\beq\label{contact condition 1}
-\na_{x,z} p\cdot {\bf n}_b|_{\Ga_b}=g{\bf e}_z\cdot {\bf n}_b|_{\Ga_b}<0.
\eeq

On the other hand, when the surface tension at the free surface  $\Ga_t$ is neglected,  one has $p|_{\Ga_t}=const$ and moreover one assumes the Taylor sign condition $-\na_{n_t} p|_{\Ga_t}>0$ on $\Ga_t$, which together  imply that $\na_{x,z} p$ acts as  an inward normal vector  on $\Ga_t$. Combining this fact with \eqref{contact condition 1} at the contact point (the regularity in our paper  also  makes \eqref{contact condition 1} hold point-wisely),  we can see that the angle $\om$ between the free surface $\Ga_t$ and the bottom $\Ga_b$ satisfies
$\om\in [0,\,\pi/2)$. In this paper, we avoid the special case $\om=0$ where the boundary  has a cusp at the contact point. Consequently, we consider for some  contact angle \[\om\in (0,\f\pi2)\] through our paper.

We need to mention that, in the case of vertical walls, Alazard, Burq and Zuily proved in Prop 6.1\cite{ABZ15} that the free surface makes a right-angle with the rigid walls, as long as the Taylor sign condition is satisfied.

\section{Trace theorem for mixed boundary conditions on curvilinear corner domains}
 Generally, the trace theorem will have some extra compatible conditions when the domain boundary is piecewise smooth, which is different from the smooth domain case \cite{PG1}.  Firstly, we would like to quote  Thm 1.5.2.1 and Thm 1.5.2.4\cite{PG1} here, where the first theorem deals with the traces on part of the boundary, while the second theorem is a special case focusing on the traces for the first quadrant and the only corner point is $(0,0)$. These theorems will be used later in this paper.
\bthm{Theorem}\label{trace on Ga j}{\it Let $\Om_p$ be a bounded open subset of $\R^2$, whose boundary is a curvilinear polygon of class $C^{K,1}$; then for $j$th curve $\Ga_j$(counterclockwise) from the boundary, the mapping 
\[u\mapsto \big\{u,\,\p_{n_j}u,\dots,\,\p^l_{n_j}u\big\}|_{\Ga_j},\quad l<m-\f12,\] which is defined for $u\in \cD(\bar\Om_p)$ and $\bf n_j$ is the unit outward normal vector on $\Ga_j$, has a unique continuous extension as an operator from 
\[H^m(\Om_p)\quad \hbox{onto}\quad \prod_j H^{m-l-\f12}(\Ga_j),\qquad l\le m-1\le K.\]}\ethm
\bthm{Remark}\label{trace on Ga j remark}{\it Thm \ref{trace on Ga j}  tells us that we can deal with traces on one piece of the boundary  similarly as in the classical trace theorem.  In fact, the proof for this theorem is almost the same as the classical trace theorem, so we can have the estimates similarly
\[|\p^l_{n_j}u|_{H^{m-l-\f12}(\Ga_j)}\le C \|u\|_{m,2}\] for some constant $C$ depends on $\Om_p$.  }
\ethm
 
\bthm{Theorem}\label{trace thm PG}{\it The mapping $u\mapsto\big\{\{f_k\}^{m-1}_{k=0},\,\{g_l\}^{m-1}_{l=0}\big\}$ defined by 
\[f_k=D^k_y u|_{y=0},\quad g_l=D^l_x u|_{x=0}\] for $u\in \cD(\overline{\R^+\times\R^+})$, has a unique continuous extension as an operator from $H^m(\R^+\times\R^+)$ onto the subspace of 
\[\prod^{m-1}_{k=0}H^{m-k-\f12}(\R^+)\times\prod^{m-1}_{l=0}H^{m-l-\f12}(\R^+)\]
defined by \\
\noindent(a) $D^l_xf_k(0)=D^k_yg_l(0)$, $l+k<m-1$ and,\\
\noindent(b) $\int^1_0\big|D^l_xf_k(t)-D^k_yg_l(t)\big|^2\f{dt}t <+\infty$, $l+k=m-1$.
}
\ethm
\bthm{Remark}\label{oblique case}{\it Similar result on a general polygon as Thm \ref{trace thm PG} also holds, see Thm 1.5.2.8\cite{PG1}. In this case, we use the tangential and normal derivatives $\p_{\tau},\p_{n}$ instead of $\p_x, \p_y$ above, and the parameter $t$ would be chosen as the arc length parameter w.r.t. the corner point.
}
\ethm
 One can see in this theorem that there are compatible conditions near the corner point, which is different from the classical trace theorem. In fact, the space $\tilde H^\f12$ used in this paper originates from the second condition, see Thm \ref{Thm 1.5.2.3}  in Appendix A.
 
 \subsection{The transformation to straighten the domain} Before stating our trace theorem, we would like to  introduce two transformations  on $\Om$ which  straighten the moving free boundary near or away from  the contact point $O(0,0)$. These transformations will be used in the proof of the trace theorem and in later sections.

First of all we introduce the local transformation near the point $O$ which maps the corner domain  $\cS\cap U_{\del S}$ with some constant $k>0$ into the moving domain $\Om\cap U_\del$:
\[T_S: \quad (x,\,z)\in S\cap U_{\del S}\mapsto (\bar x,\,\bar z)\in\Om\cap U_\del\]
with
\[\bar x=x+\bar\eta^{-1}(z)-\f 1k z,\quad \bar z=z-\ga\Big(x+\bar\eta^{-1}(z)-\f 1k z\Big)\]
where $\bar\eta^{-1}(z)$ is the inverse of $\bar\eta(x)=\eta(x)+\ga  x$, and the neighborhood $\cS\cap U_{\del S}$, $\Om\cap U_\del$ near the corner $(0,0)$ correspond to the  triangle $\cS_0$ by these transforms introduced in this section.  Note that the bottom $\Ga_b$ can be written as $z=-\ga x$ near the corner. We write for some constant $k>0$ here because we will let $k\rightarrow \infty$ sometimes. The corner domain $\cS$ is convenient for computations in our paper with a horizontal lower boundary.

 The idea for this transform $T_S$ is that  firstly one pushes the bottom of $\Om$ up to $z=0$ by replacing $z$ with $z+\ga x$, then one straightens  the new corner  domain by adjusting x-coordinate for each $z$.  One can check that $T_S$ maps the boundaries $\bar z=-\ga x$ and $\bar z=\eta(\bar x)$ of $\Om$ into the boundaries $z=0$ and $z=k x$ of $\cS$ respectively.   On the other hand, when $\bar \eta^{-1}$ exists, we know that it has the same regularity with $\eta$, so the transformation $T_S$ has the same regularity with $\eta$. Moreover, $T_S$ is invertible with the inverse transformation $T_S^{-1}$ as follows
\[T_S^{-1}:\quad (\bar x,\,\bar z)\in \Om\cap U_\del\mapsto (x,\,z)\in \cS\cap U_{\del S}\] where
\[x=\bar x-\bar\eta^{-1}(\ga\bar x+\bar z)+\f 1k (\ga \bar x+\bar z),\quad z=\ga\bar x+\bar z.\]

Secondly, we also need to construct another transform straightening the rest part of $\Om$. In fact, we give the transform $T_R$ which maps a flat strip $R$ to the rest part of $\Om$ here:
\[T_R:\quad (x,\,z)\in R\mapsto (\bar x,\,\bar z)\in \Om\]
with
\[\bar x= x,\quad \bar z=\eta(x)z+l(x)(1-z),\] where $R=\{(x,\,z)|\, x\ge x_\del,\,0\le z\le 1\}$ is a flat strip, and $x_\del>0$ is a constant fixed by  $U_\del$.  The inverse transform $T_R$ is
\[T_R^{-1}:\quad (\bar x,\,\bar z)\in \Om\cap \{\bar x\ge x_\del\}\mapsto (x,\,z)\in R\] where
\[x=\bar x,\quad z=\f{\bar z-l(\bar x)}{\eta(\bar x)-l(\bar x)}\]

A few words about $\del$:  In the transformations above, we need $\eta(x)$ is invertible on $\Om\cap U_\del$ with respect to $x$, which would be true since we will choose $\del>0$ small enough.

Sometimes we also need to use the linear transformation 
\[
T_0:\quad X=(x,\,z)\in \cS_0\mapsto \tilde X=(\tilde x,\,\tilde z)=P^t_0 X\in\cS
\] 
with 
\[
P_0=\left(\begin{matrix}1+\ga d_0 & \ga\\ d_0 & 1\end{matrix}\right),
\]
 which maps the corner domain $\cS_0$ to the corner domain $\cS$. Notice that $\cS_0$ is the straighten version of the original domain $\Om$ and therefore it keeps the same contact angle $\om$. 

To close this section, we consider a special case of the local transform $T_S$  when the flat corner domain $S$ is the first quadrant:
\[T_S: \quad (x,\,z)\in S=\{(x,\,z)|\,x,\,z\ge0\}\cap U_{\del S}\mapsto (\bar x,\,\bar z)\in\Om\cap U_\del\]
with
\[\bar x=x+\bar\eta^{-1}(z),\quad \bar z=z-\ga\Big(x-\bar\eta^{-1}(z)\Big).\] The transform is again invertible with the inverse transform
\[T_S^{-1}:\quad (\bar x,\,\bar z)\in \Om\cap U_\del\mapsto (x,\,z)\in S\] where
\[x=\bar x-\bar\eta^{-1}(\ga\bar x+\bar z),\quad z=\ga\bar x+\bar z.\]
\bthm{Remark}{\it  Indeed one can use the conformal mapping to transform $\Om$ near the corner point to a rectilinear corner domain, meanwhile the conformal mapping can also keep the Laplace operator $\Del$ in the related elliptic equation. However one cannot express this mapping explicitly and estimates for the conformal mapping are needed. In \cite{WuK, Wu2,Wu3}, the authors uses conformal mapping to change the Lipschitz domain into a smoother domain, therefore the singularities are shifted to the mapping.

We choose in this paper an easier and explicit transformation $T_S$ to straighten the corner. The price is that we get a more general elliptic equation and the singular functions behave worse sometimes.  We need to improve  the regularities about singular functions in the following sections.}  \ethm

\subsection{Trace theorem on the corner domain $\Om$}
Since  we consider the mixed problem for the corner domain $\Om$,  a trace theorem with mixed boundary conditions is  presented here. These are  modified versions from Thm 2.1\cite{BR1}. On the other hand, we will also need to consider the pure Dirichlet boundary conditions for some special cases, so we also present a trace theorem for Dirichlet case. 

\bthm{Theorem}\label{trace}{\it (Mixed Boundary) Let $m\ge 2$, $b=\tan\Phi$ and  functions $f\in H^{m-\f12}(\Ga_t)$, $g\in H^{m-\f 32}(\Ga_b)$ be given. 
When the two vectors ${\bf n}_b+b{\bf \tau}_b\parallel{\bf\tau}_t$, i.e. $\Phi-\om=l\pi+\f\pi2$ for some $l\in\Z$, one assumes additionally that  $F\in H^{m-\f32}(\Ga_t\cup\Ga_b)$ where
\[
F=\left\{\begin{array}{ll}
\na_{\tau_t}f,\quad \qquad x\in \Ga_t,\\
\f1{|{\bf n}_b+b{\bf \tau}_b|}g ,\quad x\in \Ga_b.
\end{array}\right.
\] Then there exists a function $v\in H^{m}(\Om)$ satisfying the following mixed boundary conditions
\beno
v|_{\Ga_t}=f,\quad \p_{n_b} v+b\p_{\tau_b}v|_{\Ga_b}=g.
\eeno

Moreover,  we have the following estimate
\[
\|v\|_{m,2}\le    C(\delta,\Ga_b,|\eta|_{H^{m-\f12}})(|f|_{H^{m-\f12}}+|g|_{H^{ m-\f32}}).
\]
 }
\ethm
\proof
{\bf Case 1: ${\bf n}_b+b{\bf \tau}_b\nparallel{\bf\tau}_t$.} 
The key point in this proof is about the trace near the corner.  Let 
\[f_\del=\left\{\begin{array}{ll}
f,\quad x\le \del_0\\
0,\quad x >2\del_0\end{array}\right.
\quad g_\del=\left\{\begin{array}{ll}
g,\quad x\le \del_0\\
0,\quad x >2\del_0.\end{array}\right.\] with some $\del_0$ depending on $\del$ and small. We can have  $f_\del\in H^{m-\f12}(\Ga_t)$ and $g_\del\in H^{m-\f32}(\Ga_b)$ cut-off from $f,g$ and supported in $\Om\cap U_\del$.  We now use the transform $T_S$ with $k=\infty$ which can change locally the domain $\Om\cap U_{\del}$ near the corner point $(0,0)$ into part of the first quadrant $S$.  We denote the new boundaries by $\widetilde{\Gamma}_t=\{x=0,\,0\le z\le \del_1\}$, $\widetilde{\Gamma}_b=\{0\le x\le \del_1,\,z=0\}$ with the same corner point $O(0,0)$ and $\del_1$ a small constant depending on $\del$. 

As a result, we change the boundary conditions for $v$ into conditions for $\tilde v=v\circ T_S$:
\beq\label{bdry condi for tilde v}
\tilde v|_{\widetilde{\Gamma}_t}=\tilde f_\del\in H^{m-\f12}(\tilde \Ga_t),\quad  a_1 \p_x\tilde v+a_2 \p_z\tilde v |_{\widetilde{\Gamma}_b}=\tilde g_\del\in H^{m-\f32}(\tilde \Ga_b),
\eeq
where $\tilde f_\del=f_\del\circ T_S$, $\tilde g_\del=(1+\ga^2)^\f12 g_\del\circ T_S$ and
\beno
a_1(z)=b-\ga-(1+\ga^2)d_\infty(z), \quad\hbox{and}\quad a_2=-(1+\ga^2),
\eeno
with $d_\infty(z)=-\f1{\ga+\eta'\big(\bar\eta^{-1}(z)\big)}$. Here $\ga+\eta'\big(\bar\eta^{-1}(z)\big)\neq 0$ since the tangent line for the  upper boundary won't be parallel to the tangent line for the lower boundary in the small neighborhood of $\del$.

Now we want to  solve $\tilde v$ from its boundary conditions when $|a_1(z)|>0$ on $[0,\del_1]$ for some $\del_1>0$ (which will be investigated in the end of the proof). Using zero  extension, we can find $\tilde f_\del\in H^{m-\f12}(\R^+)$ and $\tilde g_\del\in H^{m-\f32}(\R^+)$,  so the boundary conditions are given on $\R^+$ now. As a result, we have that
\beno
|\tilde f_\del|_{H^{m-\f12}(\R^+)}\leq C |f|_{H^{m-\f12}(\widetilde{\Gamma}_t)},\qquad |\tilde g_\del |_{H^{m-\f32}(\R^+)}\leq C|f|_{H^{m-\f32}(\widetilde{\Gamma}_t)}.
\eeno

Let's  consider the case $m=3$,  since the cases $m=2$ were proved already in Theorem 2.1 \cite{BR1}.  And for the case that $m\ge 4$, we can get the result in a  similar way.  

By using Theorem \ref{trace thm PG}, we have that for any given functions $f_0,\,g_0\in H^{\f52}(\R^+),\,f_1,\,g_1\in H^{\f32}(\R^+),\,f_2,\,g_2\in H^{\f12}(\R^+)$, if the following conditions are satisfied
\beq\label{condi from PG trace thm}
\begin{split}&g_0(0)=f_0(0),\quad g_1(0)=\p_zf_0(0), \quad\p_xg_0(0)=f_1(0),\quad\hbox{and}\\
&f_2-\p^2_z g_0,\,\p^2_x f_0-g_2,\,\p_xf_1-\p_zg_1\in \tilde H^{\f12}(\R^+)
\end{split}\eeq
then there exists $\tilde v\in H^{3}(S)$ such that
\beno
&&\tilde v|_{\widetilde{\Gamma}_t}=f_0, \p_x\tilde v|_{\widetilde{\Gamma}_t}=f_1, \p^2_x\tilde v|_{\widetilde{\Gamma}_t}=f_2,\quad\hbox{and}\\
&&\tilde v|_{\widetilde{\Gamma}_b}=g_0, \p_z\tilde v|_{\widetilde{\Gamma}_b}=g_1, \p^2_z\tilde v|_{\widetilde{\Gamma}_b}=g_2.
\eeno
From now on we use the same variable $s\in \R^+$ for $f_i(s),\,g_j(s)$ ($i,j=0,1,2$).
To get what we want, we need  $\tilde v$ to satisfy the boundary conditions \eqref{bdry condi for tilde v}, that is
\beno
f_0= \tilde f_\del,\quad a_1 g_0'+a_2g_1=\tilde g_\del.
\eeno Thus, we set first
\beno
g_1=\tilde f'_\del+h
\eeno
for some $h$ to be fixed and
\ben\label{equ:a_1}
g_0'=a_1^{-1}\big(\tilde g_\del-a_2\tilde f'_\del\big), \quad \textrm{on}\quad [0,\delta_1]
\een
with the initial data $g_0(0)=\tilde f(0)$.

Solving the equation above, we can find  $g_0\in H^{\f{5}{2}}([0, \del_1])$, which
can be extended to $g_0 \in  H^{\f{5}{2}}(\R^+)$. Moreover, we can derive that
\beno
|g_0 |_{ H^{\f{5}{2}}(\R^+)}\leq C(\delta,\Ga_b, |\eta|_{H^{\f52}})(|f|_{H^{\f{5}{2}}}+|g|_{H^{\f{3}{2}}}).
\eeno
Now we define $h\triangleq -a^{-1}_2(a_1 g'_0- \tilde g_\del +a_2 \tilde f'_\del)$ on $\R^+$.  As a result, we find $g_1\in H^\f32(\R^+)$ satisfying $g_1(0)=f'_0(0)$ from \eqref{condi from PG trace thm}, and we also have the estimate:
\beno
|g_1|_{ H^{\f{3}{2}}(\R^+)}\leq  C(\delta,\Ga_b,|\eta|_{H^{\f52}})(|f|_{H^{\f{5}{2}}}+|g|_{H^{\f{3}{2}}}).
\eeno

For $f_1$, we firstly define it by solving the differential equation
\[
f_1'= g_1'\quad\hbox{on}\quad [0,\,\del_1]
\]
with the initial data $f_1(0)= g_0'(0)$, and then extend it from $H^{\f32}([ 0,\del_1])$ into $H^{\f32}(\R^+)$.  Thus, we find  $f_1\in H^{\f32}(\R^+)$ satisfying \eqref{condi from PG trace thm}.

For $g_2$, we define it by
\beno
g_2=f''_0
\eeno
which implies that  $g_2\in H^{\f12}(\R^+)$.

For $f_2$, we define it by
\beno
  f_2=g''_0
\eeno
which implies that $f_2\in H^{\f12}(\R^+)$.  One can check easily that all the conditions in \eqref{condi from PG trace thm} are satisfied here. Meanwhile,  thanks to the proof of Thm \ref{trace thm PG} (Thm 1.5.2.4 in \cite{PG1}), we have the estimate
\[
\|\tilde v\|_{3,2}\le C\sum\limits_{i=0,1,2}\big(|f_i|_{H^{\f 52-i}(\R^+)}+|g_i|_{H^{\f 52-i}(\R^+)}\big)\] with some constant $C$ depending on $m=3$, which leads to 
\[
\|\tilde v\|_{3,2}\le  C(\delta,\Ga_b,|\eta|_{H^{\f52}})(|f|_{H^{\f{5}{2}}}+|g|_{H^{\f{3}{2}}}).
\] As a result, we can go back to the function $v_\del=\tilde v\circ T^{-1}_S\in H^3(\Om)$ which satisfies the boundary conditions
\[v_\del|_{\Ga_t}=f_\del,\quad \p_{n_b} v_\del+b\p_{\tau_b} v_\del|_{\Ga_b}=g_\del\] and the desired estimate for $v_\del$ holds. 

For the case $m\ge 4$, the proof also starts with fixing $g_0$, $g_1$ as above, and the remaining part follows similarly. 

On the other hand, considering the trace away from $(0,0)$, we set 
\[f^c_\del=f-f_\del,
\quad g^c_\del=g-g_\del\] We can have  $f^c_\del\in H^{m-\f12}(\Ga_t)$, $g^c_\del\in H^{m-\f32}(\Ga_b)$ and $f^c_\del,\,g^c_\del$ vanish near the corner. So we can use the transform $T_R$ to change the domain $\Om\cap U^c_\del$  into the strip $R$, and therefore the trace problem can be solved using a classical trace theorem as in \cite{Lannes}: One finds $v^c_\del\in H^3(\Om)$ satisfying the boundary conditions
\[v^c_\del|_{\Ga_t}=f^c_\del,\quad \p_n v^c_\del+b\p_\tau v^c_\del|_{\Ga_b}=g^c_\del\] and the desired estimate for $v^c_\del$ also holds.  The proof is omitted here. 

Placing $v_\del$ and $v^c_\del$ together, we finally arrive at the function \[v=v_\del+v^c_\del,\] which satisfies the boundary conditions and also the estimate in this theorem.  

To close the proof,  we investigate the condition $|a_1(z)|> 0$ on $[0,\del_x]$. In fact, one only needs $a_1(0)\neq 0$ and $\del_x$ small enough. The condition $a_1(0)\neq 0$ can be rewritten as 
\[a_1(0)=\tan\Phi+\f1{\tan\om}\neq 0,\] which means $\Phi-\om\neq l\pi+\f\pi2$, or  equivalently that 
the vector ${\bf n}_b+b{\bf \tau}_b$ on $\Ga_b$ is  not parallel to the tangent vector ${\bf\tau}_t$ on $\Ga_t$.  

\medskip

{\bf Case 2: ${\bf n}_b+b{\bf \tau}_b\parallel{\bf\tau}_t$.}
In the proof of case 1, we know that the equation \eqref{equ:a_1} can be not solved this time, which should be adjusted here a little bit. 

In fact, plugging $g_1=\tilde f'_\del+h$ into the boundary condition $a_1g'_0+a_2g_1=\tilde g_\del$, we find 
\[a_1g'_0+a_2 h=\tilde g_\del-a_2\tilde f'_\del.\] On the other hand,  some computations lead to 
\[|{\bf n}_b+b{\bf\tau}_b|=(1+\ga^2)^\f12\f{\sqrt{1+\eta'(0)^2}}{\ga+\eta'(0)}\ \ \hbox{on}\ \Om,\quad \hbox{while} \quad \f{\sqrt{1+(\eta')^2}}{\ga+\eta'}\na_{\tau_t}f_\del=-\big( \tilde f'_\del\big)\circ T^{-1}_S\]
so  our assumption on $F$ implies  
\beno
\tilde{g}_\delta-a_2\tilde{f}_\delta' \in H^{\f32}(\R^+),
\eeno with $\big(\tilde{g}_\delta-a_2\tilde{f}_\delta'\big)(0)=0$,  which  results in that $a_1g'_0+a_2 h\in  H^{\f32}(\R^+)$.
Moreover, we define that
\beno
g'_0=G, \quad g_{0}(0)=\tilde f_\delta(0).
\eeno
where $G$ is an arbitrary function belongs to $ H^{\f32}(\R^+)$ satisfying $G(0)=0$. Then we have $g_0\in H^\f52(\R^+)$ and $ h\in H^\f32(\R^+)$ satisfying $h(0)=0$. As a result, $g_1=\tilde f'_\del+h\in H^\f32(\R^+)$ and satisfies \eqref{condi from PG trace thm}.

For the other terms $f_1,f_2, g_2$, since they do not depend on $a_1$, we use the same definitions as in case 1. The case when $m\ge 4$ follows similarly as long as $g_0$, $g_1$ are fixed. In this way the proof is finished.

\ef

As the end of this section, we also  present a trace theorem with Dirichlet boundary conditions.
\bthm{Theorem}\label{Dirichlet boundary trace}(Dirichlet boundary){\it 
The trace operator 
\[v\mapsto v|_{\Ga_t\cup\Ga_b}\] 
is continuous from $H^m(\Om)$ to $H^{m-\f12}(\Ga_t\cup\Ga_b)$. On the other hand, this operator has a right continuous inverse: Let functions $f\in H^{m-\f12}(\Ga_t)$, $f_d\in H^{m-\f12}(\Ga_b)$ be given.  \\
\noindent (1)When $m\ge 2$, one needs 
\[f(O)=f_d(O)\] with $O$ the contact point.  \\
\noindent (2)When $m=1$, set
\[
F=\left\{\begin{array}{ll}f,\quad \hbox{on}\quad \Ga_t,\\
f_d\quad\hbox{on}\quad \Ga_b.\end{array}\right.
\]
 satisfying $F\in H^{\f12}(\Ga_t\cup\Ga_b)$.  
  
Then there exists a function $v\in H^m(\Om)$ satisfying Dirichlet boundary conditions
\[v|_{\Ga_t}=f,\quad v|_{\Ga_b}=f_d.\]
Moreover, one has the estimate
\[
\|v\|_{m,2}\le   C(|f|_{H^{m-\f12}}+|f_d|_{H^{m-\f12}})
\] with constant $C$ depending on $\del,\Ga_b$, and on $|\eta|_{H^{m-\f12}}$ when $m\ge2$, while on $|\eta|_{H^2}$ when $m=1$.
}\ethm
\begin{proof} 
Following the proof for Thm \ref{trace} line by line, we need to focus near the corner and  prove the case $m=3$. The case $m=1$ is proved  by Gagliardo(1957) for a bounded open set with Lipschitz boundary in Thm 1.5.1.3\cite{PG1}, and the case $m=2$ is proved already in  \cite{BR1}. 

 In fact, we have here 
\[
\td v|_{\td \Ga_t}=\td f_\del,\quad \td v|_{\td\Ga_b}=\td f_d,
\] where $\td f_d$ is defined similar as $\td g_\del$. This means we know $f_0,\,g_0\in H^{\f52}$ in  \eqref{condi from PG trace thm} and need to find $f_1,\,g_1\in H^{\f32}$, $f_2,\,g_2\in H^{\f12}$. Since these functions are easy to find out, so we have $\td v\in H^3(\Om)$ satisfying boundary conditions as   long as $f_0(0)=g_0(0)$ holds. Therefore the proof is done.
\end{proof}
\bthm{Remark}\label{one side Dirichlet trace}{\it
When only one Dirichlet boundary condition is given, for example, $v|_{\Ga_t}=f$ with $f\in H^{m-\f12}(\Ga_t)$, we can use Thm \ref{Dirichlet boundary trace} (define Dirichlet boundary condition on $\Ga_b$ according to $f$, for example) or some classical extension trick to find such a  $v\in H^m(\Om)$ satisfying the boundary condition without any compatible considerations of course.
Moreover, the estimate holds
\[
\|v\|_{m,2}\le C(|\eta|_{H^{m-\f12}})|f|_{H^{m-\f12}}
\] when $m\ge 2$.
Besides this, the trace operator
\[v\mapsto v|_{\Ga_t}\] 
is also continuous from $H^m(\Om)$ to $H^{m-\f12}(\Ga_t)$.
}\ethm

\section{The mixed boundary elliptic problem}
We consider in general the elliptic
problem for $u$ when proper conditions $h,\,f,\,g$ are given:
\[
\mbox{(MBVP)}\quad\left\{\begin{array}{ll}
\Delta u=h,\qquad \hbox{in}\quad \Omega\\
u\,|_{\Ga_t}=f,\qquad \p_{n_b}u\,|_{\Ga_b}=g
\end{array}\right.
\]
and we will find the elliptic estimates for this system.
Referring to the higher-order estimates, sometimes we will also need to consider
\[
\mbox{(MBVP)}_b\quad\left\{\begin{array}{ll}
\Delta u=h,\qquad \hbox{in}\quad \Omega\\
u\,|_{\Ga_t}=f,\qquad \p_{n_b}u+b\,\p_{\tau_b}u\,|_{\Ga_b}=g
\end{array}\right.
\] 
with some constant $b=\tan\Phi$. Moreover, we will reduce $\mbox{(MBVP)}$ and $\mbox{(MBVP)}_h$ to  the homogeneous boundary case most of the time:
\[
\mbox{(MBVP)}_h\quad\left\{\begin{array}{ll}
\Delta u=h,\qquad \hbox{in}\quad \Omega\\
u\,|_{\Ga_t}=0,\qquad \p_{n_b}u+b\,\p_{\tau_b}u\,|_{\Ga_b}=0.
\end{array}\right.
\]  
In this section, we firstly derive the elliptic estimates for the system $\mbox{(MBVP)}_h$ with $b=0$. We also consider about the regularity for $\mbox{(MBVP)}_h$ when $b\neq0$. After the estimates for the homogeneous system, we go back to retrieve the estimates for the system  $\mbox{(MBVP)}$ by trace theorems.

\subsection{$H^2$ decomposition and the elliptic estimate}
In this section, we firstly derive the $H^2$ estimate for $\mbox{(MBVP)}_h$ with $b=0$, then we go back to the $H^2$ estimate for $\mbox{(MBVP)}$. To begin with,  we present an existence theorem about the variational solution for $\mbox{(MBVP)}_h$, and a $H^2$ decomposition theorem for $\mbox{(MBVP)}_b$ when $b\neq0$,  which will be needed in higher-order regularities. 

The existence and uniqueness of the variational solution  for the mixed boundary system is proved already in Lemma 4.4.3.1\cite{PG1}, where $\Om_p$ is a bounded polygonal domain and the elliptic operator is $\Del$. In fact, one can extend this result naturally into our case with the corner domain $\Om$. We omit the proof here.
\bthm{Lemma}\label{variation solution}{\it  Let $h\in L^2(\Om)$ be a given function, then there exists a unique variational solution $u\in H^1(\Om)$ for the homogeneous system $\mbox{(MBVP)}_h$. Moreover, one has the variational estimate
\[\|\na u\|_2 \le C(|\eta|_{W^{1,\infty}})\|h\|_2.\]
}\ethm 

Based on the variational solution $u\in H^1(\Om)$, we  now consider about  the $H^2$ and higher-order regularities. In fact, the $H^2$ regularity has been proved in \cite{PG1,CD93II, D92,BR1}etc. with different assumptions. In  \cite{PG1},  the $H^2$ regularity was  discussed and the higher-order case was also briefly mentioned for the rectilinear polygon case. The author established the $H^2$ estimate there but the higher-order estimates were missing. The authors proved an estimate for general boundary condition with the boundary piecewise analytic in \cite{CD93II} .  In \cite{BR1}, the authors considered the general $C^{k,\sigma}$ curvilinear polygon case and a $H^2$ decomposition was derived based on Grisvard's result \cite{PG1}. 

In this paper, based on the result in \cite{BR1}, we first derive the $H^2$ estimate  on $\Om$, and then we take a further step to explore  the higher-order estimates  for $\mbox{(MBVP)}_h$ while tracing the dependence of coefficients on the upper free boundary.  

The major difference between our corner-domain case and the classical smooth-domain case is that there will be singularity appearing from the elliptic system. In order to have the elliptic estimates, the singularities must be explicitly studied. In fact, Theorem 3.2.5\cite{BR1} tells us that if $\Om_p$ is a  curvilinear polygon with a piecewise $\mathcal C^{2,\si}$($\si\ge0$) boundary, then the  variational solution $u\in H^1(\Om_p)$  can be uniquely decomposed into the sum of a $H^2(\Om_p)$ function $u_r$ and a singular part in $H^1\setminus H^2(\Om_p)$. On the other hand, we know from \cite{PG1} that the singularity only concentrates near the corner point,  so this theorem also holds for our unbounded corner domain $\Om$. Therefore a slightly modified version of Thm 3.2.5\cite{BR1} for the general inhomogeneous system $\mbox{(MBVP)}_b$ is presented here, adjusted for the purpose of our paper.

\bthm{Theorem}\label{H2 decompose}{\it Suppose that  the corner domain $\Om$ has a $\mathcal C^{2,0}$ upper boundary $\eta$. Let the constant $b=\tan\Phi$, and functions $h\in L^2(\Om), f\in H^{\f32}(\Ga_t)$, $g\in H^{\f12}(\Ga_b)$ be given. When $n_b+b\tau_b\parallel \tau_t$ i.e. $\Phi-\om=l\pi+\f\pi 2$ for some $l\in\Z$,  one requires additionally that 
\beq\label{compatible condition}
F\in H^{\f12}(\Ga_t\cup\Ga_b)\quad\hbox{with}\quad
F=\left\{\begin{array}{ll}
\na_{\tau_t}f, \qquad\quad\hbox{on}\quad \Ga_t,\\
\f1{|{\bf n}_b+b{\bf \tau}_b|}g ,\quad\hbox{on}\quad \Ga_b.
\end{array}\right.
\eeq
Then the variational solution $u\in H^1(\Om)$ for the elliptic problem $\mbox{(MBVP)}_b$ can be decomposed as
\[
u=u_r+\sum_{m\in A}c_{m}s_{m},
\] 
with $u_r\in H^2(\Om)$, $c_{m}$ real coefficients and singular functions
\[
s_{m}(r,\tht)=\beta(r)r^{-\lam_{m}}\cos{\big(\lam_{m}(\tht+\om_2)+\Phi\big)},
\]
where $(r,\,\tht)$ is the polar coordinates on $\cS_0$  straightened from $\Om$ near  the corner point using the transform  $(P^{-1}_0)^tT_S$ and $\tht\in [-\om_2,\,\om_1]$. Moreover, one has
\[
\lam_m=\f{-\Phi+(m+1/2)\pi}{\om}\qquad \hbox{for}\quad m\in \Z,
\]  
and $A=\{m\,|\,-1<\lam_m<0\}$.}
 \ethm
\begin{proof} In fact, since the singularity only concentrates near $(0,0)$, we find the decomposition for $\be u$ from Thm 3.2.5\cite{BR1} with $\be$ the cut-off function near $(0,0)$. To conclude the proof, we only need to consider about the regularity of $(1-\be)u$ away from the corner, and we can find indeed that $(1-\be)u\in H^2(\Om)$ by using the transformation $T_R$ and applying standard elliptic analysis on the strip $R$ (see for example \cite{Lannes}).
\end{proof}
\bthm{Remark}\label{no singular function} {\it  In fact, for the problem $\mbox{(MBVP)}$ with Neumann boundary condition (the constant $b=0$),  recalling that the contact angle $\om\in(0,\,\pi/2)$, one has \[A=\{m| -1<\lam_{m}<0\}=\emptyset.\] Consequently, combining Thm \ref{trace} and Thm \ref{variation solution}, there exists a unique variational solution $u$ and moreover one has $u\in H^2(\Om)$ without any singular part, so we can expect the classical $H^2$ estimate for this problem.}
\ethm
\bthm{Remark}\label{equiv comp condition}{\it In the following sections, instead of $\Om$, we will deal with the compatible condition \eqref{compatible condition} on the corner domain $\cS$ most of the time, so an equivalent form on $\cS$ is needed. In fact, a direct computation (similar as in Case 2 of the proof for Thm \ref{trace}) leads to the equivalent compatible condition 
\beq\label{compatible condition S}
F\in H^{\f12}(\Ga_t\cup\Ga_b)\quad\hbox{on $\cS$, with}\quad
F=\left\{\begin{array}{ll}
\na_{\tau_t}f, \qquad\quad\hbox{on}\quad \Ga_t ,\\
\f 1{|\cP_S{\bf n}_b+b{\bf\tau}_b|}g ,\quad\hbox{on}\quad \Ga_b.
\end{array}\right.
\eeq with $f, g$  the boundary conditions for $v$ on $\cS$:
\[
v|_{\Ga_t}=f,\quad \p^{\cP_S}_{n_b}v+b\p_{\tau_b}v|_{\Ga_b}=g.
\] Moreover, when we consider the corner domain $\cS_0$, we have similar compatible condition, replacing the matrix $\cP_S$ by $\bar\cP_S$.
}\ethm

When the contact angle $\om=\f{\pi}{2n}$ for some  $n\in \N$, we need to consider the following system with Dirichlet boundary conditions:
\[
\mbox{(DBVP)}\quad\left\{\begin{array}{ll}
\Delta u=h,\qquad \hbox{in}\quad \Omega\\
u\,|_{\Ga_t}=f,\qquad  u\,|_{\Ga_b}=f_d.
\end{array}\right.
\] The corresponding $H^2$ decomposition theorem is presented here, again slightly modified from Thm 3.2.5\cite{BR1}. There is a similar result in Remark 2.4.6 \cite{PG2}.
\bthm{Theorem}\label{H2 decompose 1}{\it Suppose that  the corner domain $\Om$ has a $\mathcal C^{2,0}$ upper boundary $\eta$ and let $h\in L^2(\Om)$, $f\in H^{\f32}(\Ga_t)$ and $f_d\in H^{\f32}(\Ga_b)$ be given. If $f(0)=f_d(0)$ is satisfied, then there exists a  solution $u\in H^2(\Om)$ for the elliptic problem $\mbox{(DBVP)}$.}
 \ethm

Next, we give the $H^2$ norm estimate for $u$. Based on Remark \ref{no singular function}, we know that the variational solution $u\in H^2(\Om)$ for the system $\mbox{(MBVP)}_h$ with $b=0$. In order to finish the $H^2$ estimate, we need to rewrite this system both near the corner and away from the corner using the transformations $T_S,\,T_R$ respectively. In fact, the solution $u$ can be decomposed into two elliptic problems for
\[u_c(\bx,\,\bz)=\be u(\bx,\,\bz),\quad u_{ac}(\bx,\,\bz)=(1-\be)u(\bx,\,\bz)\] where we denote by $(\bx,\,\bz)$ any point in $\Om$.  The elliptic problem for  $u_c$ is
\beq\label{u1 problem}
\left\{\begin{array}{ll}\Del u_c= h_{co}\qquad \hbox{on}\quad \Om\\
u_c|_{\Ga_t}=0,\quad \p_{n_b}u_c|_{\Ga_b}= g_{co}
\end{array}\right.
\eeq  
where
\[
h_{co}=\be h-[\be,\,\Del]u,\quad g_{co}=-(\p_{n_b}\be)u|_{\Ga_b}.
\] 
Similarly, the elliptic problem for $u_{ac}$ is
\beq\label{u2 problem}
\left\{\begin{array}{ll}\Del u_{ac}= h_{ac}\qquad \hbox{on}\quad \Om\\
u_{ac}|_{\Ga_t}=0,\quad \p_{n_b}u_{ac}|_{\Ga_b}= g_{ac}\end{array}\right.\eeq where
\[h_{ac}=(1-\be) h-[(1-\be),\,\Del]u,\quad g_{ac}=-\big(\p_{n_b}(1-\be)\big)u|_{\Ga_b}.\]

We will start firstly to find the $H^2$ estimate for $u_c\in H^2(\Om)$, which is different from classical elliptic analysis on domains with smooth boundaries. After applying $T_S$ on $u_c$ and denoting
\[v_c(x,z)=u_c\circ T_S(x,z),\quad\hbox{for}\quad \forall (x,z)\in \cS\] one finds the equivalent elliptic problem for $v_c$ (compactly supported in $\cS\cap U_{\del S}$) as
\beq\label{v1 problem}
\left\{\begin{array}{ll}\na\cdot \cP_S\na v_c=h_c\qquad \hbox{on}\quad \cS\\
v_c|_{\Ga_t}= 0,\quad \p^{\cP_S}_{n_b}v_c|_{\Ga_b}= (1+\ga^2)^{1/2}g_c\end{array}\right.\eeq where
\[
h_c=h_{co}\circ T_S,\quad g_c=g_{co}\circ T_S,
\]
and the coefficient matrix 
\[
\cP_S=P^t_SP_S\quad\hbox{with}\quad
P_S=(\na T^{-1}_S)\circ T_S=\left(\begin{matrix} 1+\ga \,d(z)& \ga\\ d(z) & 1\end{matrix}\right)
\] 
and 
\[
d(z)=\f 1k-\f 1{\eta'(\bar\eta^{-1}(z))+\ga}.
\]
\bthm{Lemma}\label{P_S} {\it$P_S$ is invertible with
\[P_S^{-1}=\left(\begin{matrix}1 & -\ga\\-d(z) &1+\ga\,d(z)\end{matrix}\right).\]  Moreover, the coefficient matrix $\cP_S=P^t_SP_S$ in \eqref{v1 problem} satisfies the coercive condition
\[\xi\cdot \cP_S\xi\ge C_0(\|P_S\|_{\infty})|\xi|^2\qquad\hbox{for}\quad \forall \xi\in\R^2\] for some constant $C_0$ depending on $\|P_S\|_{\infty}$.}\ethm
Instead of using the classical technique of difference quotients, the $H^2$ estimate for $v$ can be done here  by  a standard perturbation argument known as Korn's  procedure based on the $H^2$ estimate Thm 4.3.1.4\cite{PG1}: We froze the coefficient matrix $\cP_s$ at $(0,0)$ first, and perform the inverse transform at this corner point to reformulate the elliptic problem \eqref{v1 problem} on a triangle region  as a perturbation problem with respect to the classical elliptic problem with the operator $\Del$. Consequently, applying Theorem 4.3.1.4\cite{PG1} directly, we can find the desired $H^2$ estimate.

More generally, we present the elliptic estimate for the following  system.
\bthm{Proposition}\label{H2 estimate v1} {\it Let functions $h\in L^2(\cS)$, $f\in H^\f32(\Ga_t)$ and $g\in H^\f12(\Ga_b)$ be given.  Suppose that $v\in H^2(\cS)$ be a function compactly supported near the corner in $\cS\cap U_{\del S}$ and a solution of the elliptic problem
\[\left\{\begin{array}{ll}\na\cdot \cP_S\na v=h\qquad \hbox{on}\quad \cS\\
v|_{\Ga_t}= f,\quad \p^{\cP_S}_{n_b}v+b\p_{\tau_b}v|_{\Ga_b}= g
\end{array}\right.\] 
with some constant $b=\tan\Phi$. Moreover, when $\Phi-\om=l\pi+\f\pi2$ for some $l\in\Z$, one requires additionally the compatible condition \eqref{compatible condition S}. Then one can have the $H^2$ elliptic estimate 
\[\|v\|_{2,2}\le C(|\eta|_{H^2\cap W^{2,\infty}})\big(\|h\|_2+|f|_{H^{\f32}}+|g|_{H^{\f12}}+\|v\|_2\big).\]}
\ethm

As a result, applying this proposition to the system \eqref{v1 problem} for $v$ with $b=0$ and going back to the system \eqref{u1 problem} we can derive the corresponding $H^2$ estimate for $u_c=\be u$ as follows
\beq\label{u1 H2 estimate}\|u_c\|_{2,2}\le C(|\eta|_{H^2\cap W^{2,\infty}})\big(\|h\|_2+\|u\|_{1,2}\big),\eeq where the condition $n_b+b\tau_b\nparallel \tau_t$ is satisfied naturally since $b=0$ here and the contact angle $\om\in (0,\f\pi2)$. 

\begin{proof} (proof for Prop\ref{H2 estimate v1}) 
Step 1: Reduce the elliptic equation to homogenuous boundary problem. In order to apply Thm \ref{trace}, consider $v\circ T_S^{-1}$ on $\Om$ to arrive at the equivalent boundary conditions 
\[
v\circ T_S^{-1}|_{\Ga_t}=f\circ T_S^{-1},\quad (\p_{n_b}+b\p_{\tau_b})(v\circ T_S^{-1})|_{\Ga_b}=(1+\ga^2)^{-\f12}(g\circ T_S^{-1}).
\]  
Here, we notice that $f, g$ are compactly supported. After applying Thm \ref{trace} (note that the compatible condition is satisfied here by Remark \ref{equiv comp condition}),   one can derive that there exists a function $u_0\in H^2(\Om)$ compactly supported near the corner  in $\Om\cap U_{\del}$ and satisfying
\[u_0|_{\Ga_t}=f\circ T_S^{-1},\quad \p_{n_b}u_0+b\p_{\tau_b}u_0|_{\Ga_b}=(1+\ga^2)^{-\f12}(g\circ T_S^{-1})\] and
\[\| u_0\|_{2,2}\le C(|\eta|_{H^\f32})\big(|f\circ T_S^{-1}|_{H^{\f32}}+|g\circ T_S^{-1}|_{H^{\f12}}\big).\] As a consequence, we also find $ v_0=u_0\circ T_S\in H^2(S)$ compactly supported in $\cS\cap U_{\del S}$ and  satisfying the boundary conditions 
\[
v_0|_{\Ga_t}=f,\quad \p^{\cP_S}_{n_b}v_0+b\p_{\tau_b}v_0|_{\Ga_b}=g
\] 
and moreover
\[
\|v_0\|_{2,2}\le C(|\eta|_{H^\f32})\big(|f|_{H^{\f32}}+|g|_{H^{\f12}}\big).
\] 
Then, we define $\tilde v=v-v_0$ which satisfies 
\beq\label{homo problem for tv 1}
\left\{\begin{array}{ll}\na\cdot \cP_S\na \tilde v=\tilde h,\qquad \hbox{on} \quad \cS\\
\tilde v|_{\Ga_t}=0,\,\qquad \p^{\cP_S}_{n_b}\tilde v+b\p_{\tau_b} \tilde v|_{\Ga_b}=0, \end{array}\right.\eeq
where $\tilde h=h-\na\cdot\cP_S\na v_0$.

The $H^2$ norm of $v$ is done by 
\[
\|v\|_{2,2}\le \|\tilde v\|_{2,2}+\|v_0\|_{2,2}\le C\big(\|\tilde v\|_{2,2}+|f|_{H^{\f32}}+|g|_{H^{\f12}} \big).
\]
As a result, it remains to find the estimate of $\|\tilde v\|_{2,2}$.

Step 2: the homogenuous boundary problem. We simply take $k=1$ in the following computations for simplicity. We first consider the corresponding homogenuous boundary problem for a function $\tilde v\in H^2(S)$ as
\beq\label{homo problem for tv}
\left\{\begin{array}{ll}\na\cdot \cP_S\na \tilde v=\tilde h,\qquad \hbox{on} \quad \cS\\
\tilde v|_{\Ga_t}=0,\,\qquad \p^{\cP_S}_{n_b}\tilde v+b\p_{\tau_b} \tilde v|_{\Ga_b}=0 \end{array}\right.\eeq
where $\tilde v$ is also compactly supported in $\cS\cap U_{\del S}$.  Froze the coefficient matrix $\cP_S=P^T_SP_S$ at the corner point $(0,0)$ to find the corresponding matrix
\[\cP_0=P^t_0P_0,\quad P_0=P_S(0,0)=\left(\begin{matrix}1+\ga d_0 & \ga\\ d_0 & 1\end{matrix}\right)\] with $d_0=d(z)|_{(0,0)}=1-\f1{\ga+\eta'(0)}$.

After introducing a new function  $w(X)=\tilde v(P^t_0 X)$ with $X=(x,z)^t$ and recalling the linear transform
\[T_0:\quad X=(x,\,z)\in \cS_0\mapsto \tilde X=(\tilde x,\,\tilde z)=P^t_0 X\in\cS,\] we know that the elliptic problem \eqref{homo problem for tv} for $\tilde v$ on $\cS$ is changed into the following elliptic problem for $w$ on $\cS_0$ with the same contact angle $\om$ as for $\Om$:
\beq\label{inhomo problem for w}
\left\{\begin{array}{ll}\Del w=\tilde h\circ T_0+\big(\na\cdot(\cP_0-\cP_S)\na \tilde v\big)\circ T_0\qquad \hbox{on} \quad \cS_0\\
w|_{\Ga_t}=0,\,\qquad \p_{n_b}w+b\p_{\tau_b} w|_{\Ga_b}=g_w\end{array}\right.
\eeq 
with 
\[
g_w:=(1+\ga^2)^{-1/2}\big(\p^{\cP_0-\cP_S}_{n_b}\tilde v|_{\tilde z=0}\big)\circ T_0.
\]
For the moment,  we can also take $\cS_0$ as a bounded triangle domain as introduced in the notation part, since $w$ is again compactly supported.  We have $w|_{\Ga_\del}=0$.

Step 3: Estimate for $w$ on $\cS_0$. We want to apply the $H^2$ estimate from Theorem 4.3.1.4 \cite{PG1} on the system \eqref{inhomo problem for w} for $w$ on the triangle  $S_0$. Since this theorem requires homogeneous boundary conditions again, we first need to homogenize the boundary conditions in \eqref{inhomo problem for w}. Applying Thm \ref{trace}  one can always find a function $w_t\in H^2(\cS_0)$ compactly supported in  $\cS_0$ and  satisfying
\[w_t|_{\Ga_t}=0,\quad \p_{n_b}w_t+b\p_{\tau_b}w_t|_{\Ga_b}=g_w\] where $w_1$ is controlled by the boundary condition:
\[
\|w_t\|_{2,2}\le C|g_w|_{H^{\f12}}.
\]  Notice that in Thm \ref{trace}, when $\Phi-\om=l\pi+\f\pi 2$ we need the compatible condition for the boundary conditions, i.e. we need $g_w\in \tilde H^\f12(\Ga_b)$ here by Thm \ref{trace thm PG} and Remark \ref{oblique case}. This condition is indeed satisfied, since we know from $\tilde v$ system that $\p_{n_b}\tilde v\in \tilde H^\f12(\Ga_b)$ because $\tilde v\in H^2(\cS)$ and $\tilde v|_{\Ga_t}=0$ (where we use Thm \ref{trace thm PG} and Remark \ref{oblique case} again), and  $\p_{\tau_b}\tilde v\in \tilde H^\f12(\Ga_b)$ from the Neumann boundary condition.
On the other hand, we can write $w_t|_{\Ga_\del}=0$ since $w_1$ is compactly supported near the corner.

Defining $w_h=w-w_t$, we arrive at the homogeneous problem for $w_h$ as follows
\[\left\{\begin{array}{ll}\Del w_h=-\Del w_t+\tilde  h\circ T_0+\big(\na\cdot(\cP_0-\cP_S)\na \tilde v\big)\circ T_0\qquad \hbox{on} \quad \cS_0\\
w_h|_{\Ga_t}=0,\,\qquad \p_{n_b} w_h+b\p_{\tau_b}w_h|_{\Ga_b}=0,\quad
w_h|_{\Ga_\del}=0\end{array}\right.\] Which is a standard elliptic problem concerning Laplacian $\Del$ on the triangle $S_0$. Consequently we can apply directly Thm 4.3.1.4\cite{PG1} to find
\[\|w_h\|_{2,2}\le C\big(\|\Del w_h\|_2+\|w_h\|_2\big)\] with some constant $C$ depending only on $\cS_0$, which infers that
\beno
\|w_h\|_{2,2}&\le & C\big( \|\Del w_t\|_2+\|\tilde h\circ T_0\|_2+\|(\na\cdot(\cP_0-\cP_S)\na u)\circ T_0\|_2+\|w_h\|_2\big)\\
&\le& C\big(\|\tilde h\|_2+\|(\cP_0-\cP_S)\na u\|_{1,2}+\|w_h\|_2\big)\\
&\le& C\big(\|\tilde h\|_2+\|\na \cP_S\|_\infty\|u\|_{1,2}+\|\cP_0-\cP_S\|_\infty\|u\|_{2,2}+\|w\|_2\big).
\eeno Going back to system \eqref{homo problem for tv} for $\tilde v$ we have
\beno\|\tilde v\|_{2,2}&\le& C\|w\|_{2,2}\\
&\le & C\big(\|\tilde h\|_2+\|\na \cP_S\|_\infty\|\tilde v\|_{1,2}+\|\cP_0-\cP_S\|_\infty\|\tilde v\|_{2,2}+\|\tilde v\|_2\big)\eeno

Since $\cS_0$ is a triangle near the corner point $(0,0)$ with a small diameter $\del>0$, one finds
\[\|\cP_0-\cP_S\|_\infty\le C_\del\|\na \cP_S\|_\infty\] where $C_\del$ is also a small constant. Plugging this inequality into the above estimates for $\tilde v$ results in
\beq\label{H2 estimate for tv}\|\tilde v\|_{2,2}\le C\big(\|\tilde h\|_2+\|\tilde v\|_2\big)\eeq where $C=C(|\eta|_{H^2\cap W^{2,\infty}})$ is a constant.

Combining the estimates $\tilde v$ and $v_0$, we get the lemma proved.

\end{proof}

In order to finish the $H^2$ estimate for the solution $u$ to $\mbox{(MBVP)}_h$ with $b=0$, we need  the analysis away from the corner  on $u_{ac}$,  which is a standard procedure as in \cite{Lannes}.  We apply $T_R$ on $u_{ac}$ and set
\[v_{ac}(x,z)=u_{ac}\circ T_R(x,z)\] to rewrite \eqref{u2 problem} into the  problem for $v_{ac}$ as
\beq\label{v2 problem}
\left\{\begin{array}{ll}\na\cdot \cP_R\na v_{ac}=h_{ac}\qquad \hbox{on}\quad R\\
v_{ac}|_{z=0}= 0,\quad \p^{\cP_R}_{n}v_{ac}|_{z=-1}= (1+\ga^2)^{-1/2}g_{ac},\quad
v_{ac}|_{x=x_\del}=0\end{array}\right.\eeq where
\[\cP_R=P^t_RP_R,\quad h_{ac}=\bar h_{ac}\circ T_R,\quad g_{ac}=\bar g_{ac}\circ T_R,\]
and the coefficient matrix
\[P_R=(\na T^{-1}_R)\circ T_R=\left(\begin{matrix} 1& -\f{(\eta'-l')z+l'}{\eta-l}\\ 0 & \f1{\eta-l}\end{matrix}\right).\]
\bthm{Lemma}\label{P_R} {\it $P_R$ is invertible with
\[P_R^{-1}=\left(\begin{matrix}1 & (\eta'-l')z+l'\\0&\eta-l\end{matrix}\right).\]  Moreover, the coefficient matrix $\cP_R=P^t_RP_R$ in \eqref{v2 problem} satisfies the coercive condition
\[\xi\cdot \cP_R\xi\ge C^{-1}_0(\|P_R\|_{\infty})|\xi|^2\qquad\hbox{for}\quad \forall \xi\in\R^2.\]}\ethm The $H^2$ estimate on the flat strip $R$ can be done as in \cite{Lannes} and we omit the details here. We arrive at $H^2$ estimate for $u_{ac}=(1-\be)u$
\beq\label{u2 H2 estimate}\|u_{ac}\|_{2,2}\le C(\|P_R\|_{1,\infty})\big(\|h\|_2+\|u\|_{1,2}\big).\eeq 
Summing \eqref{u1 H2 estimate} and \eqref{u2 H2 estimate} up, one has the $H^2$ estimate for the solution $u$ of $\mbox{(MBVP)}_h$ with $b=0$ as follows
\beq\label{u H2 homo estimate}\|u\|_{2,2}\le C\big(\|h\|_2+\|u\|_{1,2}\big)\le  C\|h\|_2,\eeq 
where the constant $C=C(|\eta|_{H^2\cap W^{2,\infty}})$, and the  $L^2$ estimate for $u$ is derived from the Poincar\'e inequality and the estimate from Lemma \ref{variation solution} 
\[
\|u\|_2\le C\|\na u\|_2\le C(\|P_S\|_\infty)\|h\|_2,
\] since $u|_{\Ga_t}=0$ in the system $\mbox{(MBVP)}_h$ with $b=0$.

Consequently, going back to the corresponding inhomogeneous system $\mbox{(MBVP)}$ and combing \eqref{u H2 homo estimate}, Thm \ref{trace} and Remark \ref{no singular function}, the  solution $u\in H^2(\Om)$ for $\mbox{(MBVP)}$ has the following  $H^2$ estimate :
\beq\label{u H2 estimate}
\|u\|_{2,2}\le C(|\eta|_{H^2\cap W^{2,\infty}})\big(\|h\|_2+|f|_{H^\f32}+|g|_{H^\f12}\big).
\eeq

\subsection{$H^3$ decomposition and regularity} We will meet with singular functions when we consider higher-order regularities, which makes the higher-order estimates   different from classical estimates. In order to have a full understanding of the singular functions, we will focus on $H^3$ decomposition first, where we only have one singular function to start with.  After this part,  we can deal with higher-order cases with more singular functions. 

Similarly as in the $H^2$ estimate, we consider about  system $\mbox{(MBVP)}_h$ with $b=0$ firstly, and then we establish the corresponding result for system $\mbox{(MBVP)}$.
In this section, if $f$ is defined on $\cS$, we always denote by $\bar f=f\circ T_0$ the corresponding function on $\cS_0$. 

\subsubsection{Tangential derivative near the corner}
In order to have higher-order estimates, instead of system \eqref{v1 problem} for $v$, we consider here a more general system for $v$
\beq\label{v1 problem general}
\left\{\begin{array}{ll}\na\cdot\cP_S\na v=h\qquad\hbox{on}\quad \cS\\
v|_{\Ga_t }=f,\qquad \p^{\cP_S}_{n_b}v+b\p_{\tau_b} v|_{\Ga_b }= g\end{array}\right.
\eeq 
with $h\in H^1(\Om)$ , $f\in H^\f52(\Ga_t)$ and $g\in H^\f32(\Ga_b)$ given. Recall  that the constant $b=\tan \Phi$  and
\[
\cP_S=P_S^tP_S=\left(\begin{matrix}(1+\ga d)^2+d^2 & \ga+(1+\ga^2)d\\
\ga+(1+\ga^2)d & 1+\ga^2\end{matrix}\right)
\]
with $d=d(z)=1/k-1/(\ga+\eta'(\bar\eta^{-1}(z)))$ and $d_0=d(0)$ defined before.

Checking from Thm \ref{trace},  we already have $v\in H^2(\cS)$, while  noticing in the case when $\Phi-\om=l\pi+\f\pi2$ for some $l\in\Z$,   we requires additionally the compatible condition \eqref{compatible condition}.

Concerning $H^3$ estimates for $v$, we will take some derivative on $v$ to derive a new system for the derivative.  In fact,  recalling in Section 2, we express the unit tangent vector of the upper boundary $z=kx$ of $\cS$ as
\[{\bf\tau}_t=(\tau_1,\tau_{2})^t=-(1,\,k)^t/\sqrt{1+k^2}, \] and we define the derivative of $v$ along this tangential direction 
\[w_1=\na_{\tau_t} v=\tau_1\p_xv+\tau_2\p_z v,\] which is defined on $\cS$ indeed.

\bthm{Lemma}\label{tangent system}{\it
When $\Phi-\om\neq l\pi+\f12\pi$ ($l\in\Z$), one finds $w_1\in H^1(\cS)$ satisfying the following elliptic system
\[
\left\{\begin{array}{ll}\na\cdot\cP_S\na w_1=l_1\qquad\hbox{on}\quad S\\
w_1|_{\Ga_t}=\na_{\tau_t} f,\qquad \p^{\cP_S}_{n_b}w_1+\widetilde{b}\,\p_{\tau_b} w_1|_{\Ga_b}=R_1\end{array}\right.
\]
with
\beno
&&\widetilde{b} = \tan{(\Phi-\om)},\quad l_1=-\na\cdot(\na_{\tau_t} \cP_S)\na v+\na_{\tau_t} h,\\
&&R_1=\rho_1\p_xv|_{\Ga_b}-\tau_2h|_{\Ga_b}+\rho_2\p_x g.
\eeno
 Moreover, one has $\rho_1=\tau_2\big(\ga+(1+\ga^2)d'(0)\big)$ and 
\[
\rho_2=\tau_1-\f{\big[ (\widetilde{b}-b)  \tau_1+\big((1+\ga d_0)^2+d^2_0 \big) \tau_2 \big]}{\ga+(1+\ga^2)d_0-b}.
\]  }
\ethm

\begin{proof}
After applying $\na_{\tau_t}$ on the elliptic equation for $v$  in \eqref{v1 problem general}, it's straightforward to find the elliptic  equation for $w_1$.  Meanwhile, the boundary condition $w_1|_{\Ga_t}=\na_{\tau_t} f$ is derived in the same way. We need to focus on the boundary condition for $w_1$ on the bottom $\Ga_b$ i.e. when $z=0$.

In fact, the condition for $v$ on $z=0$ can be rewritten as
\[
g=\p^{\cP_S}_{n_b}v+b\p_{\tau_b}v|_{z=0}=-\big(\ga+(1+\ga^2)d_0-b\big)\p_xv|_{z=0}-(1+\ga^2)\p_zv|_{z=0},
\]
which means
\[
(\al-b)\p_xv|_{z=0}+\beta \p_zv|_{z=0}=-g,\qquad\hbox{with}\quad \al=\ga+(1+\ga^2)d_0,\quad \beta=1+\ga^2.
\]
We want to look for a condition on $z=0$ for $w_1$ as
\[
\p^{\cP_S}_{n_b} w_1+\widetilde{b} \p_{\tau_b} w_1|_{z=0}= \hbox{terms of} \, \p v,\,v
\]
for some $\widetilde{b}$. In order to do this, computing the left-hand side
\beno
\hbox{left}&=&\p^{\cP_S}_{n_b} w_1+\widetilde{b} \p_{\tau_b} w_1|_{z=0}\\
&=& - \al \p_xw_1|_{z=0}-\beta \p_zw_1|_{z=0}+\widetilde{b}\p_x w_1|_{z=0},
\eeno
and substituting $w_1=\tau_1\p_xv+\tau_2\p_zv$ into the expressions above to arrive at
\[
\hbox{left}= (\widetilde{b} -\al)\tau_1 \p^2_xv|_{z=0}+\big((\widetilde{b}-\al )\tau_2 -\beta\tau_1\big)\p_x\p_zv|_{z=0}-\beta\tau_2\p^2_zv|_{z=0}.
\]
On the other hand, we know from the elliptic equation in  \eqref{v1 problem general} that
\[
\p^2_zv|_{z=0}=-\f1{\beta}\big((1+\ga d_0)^2+d^2_0\big)\p^2_xv|_{z=0}-\f{2\al}{\beta}\p_x\p_zv|_{z=0}+\bar R_1
\]
with  $\bar R_1=-\f1{\beta}\big(\ga+(1+\ga^2)d'(0)\big)\p_xv|_{z=0}+\f 1\beta h|_{z=0}$, which results in that
\ben\label{neumann condition for w1}
\hbox{left}&= &\big[(\widetilde{b} -\al)\tau_1+\big((1+\ga d_0)^2+d^2_0 \big) \tau_2 \big]\p^2_xv|_{z=0}\nonumber\\
&&\quad+\big( \widetilde{b} \tau_2+\al\tau_2-\beta\tau_1\big)\p_x\p_zv|_{z=0}-\beta\tau_2\bar R_1\nonumber\\
&=&\big[ (\widetilde{b}-b)  \tau_1+\big((1+\ga d_0)^2+d^2_0 \big) \tau_2 \big]\p^2_xv|_{z=0} \nonumber\\
&&\quad+\big( \widetilde{b} \tau_2+\al\tau_2 \big)\p_x\p_zv|_{z=0}-\beta\tau_2\bar R_1+\tau_1\p_x g.
\een
Since $\bar R_1$ contains only lower order terms of $v$ (order 1 and order 0), we  choose  $\widetilde{b}$ to satisfy
\[
\big[ (\widetilde{b}-b) \tau_1+\big((1+\ga d_0)^2+d^2_0\big) \tau_2 \big]\beta=\big(  \widetilde{b} \tau_2+\al\tau_2 \big) (\al-b)
\]
to make the 2nd-order terms of $v$ in \eqref{neumann condition for w1}  turn into term of $g$ finally, which results in  
\[
\widetilde{b}=\f{1+b\big((1+\gamma^2)d_\infty(0)+\gamma\big)}{ \ga+(1+\ga^2)d_\infty(0)-b}
\] 
with $d_\infty(0)=d(0)|_{k=\infty}=-\f1{\ga+\eta'(0)}$.

To finish the proof, we need that $\ga+(1+\ga^2)d_\infty(0)-b\neq 0$, which means equivalently $\Phi-\om\neq l\pi+\f12\pi$ by a direct computation since $b=\tan\Phi$ and 
\[
\ga+(1+\ga^2)d_\infty(0)=\f{\ga\eta'(0)-1}{\ga+\eta'(0)}=-\f1{\tan\om}.
\] 
\end{proof}
Sometimes we also need to consider the case when $\Phi-\om= l\pi+\f12\pi$. The following lemma tells us that in this case the boundary condition on $\Ga_b$ for $w_1$ system becomes Dirichlet condition.
\bthm{Lemma}\label{tangent system special}{\it
When $\Phi-\om= l\pi+\f12\pi$ for some $l\in\Z$ and  if the compatible condition \eqref{compatible condition S} is satisfied, then one has  $w_1\in H^1(\cS)$ satisfying the following elliptic system
\[
\left\{\begin{array}{ll}\na\cdot\cP_S\na w_1=l_1\qquad\hbox{on}\quad S\\
w_1|_{\Ga_t}=\na_{\tau_t} f,\qquad w_1|_{\Ga_b}=\f{k}{(1+\ga^2)\sqrt{1+k^2}}g\end{array}\right.
\]}
\ethm
\begin{proof}
The proof lies in the boundary condition on $\Ga_b$ i.e. $z=0$ again. In fact,  from \eqref{v1 problem general} one has for $v$ on  $z=0$  that
\[g=-\big(\ga+(1+\ga^2)d_0-b\big)\p_xv|_{z=0}-(1+\ga^2)\p_zv|_{z=0}.\] Since we have $\Phi=\om+m\pi+\f12\pi$, we can rewrite
\[b=\tan\Phi=-\cot\om=-\f{1-\ga\eta'(0)}{\ga+\eta'(0)},\]
as a result we find
\beno
\ga+(1+\ga^2)d_0-b&=& \ga-(1+\ga^2)\Big(\f1k-\f1{\ga+\eta'(0)}\Big)-\tan\Phi\\
&=& \ga+\f1k(1+\ga^2)-\f{1+\ga^2}{\ga+\eta'(0)}+\f{1-\ga\eta'(0)}{\ga+\eta'(0)}\\
&=&\f1k(1+\ga^2).
\eeno The computations above lead to 
\[g=-\f1k(1+\ga^2)\p_xv|_{z=0}-(1+\ga^2)\p_zv|_{z=0}=\f1k\sqrt{1+k^2}(1+\ga^2)w_1|_{z=0}\] since $w_1=\na_{\tau_t}v=\f1{\sqrt{1+k^2}}(\p_xv+k\p_zv)$.

To close the proof, when the condition \eqref{compatible condition} is satisfied, we know that $v\in H^2(\cS)$ and so $w_1\in H^1(\cS)$. As for classical traces, the trace for $w_1$ on the boundary $\Ga_t\cup\Ga_b$ should be in $H^\f12$ according to Thm \ref{Dirichlet boundary trace}, which leads to the compatible condition from Thm \ref{Dirichlet boundary trace} for $\na_{\tau_t}f$, $\f{k}{(1+\ga^2)\sqrt{1+k^2}}g$. In fact, the condition  is exactly the same one as \eqref{compatible condition S}.
\end{proof}

\medskip

Now we focus on the equation of $w_1$ equation. We are ready to get back to $v$ system \eqref{v1 problem} and consider the $H^3$ decomposition for $\mbox{(MBVP)}_h$ with $b=0$ near the corner. Setting $v=v$ with  $b=0$ and applying Thm \ref{H2 decompose} on the system for $w_1=\na_{\tau_t}v$  in Lem \ref{tangent system} (so one has $\tilde b=-\tan\om$ which means $\Phi=-\om$ here),  one finds that there exist a real coefficient  $\bar c_s$ and a singular function $s(r,\tht)$ such that $w_1\in H^1(\cS)$ has the following $H^2$ decomposition
\beq\label{tangent decompose}w_1=w_r+\chi_1(\om)\bar c_s\,s\eeq where $w_r\in H^2(S)$ is the regular part of $w_1$,  and the singular function 
\[
s(r,\theta)=\beta(r)\,r^{-\lam-1}\cos{\big((\lam+1)(\theta+\om_2)-\om\big)}\in H^1\setminus H^2(\cS_0)
\] 
 is defined on the corner domain $\cS_0=\{(x,z)\,|-\ga x\le z\le \eta'(0) x\}$. Recall that $\cS=\cS_0\circ T_0$ from Prop \ref{H2 estimate v1} with $T_0$ a linear transform. Here one has 
 \[
 \lam=-\f\pi{2\om}\in (-2,-1),\,\hbox{when}\, \om\in(\f\pi4,\f\pi2)\quad\hbox{and}\, \chi_1(\om)=\left\{\begin{array}{ll}1,\quad\om\in(\f\pi 4,\f\pi 2)\\
 0,\quad\hbox{otherwise}.\end{array}\right.
 \] 
Moreover, recall that $\be(r)$ is a cut-off function.

In fact, the singular function is obtained  by taking $m=-1$ in Theorem \ref{H2 decompose}   when  the contact  angle $\om\in (\pi/4,\,\pi/2)$. In case when $\om\in (0,\f\pi 4]$, there is no singular part. Similarly, recall that  there is no singular part in $H^2$ case neither.

Since $w_1=\na_{\tau_t} v$ is the derivative of $v$,  we want also to write the singular part $s(r,\theta)$ in a similar way to retrieve the singular part for $v$ itself.  Indeed, after some computations we find a new function $\bar S(r,\tht)$ defined on on the corner domain  $\cS_0$
\[\bar S(r,\theta)=\beta(r) r^{-\lam}\cos{\big(\lam(\theta+\om_2)\big)}\]
or equivalently $S$ on domain $\cS$
\[
S=\bar S\circ T^{-1}_0.
\] satisfying
\[\na_{\tau_t}  S=-\lam\big|(P_0^{-1})^t{\tau_t}\big|\,s\circ T^{-1}_0+[\na_{\tau_t},\be(r)\circ T^{-1}_0]r^{-\lam}\cos\big(\lam(\tht+\om_2)\big)\circ T^{-1}_0\quad\hbox{on}\quad \cS, \] where the constant vector $\tau_t$ is extended all over $\cS$ naturally indeed.
Moreover, one can tell directly that  $\bar S\in H^2\setminus H^3(\cS_0)$ and satisfies
\[
\left\{\begin{array}{ll}\Del \bar S=[\Del,\be(r)]r^{-\lam}\cos{\big(\lam(\theta+\om_2)\big)},\qquad\hbox{on}\quad \cS_0\\
\bar S|_{\Ga_t}=0,\qquad \p_{n_b}\bar S|_{\Ga_b}=0\end{array}\right.
\] where one has $\p_{n_b}=-\f1r\p_\tht$ on $\Ga_b$.
 
As a result, we can rewrite \eqref{tangent decompose} as
\[w_r=\na_{\tau_t}\big(v- \chi_1(\om)c_s\,S\big)+c_s[\na_{\tau_t},\be(r)\circ T^{-1}_0]r^{-\lam}\cos\big(\lam(\tht+\om_2)\big)\circ T^{-1}_0\] 
where $c_s$ is a constant coefficient\[c_s=-\f1{\lam\big|(P_0^{-1})^t{\tau_t}\big|} \bar c_s.\]

\subsubsection{$H^3$ decomposition} Denoting
\[v_r=v_c-\chi_1(\om) c_s\,S\] and recalling $w_r=\na_{\tau_t}v_r\in H^2(\cS)$ from the previous subsection,  now we are in a position to finish the $H^3$ decomposition for $v_c$.  Classical analysis would be to find the normal derivative w.r.t $\na_{\tau_t}$ from the elliptic equation for $v_r$, while noticing that $v_c$ and $S$ satisfy different elliptic equations:
\[\na\cdot\cP_S\na v_c=h_1,\quad \hbox{while}\quad \na\cdot\cP_0\na S=([\Del,\be(r)]r^{-\lam}\cos{\big(\lam(\theta+\om_2)\big)})\circ T^{-1}_0\qquad \hbox{on} \quad \cS.\]  In order to find the elliptic equation for $v_r$, we need first to show that $\na\cdot\cP_S\na S$ is also in $H^1(\cS)$.
\bthm{Lemma}\label{S elliptic P_S system}{\it
For the given function  $\bar S(r,\theta)=\beta(r) r^{-\lam}\cos{\big(\lam(\theta+\om_2)\big)}$ on $\cS_0$, one has that \[\na\cdot\bar\cP_S\na \bar S\in H^1(\cS_0),\] where $\bar\cP_S=(P_0^{-1})^t\cP_SP_0^{-1}$ is the corresponding coefficient  matrix of $\cP_S$ on the domain $\cS_0$, and $P_0=\left(\begin{matrix}1+\ga d_0 & \ga\\ d_0 & 1 \end{matrix}\right)$ comes from $\cP_0=\cP_S|_{z=0}=P_0^tP_0$.}\ethm

\begin{proof} Recalling the linear transformation $T_0$ from $\cS_0$ to $\cS$ and the notation $d_0=d(0)$, we can find the expression for $\bar \cP_S$ on the corner domain $\cS_0$ based on $\cP_S$ as follows
\[
\bar\cP_S=\left(\begin{matrix}1+2\ga(d-d_0)+(1+\ga^2)(d-d_0)^2 & (1-\ga^2)(d-d_0)-\ga(1+\ga^2)(d-d_0)^2\\ (1-\ga^2)(d-d_0)+\ga(1+\ga^2)(d-d_0)^2 & 1-2\ga (d-d_0)+\ga^2(1+\ga^2)(d-d_0)^2\end{matrix}\right).
\]
 On the other hand, we can also rewrite $d_0$ as $d(r,\theta)|_{r=0}$ since we use the polar coordinates here. A direct computation  leads to that
\beno
 d(r,\theta)-d(0,\theta)&=&r\int^1_0\p_rd(s\,r,\,\theta)ds\\
&=&-r\int^1_0 \f{(\ga\cos \theta+\sin\theta)\p_x^2\eta\big(\bar\eta^{-1}(s\,r(\ga\cos \theta+\sin\theta))\big)}{\Big(\ga+\eta'\big(\bar\eta^{-1}(s\,r(\ga\cos \theta+\sin\theta))\big)\Big)^3}ds.
\eeno 
Therefore one can write  
\[
d(r,\theta)-d(0,\theta)=r\,l(r,\,\theta),
\]
 where $l(r,\theta)$ is a bounded function for $(r,\,\theta)$.

The computations above tells us that we can find an extra $r$ from $\bar\cP_S-Id$, so we rewrite
\[\na\cdot\bar\cP_S\na \bar S=\Delta\bar S+\na\cdot(\bar\cP_s-Id)\na\bar S\] where we know already $\Del\bar S\in H^\infty(\cS_0)$.  Checking term-by-term, we find that $\na\cdot\bar\cP_S\na \bar S$ is of order $r^{-\lam-1}$ with $\om\in (\f\pi 4,\f\pi 2)$ and $\lam\in(-2,-1)$ and hence it is indeed in $H^1(\cS_0)$.
\end{proof}

Now we can finish our $H^3$ decomposition for $v_c$. Indeed, we have the elliptic equation for $v_r$ as   \[\na\cdot\cP_S\na v_r=h_1-c_s\na\cdot\cP_S \na S\in H^1(\cS)\] from Lemma \ref{S elliptic P_S system}. Combining this equation with  the tangential derivative $\na_{\tau_t} v_r\in H^2(\cS)$, a standard analysis shows that  the corresponding normal derivative $\na_{n_t}v_r$ is also in $H^2(S)$ with respect to $\na_{\tau_t} v_r$. This leads to the $H^3$ decomposition 
\[v=v_r+c_s S\] with
\[ v_r\in H^3(\cS),\quad \,S=\beta(r) r^{-\lam}\cos\big(\lam(\theta+\om_2)\big)\circ T^{-1}_0\in H^2\setminus H^3(\cS).\]

\subsubsection{Estimate for the singular coefficient $c_s$ and the regular part $v_r$} We haven't clarified  the singular coefficient $c_s$ so far. In order to have the $H^3$ estimate for $v_r$ we have to deal with $c_s$ first. In this section, the computations are done in the domain $\cS_0$ most of the time.

Now we take the corner domain $\cS_0$ as the  bounded triangle as mentioned before
\[\cS_0=\{(x,z)|\, 0\le x\le \del,\,-\ga x\le z\le \eta'(0)x\}\] since $\bar v=v\circ T_0$ and $\bar S$ are both  compactly supported near the corner and we only need to focus on the neighborhood of the corner. We will keep track of the small constant $\del$ through this section, which is the size of the corner neighborhood and plays an important role in the following estimates.

We are going to evaluate $c_s$ by using the system of $w_1=\na_{\tau_t} v_c$ in Lem \ref{tangent system}  with $b=0$, which is written in form of $\bar w_1=w_1\circ T_0$ on $\cS_0$:
\beq\label{system of w1}
\left\{\begin{array}{ll}\na\cdot\bar\cP_S\na \bar w_1=\bar l_1\qquad\hbox{on}\quad \cS_0\\
\bar w_1|_{\Ga_t}=0,\qquad \p^{\bar\cP_S}_{n_b}\bar w_1+b_1\p_{\tau_b} \bar w_1|_{\Ga_b}=\bar R_1\end{array}\right.
\eeq
where the constant $b_1=-\tan\om$ (instead of $\widetilde b$), since one has here $b=\tan\Phi=0$, and $\bar l_1=l_1\circ T_0$, $\bar R_1=R_1\circ T_0$.

Recalling that we have the decomposition
\[\bar w_1=\bar w_r+\bar c_s s,\quad\hbox{with}\quad w_r\in H^2(\cS_0),\,s\in H^1(\cS_0)\] with $\bar c_s=-\lam\big|(P_0^{-1})^t{\tau_t}\big| c_s$,  the idea of evaluating the singular coefficient $c_s$ is to apply the $L^2$ inner product on $\bar l_1$ with some chosen function $w$ to have an equality from the elliptic equation of $\bar w_1$.  Then we plan to put the $c_s$ part to the  left side to have some estimate for the right side. Obviously there will be $\bar w_r=\na_{\tau_t}\bar v_r$ part in the equality, which  also needs to be handled. Therefore we will choose a special  function $w$ such that in the equality the $w_r$ part almost vanishes.  In this way, we can evaluate $c_s$ first and then we use the estimate of $c_s$ to finish the estimate of $v_r$. The estimates in this part follows the idea of \cite{PG2}, and similar ideas can also be found in \cite{DNBL} etc..

First of all, we will look for this special function $w\in L^2(\cS_0)$. We  denote $\Del(V^2_{b_1})\subset L^2(\cS_0)$ as the image of $\Del$ on the space $V^2_{b_1}(\cS_0)$ (see the notations part). And we denote by  $\cN_1$ the nonempty orthogonal space  $\cN_{b_1}$ of $\Del(V^2_{b_1})$ in $L^2(\cS_0)$.  We'll choose  the special function $w$  from $\cN_1$. 

In fact we know from Lem \ref{space N} in Appendix B that when $\om\neq \f\pi 4$ (we have $b_1=-\tan \om$ here to find $\Phi=-\om$), one finds $w\in \cN_1$ iff $w\in L^2(\cS_0)$ and  satisfies
\[
\left\{\begin{array}{ll}\Del w=0\qquad\hbox{on}\quad \cS_0\\
w|_{\Ga_t}=0,\quad w|_{\Ga_\del}=0,\quad \p_{n_b}w-b_1\p_{\tau_b} w|_{\Ga_b}=0.\end{array}\right.
\] 
where the boundary conditions are in the following sense according to Thm \ref{Thm 1.5.2}
\[w|_{\Ga_t},\,w|_{\Ga_\del}\in (\tilde H^\f12)^*,\ \hbox{and}\ \p_{n_b} w,\,\p_{\tau_b} w|_{\Ga_b}\in (\tilde H^\f32)^*,\]where for the part $\p_{\tau_b} w|_{\Ga_b}$, the conclusion follows from $w|_{\Ga_b}\in (\tilde H^\f12)^*$ and the definition of the dual space.

We choose \[w=w_s-w_v,\] where
\[w_s(r,\theta)=\beta(r)\,r^{\lam+1}\cos{\big((\lam+1)(\theta+\om_2)-\om\big)}\in L^2(\cS_0)\setminus H^1(\cS_0)\]
defined from $s(r,\theta)$, and we know by a direct computation that $w_s$ satisfies the following system
\[
\left\{\begin{array}{ll}\Del w_s=[\Del,\be(r)] r^{\lam+1}\cos{\big((\lam+1)(\theta+\om_2)-\om\big)}\qquad \hbox{on}\quad \cS_0,\\
w_s|_{\Ga_t}=0,\quad w_s|_{\Ga_\del}=0, \quad \p_{n_b}w_s-b_1\p_{\tau_b}w_s|_{\Ga_b}=(\sin\om)\be'(r)r^{\lam+1}.\end{array}\right.
\]
On the other hand,   let $w_v\in V=\{u\in H^1(\cS_0)|\, u|_{\Ga_t}=u|_{\Ga_\del}=0\}$ be the variational solution of the same system as $w_s$, i.e. $w_v$ is defined by the variational equation for $\forall u\in V$:
\beq\label{var eqn for w_v} 
\int_{\cS_0}\na w_v\cdot\na u-b_1\int_{\Ga_b}\p_{\tau_b} w_v\,u=-\int_{\cS_0}\Del w_s \, u+\sin\om\int_{\Ga_b}\be'(r)r^{\lam+1}\, u.
\eeq
Similar as Lem \ref{variation solution} or applying Lax-Milgram Theorem directly, we find  that $w_v\in V$ exists and is unique. Note that on the boundary one has  $w_v,\,u|_{\Ga_b}\in \tilde H^\f12$ from Thm \ref{Thm 1.5.2.3} and the condition $w_v,\,u|_{\Ga_t,\Ga_\del}=0$,  and $\int_{\Ga_b}\p_{\tau_b} w_v\,u$ is well defined since $\p_{\tau_b}w_v\in (\tilde H^\f12)^*$ from Thm \ref{Remark 1.4.4.7}.

Therefore,  $w_v$ from the variational equation \eqref{var eqn for w_v} satisfies the estimate
\[
\|\na w_v\|_2\le C(\|\Del w_S\|_2+|\be'(r)r^{\lam+1}|_{L^2})
\le \del^{\lam+1}C
\] 
remembering that the size of $\cS_0$ (as a triangle) is $\del$. Applying the Poincar\'e inequality leads to the $L^2$ estimate for $w_v$ as follows
\beq\label{w_v estimate} \|w_v\|_{2}\le \del C\|\na w_v\|_{1,2}\le \del^{\lam+2}C\eeq for some constant $C$ depending on $\lam$, where we note that there is an extra $\del$ from the Poincar\'e inequality due to the size of $\cS_0$.

As a result, a direct computation tells us that  
\[
w=w_s-w_v\in\cN_1,\quad\hbox{when}\quad\om\neq\f\pi 4
\]
 and $w$ is nontrivial, since $w_s\in L^2\setminus H^1(\cS_0)$ and $w_v\in H^1(\cS_0)$ are nontrivial. Combining the expression of $w_s$ and the estimate \eqref{w_v estimate} for $w_v$  we finally arrive at the $L^2$ estimate for $w$:
\beq\label{w estimate L2}\|w\|_2\le \del^{\lam+2}C\eeq
For the moment, we are ready to evaluate $c_s$. The proposition below gives us the estimate for the singular coefficient $c_s$ and also the $H^3$ estimate for $v_r$.

\bthm{Proposition}\label{prop:H3}\label{H3 estimate} {\it 
Let  $h_c\in H^1(\cS)$ and $g_c\in H^\f32(\Ga_b)$ in the system (\ref{v1 problem}) for $v$. Moreover, when $\om=\f\pi 4$, one requires additionally the compatible condition: $R_1\in \tilde H^\f12(\Ga_b)$ from Lem \ref{tangent system} with $b=0, h=h_c, g=(1+\ga^2)^\f12 g_c$.\\
 Then there exists a constant $c_s$ such that the solution of (\ref{v1 problem}) has the following $H^3$ decomposition
\beno
v_c=v_r+\chi_1(\om)c_s S,
\eeno
where 
\[|c_s|\le \del^{\lam+2}C\big(\|h_c\|_{1,2}+|g_c|_{H^\f32}\big),\] and
\[
\|v_r\|_{3,2}\le C\big(\|h_c\|_{1,2}+|g_c|_{H^\f32}\big),
\] 
where $C=C(\ga,k,\lam,|\eta|_{H^3})$ and $\del>0$ small enough is the diameter of the corner neighborhood. }
\ethm

\begin{proof}
{\it Step 1}: the estimate for $c_s$.  {\it In this part, the contact angle $\om\in (\f\pi 4,\f \pi 2)$}.  Applying $L^2$ inner product on the equation in \eqref{system of w1} with $w\in\cN_1$, one has
\beq\label{c_s eqn}\begin{split}
\int_{\cS_0}w\,\bar l_1&= \int_{\cS_0}w\,\na\cdot \bar \cP_S\na\bar w_1\\
&=\int_{\cS_0}w\,\na\cdot \bar \cP_S\na\bar w_r+\bar c_s\int_{\cS_0}w\,\na\cdot\bar\cP_S\na s
:= A_1+\bar c_sA_2.
\end{split}
\eeq 
We'll deal with these terms one by one and we are going to show that $A_1$ is a  small term containing powers of $\del$.
\medskip

\noindent -$A_1$ estimate. In $A_1$, A direct computation from \eqref{system of w1} tells us that  $\bar w_r$ satisfies the system
\beq\label{w_r system}
\left\{\begin{array}{ll}
\na\cdot \bar\cP_S\na \bar w_r=\bar l_1-\bar c_s\na\cdot\bar\cP_S\na s\qquad\hbox{on}\quad \cS_0\\
\bar w_r|_{\Ga_t}=0,\quad \p_{n_b}^{\bar\cP_S}\bar w_r+b_1\p_{\tau_b} \bar w_r|_{\Ga_b}=\bar R_1-\bar c_s(\p^{\bar\cP_S}_{n_b}  s-\p_{n_b} s)|_{\Ga_b}+\bar c_s (\sin\om)\be'(r)r^{-\lam-1}
\end{array}\right.
\eeq
while noticing that 
\[
\p_{n_b}s+b_1\p_{\tau_b}s|_{\Ga_b}=-(\sin\om)\be'(r)r^{-\lam-1}
\]
 is used on the boundary. We have $\bar w_r\in H^2(\cS_0)$ from the decomposition of $\bar w_1$. Although $\bar w_r$ is not in $V^2_{b_1}$ yet, but we can find another function $\bar w_N\in V^2_{b_1}$ based on $\bar w_r$.

Here, we notice that $\bar R_1-\bar c_s(\p^{\bar\cP_S}_{n_b}  s-\p_{n_b}s)|_{\Ga_b}+\bar c_s (\sin\om)\be'(r)r^{-\lam-1}-(\p^{\bar\cP_S}_{n_b}  \bar w_r-\p_{n_b} \bar w_r)|_{\Ga_b}$ is compactly supported. Thus, by trace theorem, there exists a $\bar w_{T}\in H^2(\cS_0)$ satisfying the  boundary conditions $ \bar w_T|_{\Ga_t}=0$ and 
 \[
\p_{n_b}\bar w_T+b_1\p_{\tau_b} \bar w_T|_{\Ga_b}=\bar R_1-\bar c_s(\p^{\bar\cP_S}_{n_b}  s-\p_{n_b}s)|_{\Ga_b}+\bar c_s (\sin\om)\be'(r)r^{-\lam-1}-(\p^{\bar\cP_S}_{n_b}  \bar w_r-\p_{n_b} \bar w_r)|_{\Ga_b},
 \] 
and $\bar w_{T}$ is also compactly supported near the corner to have  $\bar w_T|_{\Ga_\del}=0$. Applying Theorem \ref{trace} and remembering $\om\in (\f\pi 4,\f\pi 2)$ here, we have
\[ 
\|\bar w_T\|_{2,2}\le C\big(|\bar R_1|_{H^\f12}+\bar c_s|\p^{\bar\cP_S}_{n_b} s-\p_{n_b} s|_{H^\f12}+|\p^{\bar\cP_S}_n  \bar w_r-\p_{n_b} \bar w_r|_{H^\f12}+\bar c_s|\be'(r)r^{-\lam-1}|_{H^\f12}\big).
\] 
Plugging the expression of $\bar R_1$ from Lem \ref{tangent system}(with $v=v$, $h=h_c$, $g=(1+\ga^2)^{\f12}g_c$, $\widetilde b=b_1$ and $b=0$),  and applying Remark \ref{trace on Ga j remark} lead to
\beq\label{trace w_T}
\begin{split}\|\bar w_T\|_{2,2}\le& C\Big(\|\bar v\|_{2,2}+\|\bar h_c\|_{1,2}+|\na\bar g_c|_{H^\f12}+\bar c_s\|(\bar \cP_S-Id)\na s\|_{1,2}\\
&\quad +\|(\bar \cP_S-Id)\na\bar w_r\|_{1,2}+\bar c_s|\be'(r)r^{-\lam-1}|_{H^\f12}\Big)
\end{split}\eeq
Writing 
\[
\bar w_N=\bar w_r-\bar w_T\in H^2(\cS_0),
\] 
we have that
\[
\bar w_N|_{\Ga_t}=\bar w_N|_{\Ga_\del}=0,\quad \p_{n_b}\bar w_N+b_1\p_{\tau_b}\bar w_N|_{\Ga_b}=0,
\] 
which means we find  the function $\bar w_N\in V^2_{b_1}$ finally (also compactly supported near $(0,0)$ indeed) and $\bar w_N$ satisfies
\[\Del \bar w_N=\Del\bar w_r-\Del\bar w_T=\na\cdot\bar\cP_S\na\bar w_r-\na\cdot(\bar\cP_S-Id)\na\bar w_r-\Del\bar w_T,\] which can be rewritten as
\[\na\cdot\bar\cP_S\na\bar w_r=\Del \bar w_N+\na\cdot(\bar\cP_S-Id)\na\bar w_r+\Del \bar w_T.\]

Consequently,  one can expand $A_1$ as
\beno A_1&=&
\int_{\cS_0}w\na\cdot \bar \cP_S\na\bar w_r\\
&=&\int_{\cS_0}w\Del \bar w_N +\int_{\cS_0}w\na\cdot(\bar\cP_S-Id)\na\bar w_r
+\int_{\cS_0}w\Del\bar w_T
\eeno 
where $\int_{\cS_0} w\Del \bar w_N =0$ since when $\om\in (\f\pi 4,\f\pi 2)$, we already have $w\in \cN_1$, $\bar w_N\in V^2_{b_1}$. So we arrive at the following estimate combining \eqref{trace w_T}
\beno
|A_1|&\le & \big(\|(\bar\cP_S-Id)\na\bar w_r\|_{1,2}+\|\bar w_T\|_{2,2}\big)\|w\|_2\\
&\le& C\Big(\|\bar v_c\|_{2,2}+\|\bar h_c\|_{1,2}+|\na\bar g_c|_{H^\f12}+\bar c_s\|(\bar \cP_S-Id)\na s\|_{1,2}+\|(\bar \cP_S-Id)\na\bar w_r\|_{1,2}\\
&&\qquad+\bar c_s|\be'(r)r^{-\lam-1}|_{H^1}\Big)\|w\|_2,
\eeno 
which together with the expression of $\bar R_1,\,s$ and \eqref{w estimate L2} implies
\beno
|A_1|&\le & C\Big(\del \|\p^3\bar v_r\|_2+\|\bar v_r\|_{2,2}+\|\bar h_c\|_{1,2}+|\na\bar g_c|_{H^\f12}+\bar c_s\|(\bar \cP_S-Id)\na s\|_{1,2}\\
&&\qquad +\bar c_s|\be'(r)r^{-\lam-1}|_{H^1}+\bar c_s\|\bar S\|_{2,2}\Big)\|w\|_2\\
&\le & C\big(\del^{\lam+3}\|\p^3\bar v_r\|_2+\del^{\lam+2}\|\bar v_r\|_{2,2}+\del^{\lam+2}\|\bar h_c\|_{1,2}+\del^{\lam+2}|\na\bar g_c|_{H^\f12}+\del \,\bar c_s\big)
\eeno 
where we use that $\bar w_r=\na_{\bar \tau_t}\bar v_r$, and  moreover \[
\|\bar S\|_{2,2}\le \del^{-\lam-1} C,\quad  |\be'(r)r^{-\lam-1}|_{H^1}\le \del^{-\lam-1} C,
\] and 
\[
\|\bar \cP_S-Id\|_{\infty}\le \del C,\quad \|(\bar\cP_S-Id)\na s\|_{1,2}\le \del^{-\lam-1} C.
\] The constant $C$ above depends on $|\eta|_{W^{2,\infty}}$.

\medskip

\noindent -$A_2$ estimate.   Rewriting $A_2$ as
\[A_2=\int_{\cS_0}w\na\cdot\bar\cP_S\na s=\int_{\cS_0}w\Del s+\int_{\cS_0}w\na\cdot(\bar\cP_S-Id)\na s,\] we will show that the first integral is equal to a constant with an explicit expression, and the second integral is a small part since $\bar\cP_S-Id$ is small.

In fact, the second integral of $A_2$ has  a similar estimate as $A_1$:
\[
\Big|\int_{\cS_0}\na\cdot(\bar\cP_S-Id)\na s\, w\Big|\le  \|(\bar\cP_S-Id)\na s\|_{1,2}\|w\|_2\le \del C(|\eta|_{W^{2,\infty}}).
\] 
Now we need to deal with the first integral carefully. Recalling that $w=w_s-w_v$, we have
\[
\int_{\cS_0}w\Del s=\int_{\cS_0}w_s\,\Del s-\int_{\cS_0}w_v\,\Del s,
\]
where  we know from Green's Formula and the variational equation \eqref{var eqn for w_v} of $w_v$ that
\beno
\int_{\cS_0}w_v\,\Del s&=& -\int_{\cS_0}\na w_v\cdot\na s+\int_{\Ga_b}w_v\,\p_{n_b} s\\
&=&-b_1\int_{\Ga_b}\p_{\tau_b} w_v \,s+\int_{\cS_0}\Del w_s\,s-\sin \om\int_{\Ga_b}\be'(r) r^{\lam+1}\,s+\int_{\Ga_b}w_v\,\p_{n_b} s
\eeno 
while noticing that $s\in V$ and so the integrals above make sense. Moreover, recalling that the singular function $s(r,\tht)=\be(r)r^{-\lam-1}\cos\big((\lam+1)(\tht+\om_2)-\om\big)$  satisfies
\[
\p_{n_b}s+b_1\p_{\tau_b} s|_{\Ga_b}=-(\sin\om) \be'(r)r^{-\lam-1},
\]
 and
\[
\int_{\Ga_b}\p_{\tau_b} w_v\,s=-\int_{\Ga_b} w_v\,\p_{\tau_b} s,
\] 
holds,  one can see that these two boundary integrals  together vanish:
\[
-b_1\int_{\Ga_b}\p_{\tau_b} w_v \,s+\int_{\Ga_b}w_v\,\p_{n_b} s=-\sin\om\int_{\Ga_b} w_v\,\be'(r)r^{-\lam-1}.
\]
 On the other hand, one can compute directly that
\[
\int_{\Ga_b}\be'(r) r^{\lam+1}\,s=\cos\om\int^\infty_0 \be'(r)r^{\lam+1}\,\be(r) r^{-\lam-1}dr =-\f12 \cos\om
\] 
since $\be(0)=1$. As a result, one arrives at
\[
\int_{\cS_0}w_v\,\Del s=\int_{\cS_0}\Del w_s\, s+\f12 \sin\om\cos\om-\sin\om\int_{\Ga_b} w_v\,\be'(r)r^{-\lam-1},
\] 
which finally leads to that
\[
\int_{\cS_0}w\,\Del s=\int_{\cS_0}w_s\Del s-\int_{\cS_0}\Del w_s\, s-\f12 \sin\om\cos\om+\sin\om\int_{\Ga_b} w_v\,\be'(r)r^{-\lam-1}.
\]
Let $B=\int_{\cS_0}w_s\Del s-\int_{\cS_0}\Del w_s \,s$. Since one can have similarly as before that 
\[
\left|\int_{\Ga_b} w_v\,\be'(r)r^{-\lam-1}\right|\le\|w_v\|_{1,2}|\be'(r)r^{-\lam-1}|_{L^2}\le \del C, 
\] we only need to deal with the remainder $B$ now.
Plugging into the expressions of $w_s,\,s$ and denoting $\al=(\lam+1)(\tht+\om_2)-\om$, a direct computation leads to
\beno
B&=& \int_{\cS_0}\Big[\be(r)r^{\lam+1}\cos\al\Del\big(\be(r)r^{-\lam-1}\cos \al\big)-\Del\big(\be(r)r^{\lam+1}\cos\al\big)\be(r)r^{-\lam-1}\cos \al\Big]\\
&=&\int^{\om_1}_{\om_2}\cos^2\al d\tht\int^\infty_0\Big[\be(r)r^{\lam+1}\big(\p^2_r+r^{-1}\p_r-(\lam+1)^2r^{-2}\big)(\be(r)r^{-\lam-1})\\
&&\quad -\big(\p^2_r+r^{-1}\p_r-(\lam+1)^2r^{-2}\big)(\be(r)r^{\lam+1})\be(r)r^{-\lam-1}\Big]rdr\\
&=& \int^{\om_1}_{\om_2}\cos^2\al d\tht\int^\infty_0\Big[\be(r)r^{\lam+1}\p_r\big(r\p_r(\be(r)r^{-\lam-1})\big)\\
&&\quad -\p_r\big(r\p_r(\be(r)r^{\lam+1})\big)\be(r)r^{-\lam-1}\Big]dr\\
&=&  \int^{\om_1}_{\om_2}\cos^2\al d\tht\Big[ \be(r)r^{\lam+1}\,r\p_r(\be(r)r^{-\lam-1})-r\p_r(\be(r)r^{\lam+1})\,\be(r)r^{-\lam-1}\Big]\Big|^\infty_0\\
&=& -2(\lam+1)\int^{\om_1}_{\om_2}\cos^2\al \,d\tht= -(\lam+1)\om-\sin\om\cos\om,
\eeno
 therefore one has that
\[
\int_{\cS_0}w\,\Del s=-c_0 +O(\del),
\] 
where $c_0=(\lam+1)\om+\f32\sin\om\cos\om>0$ is a constant and $O(\del)$ is a small term of order $\del$.

Consequently,  plugging the estimates for $A_1$ and $A_2$ and the expression of $\bar l_1$ from Lem \ref{tangent system}(with $v=v_c,\,h=h_c$) back into \eqref{c_s eqn} we derive that
\beno 
c_0\,|\bar c_s|&\le& |A_1|+
\|\bar l_1\|_2\|w\|_2+\del C\,\bar c_s\\
&\le& C\Big(\del^{\lam+3}\|\p^3\bar v_r\|_2+\del^{\lam+2}\|\bar v_r\|_{2,2}+\del^{\lam+2}\|\bar v\|_{2,2}+\del^{\lam+2}\|\bar h_c\|_{1,2}\\
&&\qquad +\del^{\del+2}|\na\bar g_c|_{H^\f12}+\del \bar c_s\Big)
\eeno 
where one has
\[
\|\bar v_r\|_{2,2}\le \|\bar v_c\|_{2,2}+\bar c_s\|\bar S\|_{2,2}\le \|\bar v\|_{2,2}+\del^{-\lam-1}C\,\bar c_s.
\]
Noting that $c_0$ is independent of $\del$, so for some fixed $\del$ small enough we arrive at
\[
\f12c_0\, \bar c_s\le C\big(\del^{\lam+3}\|\p^3\bar v_r\|_2+\del^{\lam+2}\|\bar v_c\|_{2,2}+\del^{\lam+2}\|\bar h_c\|_{1,2} +\del^{\del+2}|\na\bar g_c|_{H^\f12}\big).
\] 
Therefore plugging $H^2$ estimate of $v_c$ from Prop \ref{H2 estimate v1} (with $b=0$, $h=h_c$, $g=(1+\ga^2)^\f12 g_c$) into the inequality above and remembering $c_s=C\, \bar c_s$ for some constant $C$,  we conclude the estimate for the singular coefficient $c_s$ as follows
\beq\label{c_s estimate}
c_s\le C\big(\del^{\lam+3}\|\p^3\bar v_r\|_2+\del^{\lam+2}\|\bar h_c\|_{1,2}+\del^{\lam+2} |\bar g_c|_{H^\f32}\big)\eeq with $C=C(\ga,\lam,|\eta|_{H^3\cap W^{2,\infty}})$.

\noindent {\it Step 2}: the estimate for $v_r$.  In fact,  we have $\bar v_r=\bar v_c-c_s \bar S\in H^3(\cS_0)$ and satisfies the system based on \eqref{v1 problem} when $\om\in (\f\pi4,\f\pi2)$
\beq\label{v_r system}\left\{\begin{array}{ll}
\na\cdot \bar\cP_S\na \bar v_r=\bar h_c-c_s\na\cdot(\bar\cP_S-Id)\na\bar S-c_s\Del\bar S\qquad\hbox{on}\quad \cS_0\\
\bar v_r|_{\Ga_t}=0,\quad \p_{n_b}^{\bar\cP_S}\bar v_r|_{\Ga_b}=\bar g_c-c_s(\p^{\bar\cP_S}_{n_b} \bar S-\p_{n_b}\bar S)|_{\Ga_b}
\end{array}\right.
\eeq while noticing $\p_{n_b}\bar S|_{\Ga_b}=0 $ with $\bar S=\be(r)r^{-\lam}\cos\big(\lam(\tht+\om_2)\big)$.

Applying Prop \ref{H2 estimate v1} to the system \eqref{w_r system} of $\bar w_r=\na_{\tau_t}\bar v_r$, we  can find the $H^2$ estimate for $\bar w_r$.  After that,  we can go through a standard elliptic procedure to find  that
\beno
\|\bar v_r\|_{3,2}&\le& C\Big(\|\bar h_c\|_{1,2}+|\bar g_c|_{H^\f32}+c_s\|\na\cdot(\bar\cP_S-Id)\na\bar S\|_{1,2}+c_s\|\Del\bar S\|_{1,2}\\
&&\quad+c_s\big|(\p^{\bar\cP_S}_{n_b}-\p_{n_b})\bar S|_{\Ga_b}\big|_{H^\f32}\Big),
\eeno
which leads to by Thm \ref{trace on Ga j} and a careful computation on $\bar S$ as before that
\beno
\|\bar v_r\|_{3,2}&\le& C \big(\|\bar h_c\|_{1,2}+|\bar g_c|_{H^\f32}+ c_s\|(\bar\cP_S-Id)\na\bar S\|_{2,2}+ c_s\|\Del\bar S\|_{1,2}\big)\\
&\le& C\big( \|\bar h_c\|_{1,2}+|\bar g_c|_{H^\f32}+ \del^{-\lam-1}\,c_s\big).
\eeno 
Plugging the estimate \eqref{c_s estimate} for $c_s$ from Step 1 into the inequality above one finds that
\[
\|\bar v_r\|_{3,2}\le C\big(\|\bar h_c\|_{1,2}+|\bar g_c|_{H^\f32}+\del^{2}\|\p^3\bar v_r\|_2\big)
\]
 which implies that when $\del>0$ is small enough
\[
\|\bar v_r\|_{3,2}\le C\big(\|\bar h_c\|_{1,2}+|\bar g_c|_{H^\f32}\big)
\] 
 with $C=C(\ga,\lam,|\eta|_{H^3\cap W^{2,\infty}})$.
Plugging this inequality back into \eqref{c_s estimate} we also conclude the estimate for the singular coefficient $c_s$.

When $\om\in (0,\f\pi4]$, there is no singular part in $v$, so we take $v_r=v_c$ in the estimate above to have the $H^3$ estimate for $v_c$ directly.  Moreover, when $\om=\f\pi 4$, we need the compatible condition when having the $H^2$ estimate in the system \eqref{w_r system} for $\bar w_r$ from Prop \ref{H2 estimate v1} and Remark \ref{equiv comp condition}, which leads to 
\[
\bar R_1-\bar c_s(\p^{\bar\cP_S}_{n_b}  s-\p_{n_b} s)|_{\Ga_b}+\bar c_s (\sin\om)\be'(r)r^{-\lam-1}\in \tilde H^\f12(\Ga_b)
\] where  Thm \ref{Thm 1.5.2.3}  in Appendix A is applied also.  Since $\p^{\bar\cP_S}_{n_b}  s-\p_{n_b} s$ and $\be'(r)r^{-\lam-1}$ are in $\tilde H^\f12(\Ga_b)$ already, we find the compatible condition for $\bar w_r$ as $\bar R_1\in \tilde H^\f12(\Ga_b)$.  
\end{proof}

Going back to $\mbox{(MBVP)}$, we finally establish the $H^3$ estimate for the solution $u$ from $\mbox{(MBVP)}$.
\bthm{Proposition}\label{H3 decomp for MBVP}
{\it Let $h\in H^1(\Om)$, $f\in H^\f52(\Ga_t)$, $g\in H^\f32(\Ga_b)$ and $\eta \in H^{3}(\R^+)$ in $\mbox{(MBVP)}$. Moreover, when $\om=\f\pi 4$, one requires additionally the compatible condition $L_1(\p_{\tau_b}(u-u_b), h-\Del u_b)\in H^\f12(\Ga_b)$  where $L_1$ is a linear function of $\p_{\tau_b}(u-u_b)$,  $h-\Del u_b$ and $u_b\in H^3(\Om)$ satisfying the boundary condition of $\mbox{(MBVP)}$. \\
Then there exists a constant $c_s$ such that the solution $u$  of $\mbox{(MBVP)}$ has the following $H^3$ decomposition
\[
u=u_r+\chi_1(\om) c_s  \bar S\circ T^{-1}_0\circ (T_S)^{-1},
\]
where $\bar S(r,\theta)=\beta(r) r^{-\lambda}\cos(\lambda(\theta+\om_2))$ and $\lambda=-\f{\pi}{2\om}$. Moreover, one  has  estimates for the singular coefficient $c_s$ and the regular part $u_r$:
\[ 
\del^{-\lam-2}|c_s|+\|u_r\|_{3,2}\le C\big(\|h\|_{1,2}+|f|_{H^\f52}+|g|_{H^\f32}\big)
\]
where $C=C(\Ga_b, k, \lam,|\eta|_{H^3})$ and $\del>0$ small enough is the diameter of the corner neighborhood. }
\ethm
\begin{proof} As a first step, we  remove  the boundary conditions from $\mbox{(MBVP)}$. Applying Thm \ref{trace},  there exists  $u_b\in H^3(\Om)$ satisfying 
\[ u_b|_{\Ga_t}=f,\quad \p_{n_b}u_b|_{\Ga_b}=g.\] Letting $u_h=u-u_b$, we find $u_h\in H^2(\Om)$ be the solution to $\mbox{(MBVP)}_h$. Consequently, we have the system \eqref{v1 problem} for $v_c=\be u_h\circ T_S$ with  
\[h_c=\Big(\be(h-\Del u_b)-[\be,\Del]u_h\Big)\circ T_S,\quad g_c=-\big(\p_{n_b}\be \,u_h|_{\Ga_b}\big)\circ T_S.\]
Using Prop \ref{prop:H3} for $v_c$, one can derive the result for $v_c$ and then go back to the estimate for $u_c=\be u_h$. 
On the other hand, for $u_{ac}=(1-\be)u_h$ system, we can use $T_R$ to transform it into an elliptic system on a flat strip, which can be estimated by a standard procedure (see \cite{Lannes} for example) and therefore omitted here. Combining these two parts together, we find the $H^3$ decomposition and estimate for $u$.  Note that the singular part appears only  near the corner, so there is no singular part in $u_{ac}$. Moreover, the singular function $\bar S$ doesn't depend on the right side of the boundary conditions $f,\,g$, while the singular coefficient $c_s$ is effected by the right side of the system.

In the end, the compatible condition on $R_1$ needs to be stressed here. In fact, with the expressions of $h_c, g_c$, we know that $g_c\in \tilde H^\f12(\Ga_b)$, since there is $\p_{n_b}\be$ in $g_c$ to make $g_c$ vanish near the corner. Similarly we also have $[\be,\Del]u_h\circ T_S|_{\Ga_b}\in\tilde H^\f12(\Ga_b)$ from $h_c$. As a result, the compatible condition $R_1\in \tilde H^\f12(\Ga_b)$ becomes
\[\rho_1\p_xv_c-\tau_2 \be(h-\Del u_b)\circ T_S|_{\Ga_b}\in \tilde H^\f12(\Ga_b),\] which can be rewritten as the condition in $\Om$:
\[\rho_1(1+\ga^2)^\f12\be \p_{\tau_b}u_h-\tau_2\be (h-\Del u_b)\in \tilde H^\f12(\Ga_b)\] since $\p_x v_c=\p_{\tau_b}v_c$ on $\Ga_b$ in domain $\cS$. Consequently,  we can note the compatible condition as $L_1(\p_{\tau_b}(u-u_b), h-\Del u_b)\in H^\f12(\Ga_b)$,  where $L_1$ is a linear function of  $\p_{\tau_b}(u-u_b)$ and $ h-\Del u_b$.
\end{proof}

\medskip

\subsection{Higher-order decompositions and regularities}
After studying carefully the first singular function $\bar S=\be(r) r^{-\lam}\cos\big(\lam(\tht+\om_2)\big)\in H^2\setminus H^3(\cS_0)$  and its coefficient $c_s$ from $H^3$ decomposition, we are ready to extend our theory to higher-order cases.

\subsubsection{$H^4$ decomposition and estimate}
We assume that $h\in H^2(\cS)$, $f\in H^\f72(\Ga_t)$ and $g\in H^\f52(\Ga_b)$ in $\mbox{(MBVP)}$ in this part. The idea of $H^4$ decomposition is to take $\p^2_{\tau_t}$ on system  \eqref{v1 problem} for $v_c$ and study the tangential system as we did in $H^3$ decomposition. After that, we go back to $\mbox{(MBVP)}$ as before.

Recall that in $H^3$ case we showed  $\na\cdot\bar \cP_s\na \bar S\in H^1$  in Lemma \ref{S elliptic P_s system} thanks to the fact that $\bar\cP_s$ is close to $Id$, which enables us to finish the $H^3$ decomposition. Compared to the rectilinear polygonal problem with $\Del$ operator in \cite{PG1} where they have $\Del \bar S\in H^\infty$, we only find a very limited regularity ($H^1$ indeed) for $\na\cdot\bar \cP_s\na \bar S$. As a consequence,  in order to improve  the  regularity to $H^4$, we have to modify $\na\cdot\bar \cP_s\na \bar S$ by some asymptotic analysis.

\bthm{Lemma}\label{S_{1,r,1}}{\it  There exist functions  $S_{1,r,i}=\be(r)r^{-\lam+i}a_i(\tht)\in H^{i+2}\setminus H^{i+3}(\cS_0)$ with $i=1,2,\dots,l$ satisfying
\[\na\cdot\bar \cP_S\na \Big(\bar S+\sum_{i=1}^l\bar S_{1,r,i}\Big)\in H^{l+1}(\cS_0)\] and 
\[
\bar S_{1,r,i}|_{\Ga_t}=0,\qquad \p^{\bar\cP_S}_{n_b}\Big(\bar S+\sum_{i=1}^l\bar S_{1,r,i}\Big)|_{\Ga_b}\in H^{l+\f32}(\Ga_b)
\]
 with $a_i(\tht)$  bounded functions of $\tht\in[-\om_2,\,\om_1]$. 
}\ethm
\begin{proof} The proof is a standard approximation analysis. We prove in details only the case when $l=1$.  Firstly,  we would  expand $\bar\cP_s$ with respect to $r$.  In fact, $\bar\cP_S$ can be rewritten as
\[
\bar \cP_S=Id+\left(\begin{matrix}2\ga &1-\ga^2\\ 1-\ga^2 &-2\ga\end{matrix}\right)(\bar d-d_0)+(1+\ga^2)\left(\begin{matrix}1 &-\ga\\-\ga & \ga^2\end{matrix}\right)(\bar d-d_0)^2
\]
with $d_0=d(z)|_{z=0}=\bar d(r,\tht)|_{r=0}$ and
\[
\bar d(r,\tht)= d\circ T_0=\f 1k-\f1{\ga+\eta'\big(\bar\eta^{-1}(r(\ga\cos\tht+\sin\tht))\big)}.
\]
Moreover one has the following Taylor approximation with respect to $r$
\[
\bar d(r,\tht)=d_0+d_1(\tht)\,r+O(r^2)\quad\hbox{with}\quad d_1(\tht)=-\f{\eta''(0)(\ga\cos\tht+\sin\tht)}{(\ga+\eta'(0))^3},
\] 
so it's enough to show that
\[
\na\cdot\Big(Id+\left(\begin{matrix}2\ga &1-\ga^2\\ 1-\ga^2 &-2\ga\end{matrix}\right) d_1(\tht) r\Big)\na(\bar S+\bar S_{1,r,1})\in H^2(\cS_0)
\] 
in order to finish the proof. Equivalently, we need that
\[
A=\na\cdot\Big(Id+\left(\begin{matrix}2\ga &1-\ga^2\\ 1-\ga^2 &-2\ga\end{matrix}\right) d_1(\tht) r\Big)\na\Big(r^{-\lam}\cos\big(\lam(\tht+\om_2)\big)+r^{-\lam+1}a_1(\tht)\Big)
\] 
is in $H^2(\cS_0)$ when we get rid of the cut-off fucntion $\be(r)$. Indeed, a direct computation leads to
\beno
A&=&\Del(r^{-\lam}\cos\al)+\na\cdot d_1(\tht)r\left(\begin{matrix}2\ga &1-\ga^2\\ 1-\ga^2 &-2\ga\end{matrix}\right)\na(r^{-\lam}\cos\al)+\Del(r^{-\lam+1}a_1(\tht))\\
&&\quad+\na\cdot d_1(\tht)r\left(\begin{matrix}2\ga &1-\ga^2\\ 1-\ga^2 &-2\ga\end{matrix}\right)\na(r^{-\lam+1}a_1(\tht))
\eeno 
where we write $\al=\lam(\tht+\om_2)$ for short. We can see that $A\in H^2$ if the following equation holds:
\beno0&=&\Del(r^{-\lam+1}a_1(\tht))+\big(\na d_1(\tht)\big)\cdot r\left(\begin{matrix}2\ga &1-\ga^2\\ 1-\ga^2 &-2\ga\end{matrix}\right)\na(r^{-\lam}\cos\al)\\
&&\quad+d_1(\tht)\na\cdot r\left(\begin{matrix}2\ga &1-\ga^2\\ 1-\ga^2 &-2\ga\end{matrix}\right)\na(r^{-\lam}\cos\al),\eeno which can be rewritten as
\beq\label{a eqn}\begin{split} a''_1(\tht)+(\lam-1)^2a_1(\tht)=&C_0(1+\ga^2)^\f32\Big[\lam(\lam+1)\cos\big((\lam+3)(\tht+\om_2)\big)\\
&\qquad -\lam^2\cos\big((\lam+1)(\tht+\om_2)\big)\Big]\end{split}\eeq with  the notation $d_1(\tht)=C_0(\ga\cos\tht+\sin\tht)$ and the polar forms $\Del=\p^2_r+r^{-1}\p_r+r^{-2}\p^2_\tht$ and $\na=(\cos\tht\p_r-r^{-1}\sin\tht\p_\tht,\,\sin\tht+r^{-1}\cos\tht\p_\tht)^t$.

Now we turn to the boundary condition $\p^{\bar\cP_s}_{n_b}(\bar S+\bar S_{1,r,1})|_{\Ga_b}\in H^\f52$. After a similar computation, we find that this condition is satisfied when
\[
\p_n(r^{-\lam+1}a_1(\tht))+C_0{\bf n}\cdot (\ga\cos\tht+\sin\tht)r\left(\begin{matrix}2\ga &1-\ga^2\\ 1-\ga^2 &-2\ga\end{matrix}\right)\na(r^{-\lam}\cos\al)\Big|_{\Ga_b}=0,
\] where $\Ga_b$ means $\tht=-\om_2$ and $\p_n=r^{-1}\p_\tht$, so we arrive at the boundary condition
\beq\label{a Ga_b} a'_1(-\om_2)=0.\eeq Moreover, we also require $\bar S+\bar S_{1,r,1}|_{\Ga_t}=0$, which implies
\beq\label{a Ga_t} a_1(\om_1)=0.\eeq Putting \eqref{a eqn}, \eqref{a Ga_b} and \eqref{a Ga_t} together, we derive the initial value problem for $a_1(\tht)$ with $\tht\in [-\om_2,\,\om_1]$, and a direct computation tells us that the solution is
\beq\label{a expression} \begin{split}a_1(\tht)=&\bar C_0\Big[\Big(2-\f{\sin 3\om}{\sin\om}\Big)\cos\big((\lam-1)(\tht+\om_2)\big)+2\cos\big((\lam+1)(\tht+\om_2)\big)\\
&\qquad -\cos\big((\lam+3)(\tht+\om_2)\big)\Big]\end{split}\eeq where we note $\bar C_0=\f18\lam C_0(1+\ga^2)^\f32$.

For the general case $S_{1,r,l}=\be(r)r^{-\lam+l}a_l(\tht)$ ($l\ge 2$), following the procedure above, we can prove by an induction argument that, the unknown function $a_l(\tht)$ satisfies
\beq\label{a_l system}
\left\{\begin{array}{ll}\Del\big(r^{-\lam+l}a_l(\tht)\big)+\na\cdot\bar\cP_{S1}\na\Big(r^{-\lam}\cos\al+\sum^{l-1}_{i=1}r^{-\lam+i}a_i(\tht)\Big)=0,\\
\p_n\big(r^{-\lam+l}a_l(\tht)\big)+\p^{\bar\cP_{S1}}_n\big(r^{-\lam}\cos\al+\sum^{l-1}_{i=1}r^{-\lam+i}a_i(\tht)\big)|_{\tht=-\om_2}=0,\\
a_l(\om_1)=0
\end{array}\right.
\eeq where we denote \[\bar\cP_{S1}=Id+\left(\begin{matrix}2\ga &1-\ga^2\\ 1-\ga^2 &-2\ga\end{matrix}\right) d_1(\tht) r.\]  Since this 2nd-order ODE is always solvable, we can find the solution $a_l(\tht)$ in a similar expansion as $a_1(\tht)$ and bounded.
\end{proof}
\bthm{Remark}\label{general modify}{\it For a general singular function $\bar S_i=\be(r)r^{-\lam_i}\cos\big(\lam_i(\tht+\om_2)\big)$, we have similar modifying functions $\bar S_{i,r,j}=\be(r) r^{-\lam_i+j}a_{i,j}(\tht)$ with some bounded function $a_{i,j}$ satisfying similar problem as  \eqref{a_l system}. On the other hand, based on the conclusion on the boundary, one can have indeed $\p^{\bar\cP_S}_{n_b}\Big(\bar S+\sum_{i=1}^l\bar S_{1,r,i}\Big)|_{\Ga_b}\in \tilde H^{l+\f32}(\Ga_b)$ from the proof.
}
\ethm

As a result, we can rewrite the $H^3$ decomposition for $v_c$ as
\[ v_c=v_{r1}+\chi_1(\om)c_s(S+S_{1,r,1})\] where $ v_{r1}=v_r-\chi_1(\om)c_s S_{1,r,1}\in H^3(\cS)$ and  $ S_{1,r,1}=\bar S_{1,r,1}\circ T^{-1}_0\in H^3(\cS)$. We'll use this modified decomposition and take $v_{r1}$ as the regular part instead of $v_r$. Similarly as in $H^3$ decomposition, we will consider the $H^4$ decomposition and estimates based on $v_{r1}$. In fact,  from \eqref {v1 problem} one finds the system for $v_{r1}$ as follows
\beq\label{v_r1 system}
\left\{\begin{array}{ll}
\na\cdot\cP_S\na v_{r1}= h_c-\chi_1(\om)c_s\na\cdot\cP_S\na(S+S_{1,r,1})\quad \hbox{on}\quad \cS,\\
v_{r1}|_{\Ga_t}=0,\quad \p^{\cP_S}_{n_b} v_{r1}|_{\Ga_b}=(1+\ga^2)^\f12 g_c-\chi_1(\om)c_s\p^{\cP_S}_{n_b}( S+ S_{1,r,1})|_{\Ga_b}.
\end{array}\right.
\eeq 
Similarly as in $H^3$ decomposition,  when $\om\neq \f{\pi}4$, we take twice the tangential derivative $\na_{\tau_t}$on \eqref{v_r1 system}  by Lemma \ref{tangent system} to  derive the system for $w_2=\na^2_{\tau_t}v_{r1}\in H^1(\cS)$ as follows
\beq\label{w_2 system}\left\{\begin{array}{ll}
\na\cdot\cP_S\na w_2= l_2\qquad \hbox{on}\quad \cS,\\
w_2|_{\Ga_t}=0,\quad \p^{\cP_S}_{n_b} w_2+b_2\p_{\tau_b} w_2|_{\Ga_b}=R_2
\end{array}\right.\eeq where $b_2=-\tan(2\om)$  and
\beno
l_2&=&-\na\cdot(\na_{\tau_t}\cP_S)\na(\na_{\tau_t} v_{r1})-\na_{\tau_t}\big(\na\cdot(\na_{\tau_t}\cP_S)\na  v_{r1}\big)\\
&&\quad+\na^2_{\tau_t}\Big(h_c-\chi_1(\om)c_s\na\cdot\cP_S\na(S+S_{1,r,1})\Big)\\
R_2&=&\rho_1\p_x(\na_{\tau_t}v_{r1})|_{\Ga_b}\\
&&\quad+\tau_2\Big(\na\cdot(\na_{\tau_t}\cP_S)\na  v_{r1}-\na_{\tau_t}\big(h_c-\chi_1(\om)c_s\na\cdot\cP_S\na(S+S_{1,r,1})\big)\Big)|_{\Ga_b}\\
&& \quad+\rho_2\rho_1\p^2_x v_{r1}|_{\Ga_b}-\rho_2\tau_2\,\p_x\Big(h_c-\chi_1(\om)c_s\na\cdot\cP_S\na(S+S_{1,r,1})\Big)|_{\Ga_b}\\
&&\quad +\rho_2\tilde\rho_2\,\Big((1+\ga^2)^\f12 \p^2_xg_c-\chi_1(\om)c_s\p^2_x\p^{\cP_S}_{n_b} (S+S_{1,r,1})|_{\Ga_b}\Big)
\eeno 
with constants $\tau_1,\tau_2, \rho_1$ the same as in Lemma \ref{tangent system} and 
\beno
&&\rho_2=\tau_1-\f1{\ga+(1+\ga^2)d_0+\tan\om}\Big(\tau_1(\tan\om-\tan2\om)+\tau_2\big((1+\ga d_0)^2+d^2_0\big)\Big), \\
&&\tilde\rho_2=\tau_1-\f1{\ga+(1+\ga^2)d_0}\Big(-\tau_1\tan\om+\tau_2\big((1+\ga d_0)^2+d^2_0\big)\Big).
\eeno

Recalling that $h_c, g_c$ from \eqref{v1 problem}, we know from Lemma \ref{S_{1,r,1}} and the assumption $h\in H^2(\cS),g\in H^\f52(\Ga_b)$ that $l_2\in  L^2(\cS)$ and $R_2\in  H^\f12(\Ga_b)$,  while noticing that $w_2\in H^1(\cS)$. Therefore we can apply Thm \ref{H2 decompose} to \eqref{w_2 system} to have $H^2$ decomposition for $w_2$ in  case when $\om\neq\f\pi4$:
\[
w_2=w_{r2}+\chi_2(\om)\bar c_2 s_2
\]
with
\[
w_{r2}\in H^2(\cS),\quad s_2\in H^1\setminus H^2(\cS),\quad\chi_2(\om)=\left\{\begin{array}{ll}1,\quad \om\in (\f\pi 6,\f\pi 4)\\
0,\quad \hbox{otherwise}.\end{array}\right.
\]
The expression of $s_2$ in domain $\cS_0$  is
\[
s_2(r,\tht)=\be(r)r^{-\lam-2}\cos\big((\lam+2)(\tht+\om_2)-2\om\big),
\] 
which appears only when $\om\in (\f\pi 6,\f\pi 4)$.  Moreover, when $\om=\f\pi 6$, the compatible condition 
\beq\label{comp H4 1}
R_2\in \tilde H^\f12(\Ga_b)
\eeq is needed, according to Thm \ref{H2 decompose} and Thm \ref{Thm 1.5.2.3} in Appendix A.

In case when $\om=\f\pi4$, we apply Lem \ref{tangent system} first on $\na_{\tau_t}v_{r1}$ and then Lem \ref{tangent system special} on $w_2=\na^2_{\tau_t}v_{r1}$ to find the system for $w_2$ as 
\beq\label{D system for  w2}\left\{\begin{array}{ll}
\na\cdot\cP_S\na w_2= l_2\quad \hbox{on}\quad \cS,\\
w_2|_{\Ga_t}=0,\quad w_2|_{\Ga_b}=g_2
\end{array}\right. 
\eeq
where 
\[
\begin{split}g_2&=-\tau_2(1+\ga^2)^{-1}\Big[\rho_1\p_xv-\tau_2\big(h_c-\chi_1(\om)c_s\na\cdot\cP_S\na(S+S_{1,r,1})\big)\\
&\qquad+\rho_2(1+\ga^2)^\f12 \p_xg_c-\rho_2\chi_1(\om)c_s\p_x\p^{\cP_S}_{n_b}(S+S_{1,r,1})\Big]_{\Ga_b}\end{split}
\] Meanwhile the compatible condition \eqref{compatible condition S} is required here,  which results in by Thm \ref{Thm 1.5.2.3} in Appendix A that 
\beq\label{comp H4 2}
g_2\in \tilde H^\f12(\Ga_b).
\eeq
We need the $H^2$ decomposition for the Dirichlet-conditions case. In fact, Thm \ref{H2 decompose 1} tells us (it is well-known already), when the contact angle $\om\in (0,\,\f\pi2)$, there is no singular part for $w_2$, which infers that $w_2\in H^2(\cS)$ in the special case when $\om=\f\pi4$.

In a word,   when the contact angle $\om\notin (\f\pi 6,\f\pi 4)$, we have $w_2\in H^2(\cS)$.

\medskip
Now it's the time to find the corresponding singular part for $v$. In fact,  a direct computation leads to that 
\[
\na^2_{\tau_t} S=\lam^2\big|(P^{-1}_0)^t\tau_t\big|^2 s_2\circ T^{-1}_0+ [\na^2_{\tau_t},\be(r)\circ T^{-1}_0]r^{-\lam}\cos\big(\lam(\tht+\om_2)\big)\circ T^{-1}_0\quad\hbox{on}\quad \cS
\] 
while recalling $\bar S=\be(r)r^{\lam}\cos\big(\lam(\tht+\om_2)\big)$ defined on $\cS_0$ and $S=\bar S\circ T^{-1}_0$ defined on $\cS$. Consequently,  we find that
\beq\label{w_r2 relation}
\na^2_{\tau_t}(v_{r1}-c_{s2} S)=w_{r2}-c_{s2}[\na^2_{\tau_t},\be(r)\circ T^{-1}_0]r^{-\lam}\cos\big(\lam(\tht+\om_2)\big)\circ T^{-1}_0\in H^2(\cS)
\eeq 
where the coefficient $c_{s2}=\lam^{-2}\big|(P^{-1}_0)^t\tau_t\big|^{-2} \bar c_2$. Setting
\[v_{r2}=v_{r1}-c_{s2} S,\] we now have  $\p^2_{\tau_t}v_{r2}\in H^2(\cS)$ already.  We can show by a standard elliptic analysis based on the system \eqref{v_r1 system} of $v_{r1}$ that
\[v_{r2}\in H^4(\cS).\] As a result, we finally derive the $H^4$ decomposition for $v$: 
\[ v_{r1}=v_{r2}+\chi_2(\om)c_{s2}S \] or equivalently
\[v_c=v_{r1}+\chi_1(\om)c_s( S+S_{1,r,1})=v_{r2}+\chi_1(\om)c_s(S+S_{1,r,1})+\chi_2(\om)c_{s2} S.\]
\bthm{Remark}\label{small angle om}{\it Note in the case when both $\chi_1,\,\chi_2$ are zero, i.e. the contact angle $\om\in (0,\,\f\pi 6)$, we find that $v_c\in H^4(\cS)$. This fact tells us that smaller contact angle implies smoother solution from (MBVP), which agrees with previous literature. The case when $\om=\f\pi 6$ is also a special case as $\om=\f\pi 4$, which will be discussed later. }\ethm
\bthm{Proposition}\label{prop:H4}($H^4$ estimate){\it Assume that $h_c\in H^2(\cS)$, $g_c\in H^\f52(\Ga_b)$ for the system \eqref{v1 problem} of $v_c$. Moreover one requires additionally compatible conditions \eqref{comp H4 1} when $\om=\f\pi 6$ and \eqref{comp H4 2} when $\om=\f \pi 4$. \\
Then one has the $H^4$ decomposition
\[
v_c=v_{r2}+\chi_1(\om)c_s(S+S_{1,r,1})+\chi_2(\om)c_{s2} S
\] 
whith the regular part $v_{r2}\in H^4(\cS)$ and  the estimate for both the singular coefficient and the regular part
\[
\del^{-\lam-3}|c_{s2}|+\|v_{r2}\|_{4,2}\le C(\|h_c\|_{2,2}+|g_c|_{H^\f52})
\] 
where the constant $C=C(\ga,k,\lam,|\eta|_{H^4})$.}
\ethm
\begin{proof}
The proof is similar as the proof of Propostion \ref{prop:H3}. So we write down here only the sketch of the proof. Throughout the proof, we will denote by $\bar f=f\circ T_0$ on $\cS_0$ corresponding to a function $f$ defined on $\cS$.

In order to evaluate $v_{r2}$, we only need to have estimate for $c_{s2}$ here since we've dealt with $c_s$ already. In this case we find the contact angle 
\[
\om\in (\f\pi 6,\f\pi 4),\quad\hbox{ and}\quad \lam=-
\f\pi{2\om}\in (-3,-2).
\]

To begin with, we need again to find a function $w$ in $\cN_2$, which is the orthogonal space of $\Del(V^2_{b_2})$ in $L^2(\cS_0)$. Recall from the notation part that
\[V^2_{b_2}=\{u\in H^2(\cS_0)|\,u|_{\Ga_t}=u|_{\Ga_\del}=0,\,\p_{n_b}u+b_2\p_{\tau_b} u|_{\Ga_b}=0\}.\] We define
\[w=w_s-w_v\] where $w_s=\be(r)r^{\lam+2}\cos\big((\lam+2)(\tht+\om_2)-2\om)\big)$ and $w_v$ the variational solution of the same system of $w_s$. We can show that $w\in \cN_2$.

On the other hand, since one has when $\om\in (\f\pi 6,\f\pi 4)$ the decomposition
\[v_c=v_{r2}+c_{s2}S\]  and  $w_2=\na^2_{\tau_t}v_{r1}=\na^2_{\tau_t} v_c$, one finds   by \eqref{w_r2 relation} that \[w_2=w_{r2}+\bar c_2 s_2, \] with $\bar c_2=\lam^2|(P^{-1}_0)^t\tau_t|^2\,c_{s2}$.

So recalling the system \eqref{w_2 system} with no singular part and $v_{r1}$ replaced by $v_c$, we have the following $L^2$ product as in the proof of Propostion \ref{prop:H3}
\beno\int_{\cS_0}\bar l_2\, w&=&\int_{\cS_0}\na\cdot\bar\cP_S\na \bar w_2\, w\\
&=&\int_{\cS_0}\na\cdot\bar\cP_S\na\bar w_{r2}\, w+\bar c_2\int_{\cS_0}\na\cdot\bar\cP_S\na s_2\, w.\eeno
Similar computations lead to the estimate for $c_{s2}$ as below
\[c_{s2}\le C\big(\del^{\lam+4}\|\p^2\bar v_{r2}\|_2+\del^{\lam+3}\|\bar v_c\|_{3,2}+\del^{\lam+3}\|\bar h_1\|_{2,2}\big).\] With this estimate one can conclude the proof by a standard elliptic analysis for $v_{r2}$.
\end{proof}
\medskip
Going back to the solution $u$ of $\mbox{(MBVP)}$ and combining Thm \ref{trace}, we have the following proposition similar as Prop \ref{H3 decomp for MBVP}. Since we have the same function $\bar S$ here for different $\om$ and $\lam$, we mark the two different versions as $\bar S_1, \bar S_2$ for clarity, and we also rewrite $c_s$ as $c_{s1}$.
\bthm{Proposition}{\label{H4:u}}
{\it Let $h\in H^2(\Om)$, $f\in H^\f72(\Ga_t)$, $g\in H^\f52(\Ga_b)$ and $\eta\in H^4(\R^+)$. One requires additionally compatible conditions $L_2\big(\p^2(u-u_b),\p(u-u_b),\p(h-\Del u_b)\big)\in \tilde H^\f12(\Ga_b)$ when $\om=\f\pi 6$ and $L_1\big(\p_{\tau_b}(u-u_b),h-\Del u_b\big)\in \tilde H^\f12(\Ga_b)$ when $\om=\f\pi4$, where $L_2$ is a linear function of $\p^2(u-u_b)$, $\p(u-u_b)$, $\p(h-\Del u_b)$ and $u_b$, $L_1$ as in Prop \ref{H3 decomp for MBVP}. \\
Then there is the $H^4$ decomposition for the solution $u$ of $\mbox{(MBVP)}$
\[
u=u_{r2}+u_{s2},
\] with $u_{r2}$ the regular part and $u_{s2}$ the singular part
\[
u_{s2}=\chi_1(\om)c_{s1} (\bar{S}_1+\bar{S}_{1,r,1})\circ T_0^{-1}\circ (T_S)^{-1}+\chi_2(\om)c_{s2}\,\bar{S}_2\circ T_0^{-1}\circ (T_S)^{-1},
\]
where 
\[
\bar{S}_1=\bar S_2=\bar S=\beta(r) r^{-\lambda}\cos(\lambda(\theta+\om_2))\quad\hbox{for different}\ \  \om \ \hbox{in}\  \lam=-\f{\pi}{2\om},
\]
  $\bar{S}_{1,r,1}(r,\theta)=\beta(r) r^{-\lambda+1}a_1(\theta)$,  and $a_1(\theta)$ is a bounded function for $\theta\in [-\om_2, \om_1]$ from Lem \ref{S_{1,r,1}}. Moreover, one has the following estimate for the regular part
\beno
\|u_{r2}\|_{4,2}\leq C\big(\|h\|_{2,2}+|f|_{H^\f72}+|g|_{H^\f52}\big),
\eeno
with the constant  $C=C(\Ga_b,k,\lam,|\eta|_{H^4})$.
}
\ethm
\begin{proof} Similarly as in Prop \ref{H3 decomp for MBVP}, we only need to deal with the compatible condition here.   In fact, we have 
\[
h_c=\Big(\be(h-\Del u_b)-[\be,\Del]u_h\Big)\circ T_S,\quad g_c=-\big(\p_{n_b}\be \,u_h|_{\Ga_b}\big)\circ T_S.
\] for the system \eqref {v1 problem} of $v_c=\be u_h\circ T_S$, where $u_b\in H^3(\Om)$ satisfying the boundary conditions 
\[ u_b|_{\Ga_t}=f,\quad \p_{n_b}u_b|_{\Ga_b}=g.\] Moreover,  $u_h=u-u_b$ is the solution to $\mbox{(MBVP)}_h$.  

When $\om=\f\pi 4$, the compatible condition $g_2\in \tilde H^\f12(\Ga_b)$ turns out to be the same as in Prop \ref{H3 decomp for MBVP}: $L_1(\p_{\tau_b}u, h-\Del u_b)\in H^\f12(\Ga_b)$  where $L_1$ is a linear function of $\p_{\tau_b}u$,  $h-\Del u_b$.

When $\om=\f\pi 6$,  we need $R_2\in\tilde H^\f12(\Ga_b)$ from  Prop \ref{prop:H4}.  Checking \eqref{w_2 system} for $R_2$, the term $\na\cdot (\na_{\tau_t}\cP_S)\na v_{r1}$ can be rewritten as 
\[
\na\cdot (\na_{\tau_t}\cP_S)\na v_{r1}=\sum_{|\al|=2}b_\al(\p\cP_S) \p^\al v_{r1}+\sum_{|\al|=1}c_{\al}(\p^2\cP_S)\p v_{r1}
\] where $b_\al(\p\cP_s)$, $c_{\al}(\p^2\cP_S)$ are functions of $\p\cP_S, \p^2\cP_S$ and $\p^2$ means 2nd-order partial derivative. In fact, when we have 
\[
\sum_{|\al|=2}b_\al(\p\cP_S|_{(0,0)}) \p^\al v_{r1}+\sum_{|\al|=1}c_{\al}(\p^2\cP_S|_{(0,0)})\p v_{r1}\in \tilde H^\f12(\Ga_b),
\] we will also find $\na\cdot (\na_{\tau_t}\cP_S)\na v_{r1}\in \tilde H^\f12(\Ga_b)$. So all the coefficients for $v_{r1}$ can be  taken as  constants in the compatible condition.

As a result, the compatible condition reads 
\[L_2\big(\p^2(u-u_b),\p (u-u_b), \p(h-\Del u_b)\big)\in \tilde H^\f12(\Ga_b), \] where  and $L_2$ is a linear function of $\p^2(u-u_b)$, $\p (u-u_b)$ and $\p(h-\Del u_b)$.

\end{proof}

\subsubsection{$H^5$ and higher-order decompositions and estimates}
When we consider the $H^{2+K}$ ($K\in \N$) decomposition and regularity, the key point is to find out the number of singular functions, which depends on the eigenvalue $\lam_m$ in Thm \ref{H2 decompose}. In fact, at this time,  we have $b_K=-\tan(K\om)$ from the  boundary conditions after taking the derivative $\na^K_{\tau_t}$ on $v$, so we find 
\beq\label{range of lambda}
\lam_m=\f{K\om+(m+\f12)\pi}{\om}=K+\f{(m+\f12)\pi}{\om}\in (-1,0)\quad\hbox{for some}\quad m\in \Z
\eeq
where the contact angle $\om\in (0,\pi/2)$.

For example in the case of $H^5$ decomposition ($K=3$), we find that when 
\[
m=-1, \,\om\in(\f\pi8,\f\pi6)\quad\hbox{or}\quad m=-2,\, \om\in(\f{3\pi}8,\f\pi2),
\] one has correspondingly 
\[
\lam_{-1}=-\f\pi{2\om}+3=\lam+3,\quad\lam_{-2}=-\f{3\pi}{2\om}+3=3\lam+3\in(-1,0)
\] 
Following the previous analysis, we have the $H^5$ decomposition for $v_c$:
\beno
v_c&=&v_{r3}+\chi_1(\om)c_{s1}(S_1+S_{1,r,1}+S_{1,r,2})+\chi_2c_{s2}(S_2+S_{2,r,1})\\
&&\quad+\chi_3(\om)c_{s3}S_3+\chi_4(\om)S_{4}
\eeno 
where $S_{i,r,j}$ ($i,j =1,2$) are  functions modifying the regularity of $\na\cdot\cP_s\na S_i$  from Lem \ref{S_{1,r,1}} and Remark \ref{general modify}, and we can tell from a similar computation as before that $S_3=S$ for different $\om$, $S_4$ the singular function corresponds to $\lam_{-2}$:
\[
S_4=\bar S_4\circ T^{-1}_0=\be(r)r^{-3\lam}\cos\big(3\lam(\tht+\om_2)\big)\circ T^{-1}_0.\] Moreover,  one defines the cut-off functions as follows
\[
\chi_3(w)=\left\{\begin{array}{ll}1,\quad \om\in (\f\pi 8,\f\pi 6)\\
0,\quad \hbox{otherwise}\end{array}\right.\quad\hbox{and}\quad \chi_4(\om)=\left\{\begin{array}{ll}1,\quad \om\in (\f{3\pi} 8,\f\pi 2)\\
0,\quad \hbox{otherwise.}\end{array}\right.
\] 
In a word, the higher-order decompositions can be done in a similar way, where the singular functions arise from the $H^2$ decomposition Thm \ref{H2 decompose} as well as the regularity-modifying Lem \ref{S_{1,r,1}}, meanwhile the singular coefficients can be estimated similarly as before. Based on the decomposition of $v_c$,  we can go back to the  solution $u$ for $\mbox{(MBVP)}$. 

To sum up, we present a general result for the $H^{2+K}$ decomposition and estimates for $\mbox{(MBVP)}$. 
\bthm{Proposition}\label{higher order decomposition u}
{\it Let $h\in H^K(\Om)$, $f\in H^{K+\f32}(\Ga_t)$, $g\in H^{K+\f12}(\Ga_b)$ and $\eta\in H^{K+2}(\R^+)$ for $K\in\N$. When $\om=\f\pi{2n}$ for some $n\in \N$, one requires additionally the compatible condition $L_{n-1}\in \tilde H^\f12(\Ga_b)$, where $L_{n-1}=L_{n-1}\big(\p^{\al}(u-u_b), \p^{n-2}(h-\Del u_b)\big)$ is a linear function of $\p^{\al}(u-u_b)$, $\p^{n-2}(h-\Del u_b)$  on $\Ga_b$ with $0<|\al|\le n-1$.

Then  one can have the $H^{2+K}$ decomposition for the solution $u$ from  $\mbox{(MBVP)}$
\[
u=u_{rK}+u_{sK}
\]
with $u_{rK}\in H^{2+K}$ the regular part and $u_{sK}\in H^2$ the singular part concentrates near the corner in the form
\[
u_{sK}=\sum_{1\le i\le m_K}\chi_i(\om)c_{si}\Big(\bar S_i+\sum_{1\le j\le n_i}\bar S_{i,r,j}\Big)\circ T_0^{-1}\circ (T_S)^{-1}
\]
where $\chi_i(\om)$ are the characteristic functions, $c_{si}$ are singular coefficients depending on  the right side of the system and $\bar S_i$ the singular functions decided by the left side of the system with the following formula
\beno
\bar S_i=\beta(r)r^{-\lam_i}\cos\big(\lam_i(\tht+\om_2)\big) \quad \hbox{on}\quad \cS_0,
\eeno
with $\be$  the cut-off function and $\lam_i=(m_i+\f12)\pi/\om$ with some $m_i\in \Z$. Here $\bar S_{i,r,j}=\be(r)r^{-\lam_i+j}a_{i,j}(\tht)$ are the regularity-modifying functions from Remark \ref{general modify} and the constants $m_K,n_i\in \N$ depend on $K,i$ respectively. Moreover, one has estimates both for the two parts as follows:
\[
\|u_{sK}\|_{2,2}+\|u_{rK}\|_{K+2,2}\le C\big(\|h\|_{K,2}+|f|_{H^{K+\f32}}+|g|_{H^{K+\f12}}\big)
\]
for some constant $C=C(\Ga_b,k,\lam_i,|\eta|_{H^{K+2}})$.
}\ethm
\begin{proof}
The proof can be done by an induction argument. In fact, we've showed the cases when $K=1,2$. For the higher order cases, the proof is similar by starting with taking $\na_{\tau_t}$.  Moreover, similarly as in the beginning of this subsection with $K=3$, one could verify the expression for $\lam_i$ by checking on $\lam_m\in (-1,0)$ there.

On the other hand, we need to notice the compatible condition when $\om=\f\pi{2n}$. One can tell from the compatible conditions before that, Since the singular functions in \eqref{comp H4 1} and \eqref{comp H4 2} are good, the condition lies on the derivatives for $v, h_c, g_c$ consequently.
\end{proof}

\bthm{Remark}\label{critical angle}{\it From the discussion above,  one can summarize on the range of $\om$ with respect to the order $K$  when no singularity appears. In fact, one finds from \eqref{range of lambda} that, for a fixed $K\in \N$, only when there exists some $m\in{\Z}^-$ satisfying
\[\om\in \Big(-\f{(2m+1)\pi}{2(K+1)},\,-\f{(2m+1)\pi}{2K}\Big)\subset (0,\,\f\pi 2),\] there will be a singular function corresponding to $m$.  As a result,  a lowerbound for $\om$ when the singularity takes place for some $K$ is $\om=\f\pi{2(K+1)}$ with $m=-1$, which means when 
\[\om\in (0,\, \f\pi{2(K+1)}),\] no singularity appears.
}\ethm

To close this section, we would like to consider again  the  special case when $\om=\f\pi{2n}$ for some $n\in \N$. In fact, we can derive the decomposition by taking the derivative $\na_{\tau_t}$  on the system \eqref{v1 problem} for $v_c$ on the corner domain $\cS$. when we take $\na^l_{\tau_t}$ on $v_c$ with $l\le n-1$, we find again mixed boundary system essentially for $\na^l_{\tau_t}v_c$,  with $b_l=-\tan(l\om)$ in the Neumann condition on $\Ga_b$. Applying Thm \ref{H2 decompose}, we need the eigenvalue 
\[
\lam_m=\f{l\om+(m+1/2)\pi}{\om}=l+2n(m+1/2)\in (-1,0)
\] for some $m\in \Z$, $n\in \N$, which turns out to be impossible. So there is no singular part in $\na^l_{\tau_t}v_c$, and consequently no singular part  for $v_c$.

when we take $\na^n_{\tau_t}$ on $v_c$, a Dirichlet-boundary system similar as \eqref{D system for  w2} appears with compatible conditions. As a result, there is no singular part in $\na^n_{\tau_t}v_c$ by Thm \ref{H2 decompose 1}, which implies that $v_c$ has no singular part.  For the case of the higher-order derivatives, the regularity follows similarly. As a summary, {\it we find no singular part in $v_c$ and the solution $u$ for $\mbox{(MBVP)}$ in $H^{2+K}$ decomposition  when the contact angle $\om=\f\pi{2n}$, i.e. we have $u\in H^{2+K}(\Om)$ with standard elliptic estimates, meanwhile, to be the price, we have $n-1$ compatible conditions. }

\medskip

\section{D-N operator estimates}
In this section, applying the theory of  elliptic decompositions and  estimates derived above, we are ready to investigate the properties of D-N operator. Recall the definition of the scaled D-N operator $G(\eta)$:
\[
G(\eta) f=\sqrt{1+|\eta'|^2}\p_{n_t}u|_{\Ga_t}.
\] with $u$ satisfying the elliptic system
\beq\label{DN elliptic system}
\left\{\begin{array}{ll}
\Delta u=0,\qquad \hbox{in}\quad \Omega\\
u|_{\Ga_t}=f,\qquad \p_{n_b} u|_{\Ga_b}=0.
\end{array}\right.
\eeq
One needs to notice that in this section,  the boundary function $f$ sometimes will be defined on $\R^+$ for convenience, while  in previous sections,  the boundary function $f$ ( and $g$) is always defined on curves $\Ga_t$ (and $\Ga_b$).  And it's well known that these two definitions are equivalent.
 
 In fact, it's well-known already that, D-N operator is an order-one elliptic operator, and it's global. When constructing energy estimates for the linearized water-wave problem, one will need  estimates for D-N operator as well as estimates for the shape derivative of D-N operator.  

We plan to give the elliptic estimate  for D-N operator first, which  is similar as in smooth-boundary case, except that there will be singular parts arising from the related elliptic problem.  

Next, we plan to derive the shape derivative  and the estimates, again with  singular parts.  In this part, as a  first step, we will derive a global expression for the shape derivative following \cite{ShZeng},  and the estimates with our elliptic theories follows. For the second step, we would like to show a specified form of the shape derivative under a special case, which tells us that the shaped derivative with a corner is similar as the shape derivative in \cite{Lannes}.

In this section, we assume $\eta\in H^{K+2}(\R^+)$ and $f\in H^{K+\f32}(\R^+)$ with some integer $K\geq 2$. Applying Prop \ref{higher order decomposition u} on \eqref{DN elliptic system}, we find  $u\in H^2(\Om)$ and  the decomposition
\beno
u=u_{rK}+ u_{sK},
\eeno
where the regular part $u_{rK}\in H^{K+2}(\Om)$ and the singular part $u_{sK} \in H^2(\Om)$.  Here, recall that $v_{sK}$ is supported in the neighborhood of corner.  Moreover, we have the estimate
\beq\label{v_rK estimate}
\|u_{rK}\|_{K+2,2}\leq C(|\eta|_{H^{K+2}})|f|_{H^{K+\f32}}.
\eeq

\subsection{Basic properties for Dirichlet-Neumann operator}
The proposition below shows that D-N operator is an order-1 operator.
\bthm{Proposition}\label{DN op estimate}
{\it If $\eta\in H^{K+2}(\R^+)$, then for all $f\in H^{K+\f32}(\R^+)$, we have that
\beno
\big|G(\eta)f-\sqrt{1+|\eta'|^2}\p_{n_t} u_{sK}|_{\Gamma_t}\big|_{H^{K+\f12}}\leq C( |\eta|_{H^{K+2}})|f|_{H^{K+\f32}}.
\eeno}
\ethm
\begin{proof}
Recalling the definition of D-N operator, we have that
\beno
G(\eta)f-\sqrt{1+|\eta'|^2}\p_{n_t} u_{sK}=\sqrt{1+|\eta'|^2}\p_{n_t} (u-u_{sK})|_{\Gamma_t}=\sqrt{1+|\eta'|^2}\p_{n_t} u_{rK}|_{\Gamma_t}.
\eeno
Applying Rmk \ref{one side Dirichlet  trace} and \eqref{v_rK estimate}, this proposition can be proved.
\end{proof}
\bthm{Remark}{\it 
 Note that when $K$ is fixed, we know from Prop \ref{higher order decomposition u} that the expression for $v_{sK}$ is a fixed linear combination of singular functions with constant coefficients depending on the right side of the elliptic system. As a result, the expression $\sqrt{1+|\eta'|^2}\p_{n_t} u_{sK}$ from  Prop \ref{DN op estimate} can be computed explicitly, which would be a linear combination of $\p^{\bar\cP_s}_{n_t} \big(\bar S_i+\sum_j \bar S_{i,r,j}\big)$ with singular coefficients $c_i$.
}\ethm
\medskip

Besides, we give some positive properties of D-N operator which can be found in \cite{Lannes}.  In \cite{Lannes}, the author considered the case  that $\eta$ is defined on $\R$, but the proof can be adjusted easily for our case where $\eta$ is defined on $\R^+$ only.
\bthm{Proposition}
{\it Let $\eta\in H^3(\R^+)$, then the following properties hold for any two functions $f,g\in H^\f12(\R^+)$:\\
1) D-N operator is self-adjoint:
\beno
(G(\eta)f, g)=(f, G(\eta)g),
\eeno
2) D-N operator is positive:
\beno
(G(\eta)f, f)\geq 0,
\eeno
3) The order-one elliptic estimate 
\beno
(G(\eta)f, g)\leq C|f|_{H^{\f12}}|g|_{H^{\f12}},
\eeno
and for any $\mu>\mu_0$ with some fixed $\mu_0$, we have that
\beno
(G(\eta)f+\mu f, f)\geq C_\mu |f|_{H^{\f12}}^2
\eeno where the constants $C,\,C_u$ depending on $|\eta|_{H^3}$.}
\ethm

\medskip

\subsection{Shape derivative of Dirichlet-Neumann operator } As already mentioned in the beginning of this section,  the D-N operator is a global one and the shape derivative is also taken globally with respect to the free surface deviation $\eta$. This seems to have some conflict with the local singularities  near the corner as well as our local straightening procedures. In fact, in our case it's  difficult to localize the related elliptic problem first and  find the `local' shape derivative applying the method in \cite{Lannes}. Instead, inspired by \cite{ShZeng},  we consider to find the shape derivative in a direct way without any straightening procedure. Based on the global expression, the estimates can be finished using the elliptic theory in previous sections.

\subsubsection{The global expression}
We will prove Thm \ref{DN op geometric} through this section. Following the idea in \cite{ShZeng}, for a given upper surface $z=\eta(x)$ with $\eta\in H^{K+2}(\R^+)$, we consider its variation  with a (time) parameter $s$ and a variational vector field ${\bf w}=(w_1,\,w_2)^t\in H^{K+2}(\Ga_t)$ on the surface, meanwhile a Lagrangian formulation is used in this part. 

Now, the new upper surface at time $s$ is denoted by 
\[
\Ga_{ts}=\{(x,z)|\,z=\eta_s(x),\,x\ge x_s\}
\] 
with  $x_s|_{s=0}=0$ and the bottom is  still $\Gamma_b$ with an extended left end:
\[
\Ga_b=\{(x,z)|\,z=l(x),\,x\ge x_s\}\quad\hbox{where}\quad l(x)=-\ga x,\,x_s\le x\le x_0.
\]  Here the contact point is $(x_s, -\gamma x_s)$. 

We denote by $X(s,x,\eta(x))$ the trajectory of any point $(x,\eta(x))$ on $\Ga_t$ at time $s$, which satisfies
\[
X(0,x,\eta(x))=\big(x,\,\eta(x)\big),\quad \p_sX(s,x,\eta(x))={\bf w}_s\big( x,\eta(x)\big),
\]  
where ${\bf w}_s|_{\Gamma_b}$ is assumed to be always tangent to $\Ga_b$ at the contact point. Moreover, ${\bf w}_s|_{s=0}=\bf w$. This means that  our contact point is assumed to move along $\Ga_b$ in the variation. The unit outward normal vector on $\Ga_{ts}$ is denoted by ${\bf n}_{ts}$. The domain corresponding to $\Ga_{ts}$ at time $s$ is denoted as $\Om_s$ with bottom $\Ga_b$.

The material derivative  is therefore defined as 
\[D_s=\p_s+\na_{{\bf w}_s}\] on the upper surface $\Ga_{ts}$. Moreover,  $D_s$ is tangent to  $\bigcup\limits_{s}\Ga_{ts}$.

Instead of  the scaled D-N operator listed before, we  prefer to use  the non-scaled D-N operator related to system \eqref{DN elliptic system}: 
\[
G(\eta)f=\p_{n_t} u|_{\Ga_t}=\na_{n_t} u|_{\Ga_t},
\] 
since we don't want to deal with the coefficient $\sqrt{1+|\eta'|^2}$ here. The notation is still the same as the scaled one for the sake of convenience. Besides, we choose to define $f$ on $\Ga_t$ here rather than  on $\R^+$.

Consequently, the variation of the D-N operator on $f$ with respect to the parameter $s$ is expressed by  
\[
D_s G(\eta_s)f_s|_{s=0}=D_s\big(\na_{n_{ts}} u_s|_{\Ga_{ts}}\big)\big|_{s=0}=D_s\na_{n_{ts}} u_s\big|_{\Ga_{ts},s=0},
\] while noting that $D_s$ is the material derivative on $\Ga_{ts}$, and  $f_s$, $u_s$ change along the material trajectory satisfying $f_s|_{s=0}=f$, $u_s|_{s=0}=u$.  

Here we will need to extend $\bf w$ on $\Ga_t$ to the whole domain $\Om$. Firstly, $\bf w$ is extended to the bottom $\Ga_b$: For each point $\big(x,\,l(x)\big)\in \Ga_b$, $\bf w$ is always tangent to $\Ga_b$ and satisfies 
\[{\bf w}\big(x,\,l(x)\big)=\left\{\begin{array}{ll}{\bf w}(0,0),\quad x\le x_0,\\
{\bf 0},\quad x\ge x_1
\end{array} 
\right.\] with some fixed small number $0<x_0<x_1$. This extension could be done since $\Ga_b$ becomes $z=-\ga x$ near $O(0,0)$. According to  Thm \ref{Dirichlet boundary trace}, we know that $\bf w$ could be extended into $\Om$ and also be controlled by the boundary
\beq\label{DN w estimate}
\|{\bf w}\|_{K+2,2}\le C(|\eta|_{H^{K+\f32}})|{\bf w}|_{H^{K+\f32}(\Ga_t)}.
\eeq Here we only use ${\bf w}\in H^{K+\f32}(\Ga_t)$ in stead of $H^{K+2}(\Ga_t)$, since we will see from the computations below that $H^{K+\f32}(\Ga_t)$ is enough already. The extension of ${\bf w}_s$ on $\Om_s$ is similar when it's needed, and $D_s$ is also tangent to $\bigcup\limits_{s}\Ga_b$.  Besides, one can see from Section 7 that only $|\eta|_{H^{K+\f32}}$  is needed instead of $|\eta|_{H^{K+2}}$.

Moreover,  although it seems that ${\bf n}_{t}$ (and ${\bf n}_b$) should be extended as well, one needs to notice that from the final expression of the shape derivative we only care about the original normal vectors on the boundary $\Ga_t$ (and $\Ga_b$) at $s=0$. Therefore, we don't consider explicitly the extension for ${\bf n}_t$ (and ${\bf n}_b$) here.

In order to find the shape derivative, we plan to expand $D_sG(\eta_s)f_s|_{s=0}$ into several  terms including $G(\eta)\big(D_sf_s|_{s=0}\big)$, therefore the shape derivative for $G(\eta)f$ is naturally given by the commutator
\[D_s G(\eta_s)f_s|_{s=0}-G(\eta)\big(D_sf_s|_{s=0}\big).\]
In fact, a direct computation shows that
\beq\label{Ds Gf 1}\begin{split}
D_sG(\eta_s)f_s|_{s=0}=&D_s\na_{n_{ts}} u_s\big|_{\Ga_{ts},s=0}=[D_s,\,\na_{n_{ts}}]u_s+\na_{n_{ts}}D_su_s\big|_{\Ga_{ts},s=0}\\
=& (D_s{\bf n}_{ts})|_{s=0}\cdot \na u-\na_{n_t}{\bf w}\cdot \na u+\na_{n_t}(D_su_s)|_{s=0}\big|_{\Ga_t},
\end{split}\eeq where one needs to deal with two terms $D_s{\bf n}_{ts}|_{s=0}$ and $\na_{n_t}(D_s u_s)|_{s=0}$. 

Firstly, remembering that ${\bf n}_{ts}$ is the unit outward normal vector on $\Ga_{ts}$, a simple computation  from \cite{ShZeng} tells us that
\beq\label{Ds nt}
D_s{\bf n}_{ts}|_{s=0}=-\big((\na {\bf w})^t{\bf n}_t\big)^T
\eeq where $^T$ means to take the tangential part on $\Ga_t$.

Secondly, we consider the elliptic system for $D_s u_s$ on $\Om_s$ (since $\Om$ also changes for time $s$). According to system \eqref{DN elliptic system} for $u$, $u_s$ satisfies the following system
\[
\left\{\begin{array}{ll}
\Delta u_s=0,\qquad \hbox{in}\quad \Om_s\\
u_s|_{\Ga_{ts}}=f_s,\qquad \p_{n_b} u_s|_{\Ga_b}=0.
\end{array}\right.
\]
So  one can see that $\Del D_su_s$ satisfies
\[
\Del D_s u_s=D_s\Del u_s+[\Del,\,D_s]u_s=2\na {\bf w}_s\cdot \na^2u_s+\Del{\bf w}_s\cdot\na u_s
\] on $\Om_s$,  where one uses the notation $\na {\bf w}_s\cdot \na^2u=\sum\limits_{i,j=1,2}\p_i  w_{sj}\p_i\p_j u$ with $\p_1=\p_x$, $\p_2=\p_z$. Besides, referring to   boundary conditions of $D_su_s$, one finds that $D_s u_s|_{\Ga_{ts}}=D_s f_s$ and 
\[\begin{split}
\p_{n_b}D_su_s|_{\Ga_b}=&[\na_{n_b},\,D_s]u_s+D_s\na_{n_b}u_s|_{\Ga_b}
=u_s\,\na_{n_b}{\bf w}_s-D_s{\bf n}_b\cdot\na u_s|_{\Ga_b}\\
=& u_s\,\na_{n_b}{\bf w}_s+\na_{{\na u_s}^T}{\bf w}_s\cdot {\bf n_b}|_{\Ga_b},
\end{split}\] where we used the boundary condition $\na_{n_b}u_s|_{\Ga_b}=0$ and the fact $D_s{\bf n}_b=-\big((\na {\bf w}_s)^t{\bf n}_b\big)^T$, which is similar as $D_s{\bf n}_{ts}$.

Consequently, we find that $D_s u_s|_{s=0}$ satisfies the elliptic system 
\[\left\{\begin{array}{ll}
\Del (D_s u_s)|_{s=0}=2\na {\bf w}\cdot \na^2u+\Del{\bf w}\cdot\na u,\quad \hbox{on}\quad \Om,\\
(D_s u_s)|_{s=0}\big|_{\Ga_{t}}=D_s f_s|_{s=0},\quad \p_{n_b}(D_su_s)|_{s=0}\big|_{\Ga_b}= u\,\na_{n_b}{\bf w}+\na_{(\na u)^T}{\bf w}\cdot {\bf n_b}|_{\Ga_b}.
\end{array}\right.\]
In order to retrieve $G(\eta)\big(D_sf_s|_{s=0}\big)$ from \eqref{Ds Gf 1}, we decompose $D_s u_s|_{s=0}$ into two parts as follows
\[D_s u_s|_{s=0}=u_1+u_2\] with $u_1$ satisfing
\beq\label{u1 elliptic system}\left\{\begin{array}{ll}
\Del u_1=0,\quad \hbox{on}\quad \Om,\\
u_1|_{\Ga_{t}}=D_s f_s|_{s=0},\quad \p_{n_b}u_1|_{\Ga_b}=0
\end{array}\right.\eeq
and $u_2$ satisfing
\beq\label{u2 elliptic system}\left\{\begin{array}{ll}
\Del u_2=2\na{\bf w}\cdot\na^2 u+\Del{\bf w}\cdot \na u,\quad \hbox{on}\quad \Om,\\
u_2|_{\Ga_{t}}=0,\quad \p_{n_b}u_2|_{\Ga_b}=u\,\na_{n_b}{\bf w}+\na_{(\na u)^T}{\bf w}\cdot {\bf n_b}|_{\Ga_b}.
\end{array}\right.\eeq
 As a result, we find 
\beq\label{Ds u term}
\na_{n_t}(D_su_s)|_{s=0}\big|_{\Ga_t}=\na_{n_t}u_1+\na_{n_t}u_2|_{\Ga_t}=G(\eta)\big(D_s f_s|_{s=0}\big)+\na_{n_t}u_2|_{\Ga_t}.
\eeq Plugging \eqref{Ds nt}, \eqref{Ds u term} into \eqref{Ds Gf 1}, we finally arrive at the expression
\[ 
D_s G(\eta_s)f_s|_{s=0}=G(\eta) \big(D_sf_s|_{s=0}\big)+\big(\na_{n_t} u_2-\na_{n_t}{\bf w}\cdot\na u-\na_{(\na u)^T}{\bf w}\cdot{\bf n}_t\big)\big|_{\Ga_t},
\] from which we conclude that the shape derivative of $G(\eta)f$ is 
\beq\label{shape deri  geometric}
D_sG(\eta_s)f_s|_{s=0}-G(\eta)\big(D_sf_s|_{s=0}\big)=\na_{n_t} u_2-\na_{n_t}{\bf w}\cdot\na u-\na_{(\na u)^T}{\bf w}\cdot{\bf n}_t\big|_{\Ga_t}.
\eeq 

Based on the expression above, We are ready to have some estimates for  the shape derivative concerning singularities near the corner. 

To begin with, recalling that in the beginning of Section 6,  we have the decomposition for $u$ satisfying system \eqref{DN elliptic system} with $f\in H^{K+\f32}(\Ga_t)$:
\[u=u_{rK}+u_{sK},\] where $u_{rK}\in H^{K+2}(\Om)$ and $u_{sK}\in H^2(\Om)$ are the regular and singular part respectively. Moreover,  the estimate holds
\beq\label{DN u estimate}
\|u_{rK}\|_{K+2,2}+\|u_{sK}\|_{2,2}\le C\big(|\eta|_{H^{K+2}}\big)|f|_{H^{K+\f32}}.
\eeq Based on the decomposition of $u$, we also need to consider about the decomposition of $u_2$.  In fact, recalling system \eqref{u2 elliptic system} for $u_2$, let 
\[u_2=\td u_{2rK}+\td u_{2sK},\] where $\td u_{2rK}$ satisfies
\beq\label{DN u2rK  system}\left\{\begin{array}{ll}
\Del \td u_{2rK}=h_1,\quad \hbox{on}\quad \Om,\\
\td u_{2rK}|_{\Ga_t}=0,\quad \p_{n_b}\td u_{2rK}|_{\Ga_b}=g_1
\end{array}\right.\eeq with
\[h_1=2\na{\bf w}\cdot\na^2 u_{rK}+\Del{\bf w}\cdot \na u_{rK},\quad g_1=u_{rK}\,\na_{n_b}{\bf w}+\na_{(\na u_{rK})^T}{\bf w}\cdot {\bf n_b}|_{\Ga_b},\]
and the remainder $\td u_{2sK}$ satisfies
\beq\label{DN u2sK  system}\left\{\begin{array}{ll}
\Del \td u_{2sK}=h_2,\quad \hbox{on}\quad \Om,\\
\td u_{2sK}|_{\Ga_t}=0,\quad \p_{n_b}\td u_{2sK}|_{\Ga_b}=g_2
\end{array}\right.\eeq with
\[h_2=2\na{\bf w}\cdot\na^2 u_{sK}+\Del{\bf w}\cdot \na u_{sK},\quad g_2=u_{sK}\,\na_{n_b}{\bf w}+\na_{(\na u_{sK})^T}{\bf w}\cdot {\bf n_b}|_{\Ga_b}.\] Since we have ${\bf w} \in H^{K+\f32}(\Ga_t)$, we could tell directly that
\[h_1\in H^K(\Om),\quad g_1\in H^{K+\f12}(\Ga_b), \quad h_2\in L^2(\Om),\quad g_2\in H^{\f12}(\Ga_b),\] and the estimates follows 
\beq\label{DN right side estimate}\begin{split}
& \|h_1\|_{K,2},\, |g_1|_{H^{K+\f12}}\le C\|u_{rK}\|_{K+2,2}\|{\bf w}\|_{K+2,2},\\
& \|h_2\|_{2},\, |g_2|_{H^{\f12}}\le C\|u_{sK}\|_{2,2}\|{\bf w}\|_{K+2,2}
\end{split}\eeq with the constant $C$ depending on $\Ga_b$.

Applying Prop \ref{higher order decomposition u} on both $\td u_{2rK}$ and $\td u_{2sK}$ , one gets  for $\td u_{2rK}$ that there also exists a decomposition 
\[
\td u_{2rK}=v_{rK}+v_{sK}
\] such that $v_{rK}\in H^{K+2}(\Om)$, $v_{sK}\in H^2(\Om)$ satisfy
\[
\|v_{rK}\|_{K+2,2}+\|v_{sK}\|_{2,2}\le C\big(|\eta|_{H^{K+2}}\big)\big(\|h_1\|_{K,2}+|g_1|_{H^{K+\f12}}\big).
\] Besides, one can also see that $\td u_{2sK}\in H^2(\Om)$ while noticing that there is no more singularity here, and the estimate follows:
\[
\|\td u_{2sK}\|_{2,2}\le C\big(|\eta|_{H^2}\big)\big(\|h_2\|_2+|g_2|_{H^\f12}\big).
\] Combining these two estimates above together with \eqref{DN w estimate}, \eqref{DN u estimate} and \eqref{DN right side estimate}, we arrive at the decomposition and the estimate for $u_2$ as 
\[u_2=u_{2rK}+u_{2sK},\]  where $u_{2rK}=v_{rK}\in H^{K+2}(\Om)$, $u_{2sK}=\td u_{2sK}+v_{sK}\in H^2(\Om)$ and satisfy the following estimate
\beq\label{DN u2 estimate}
\|u_{2rK}\|_{K+2,2}+\|u_{2sK}\|_{2,2}\le C\big(|\eta|_{H^{K+2}}\big)|f|_{H^{k+\f32}}|{\bf w}|_{H^{K+\f32}}
\eeq

Consequently, going back to the shape derivative \eqref{shape deri  geometric} and applying the decompositions for $u$ and $u_2$, some more direct computations and Rmk \ref{one side Dirichlet trace} lead to our  Thm \ref{DN op geometric} in the beginning of our paper.

\subsubsection{A special case} We prove Thm \ref{thm:D-N-S} in this section.  Here we  consider the domain $\Om$ with a straight bottom $z=-\ga x$ and assume that $\bar\eta=\eta+\ga x$ is invertible globally, which means that the local transformation $T_S$ becomes a global one.  

Under this assumption, we show that  the shape derivative is in a similar form as in \cite{Lannes}.  The result can be taken as a special and `local' expression for the shape derivative near the corner, and the computations in this part are only formal computations, since we don't have any existence theorem as Lemma \ref{variation solution} for this infinite-depth case. 

The plan of this section is as follows: Firstly, similar as \cite{Lannes}, we  straighten the domain $\Om$ into  $\cS$ and  derive the corresponding shape derivative on $\cS$.  Notice that in this section, the  function $f$ on the upper surface $\Ga_t$  is defined as $f=f(x)$ with $x\in \R^+$, which changes as we change the domain $\Om$ into $\cS$. Secondly, we go back to find the shape derivative on $\Om$, which is a little different from the one on $\cS$. 

Through out this part, we denote by $(\bar x,\,\bar z)$ the point in $\Om$, while we use $(x,z)$ for the point in $\cS$.

To begin with,  we define as before that \[v=u\circ T_S\quad\hbox{on}\quad \cS=\{(x,z) |0\leq z\leq kx\}\] and focus on the new system equivalent with system \eqref{DN elliptic system}:
\beq\label{equ:v_eta}
\left\{\begin{array}{ll}
\na\cdot\cP_{S}\na v=0\qquad\hbox{on}\quad \cS\\
v|_{\Gamma_t}=f_s \triangleq f\circ T_S,\qquad \p^{\cP_{S}}_{n_b}v|_{\Gamma_b}=0
\end{array}\right.
\eeq
where recall from Section 5.1 that
\[
\cP_{S} =\left(\begin{matrix}(1+\ga d)^2+d^2 & \ga+(1+\ga^2)d\\
\ga+(1+\ga^2)d & 1+\ga^2\end{matrix}\right)\triangleq\left(\begin{matrix}p_1 & p_2\\
p_2 & p_4 \end{matrix}
\right),
\] 
with  
\[
d=d(z)=\f 1k-\f1{\ga+\eta'( \bar\eta^{-1}(z))},\quad \bar \eta(\bar x)=\ga \bar x+\eta(\bar x)
\] with some constant $k$. In fact, one can fix the constant $k$ here as, for example $k=\ga+\eta'(0)$ or simply $k=1$.
 
We choose to use the scaled D-N operator here. A direct computation shows that  the corresponding scaled D-N operator for system \eqref{equ:v_eta} is
\beq\label{relation of two DN op}
\begin{split}
G(\eta)f_s=&\sqrt{1+k^2}\p^{\cP_S}_{n_t} v|_{\Gamma_t}\\
=&\Big(\f{k}{\ga+\eta'}\sqrt{1+|\eta'|^2}\p_{n_t} u|_{\Ga_t}\Big)\circ T_S=\Big(\f{k}{\ga+\eta'}G(\eta) f\Big)\circ T_S
\end{split}
\eeq as  long as one notices that the relationship  between $(x,z)\in \cS$ and $(\bar x,\bar z)\in \Om$ on the upper boundaries as follows:
\[\bar x=\bar\eta^{-1}(kx).\]

In the following text, we denote by $d_\eta g$  the Fr\'echet derivative at $\eta$ for some function $g(\eta)$ as \[d_\eta g\cdot h=\f d{d\eps}g(\eta+\eps h)\big|_{\eps=0},\] where $h$ would be in the same space as $\eta$.
The Fr\'echet derivative of $v$ is denoted by
\[ v_{\eta}\triangleq d_{\eta}v\cdot h,\]
which  solves
\beq\label{equ:v_a}
\left\{\begin{array}{ll}
\na\cdot\cP_{S}\na v_{\eta}=-\na\cdot d_{\eta}\cP_{S}\cdot h\na v\qquad\hbox{on}\quad \cS\\
v_{\eta}|_{\Gamma_t}=d_\eta f_s\cdot h,\qquad \p^{\cP_{S}}_{n_b}v_{\eta} |_{\Gamma_b}=-(0,1)^t \cdot(d_{\eta}\cP_{S}\cdot h\na_{x,z} v )|_{\Gamma_b}.
\end{array}\right.
\eeq where ${\bf n}_b=-(0,1)^t$ since $\Ga_b=\{ z=0\}$ for $\cS$.

Next, motivated by \cite{Lannes}, we give an explicit function solving \eqref{equ:v_a} except for the Dirichlet condition at the surface. 
\bthm{Lemma}
The function $v_{\eta}^1=B\p_x v$ solves
\beq
\left\{\begin{array}{ll}
\na\cdot\cP_{S}\na v_{\eta}^1=-\na\cdot d_{\eta}\cP_{S}\cdot h\na v\qquad\hbox{on}\quad \cS\\
v_{\eta}^1|_{\Gamma_t}=B\p_xv|_{\Gamma_t},\qquad \p^{\cP_{S}}_{n_b}v_{\eta}^1 |_{\Gamma_b}=-(0,1)^t \cdot(d_{\eta}\cP_{S}\cdot h\na_{x,z} v )|_{\Gamma_b},
\end{array}\right.
\eeq
with the quantity $B=B(z)=-\int^z_0 (d_{\eta} d\cdot h)(s) ds$.
\ethm
\begin{proof}
Here we notice that $\cP_{S}$ and $B$ depend  only on $z$. Combining the equation for $v$ from \eqref{equ:v_eta}, a direct computation leads to
\beno
\na\cdot\cP_{S}\na v_{\eta}^1=\na\cdot Q\na v,
\eeno
where
\beno
Q=\left(
  \begin{array}{cc}
   2p_2 B' & p_4B' \\
    p_4B' & 0 \\
  \end{array}
\right).
\eeno

On the other hand, we have that
\beno
d_{\eta} \cP_{S}\cdot h=\left(
  \begin{array}{cc}
   2(\gamma+(1+\gamma^2)d)d_{\eta} d\cdot h & (1+\gamma^2)d_{\eta} d\cdot h \\
    (1+\gamma^2)d_{\eta} d \cdot h& 0 \\
  \end{array}
\right).
\eeno

It is straightforward to check that
\beno
\na\cdot\cP_{S}\na v_{\eta}^1=-\na\cdot(d_{\eta} \cP_{S}\cdot h\na v).
\eeno

Next, we verify the boundary condition. Noticing that $B|_{\Gamma_b}=0$ and combining the boundary condition from \eqref{equ:v_eta}, we find that $v_{\eta}^1$ satisfies the boundary condition 
\beno
\p^{\cP_{S}}_{n_b}v_{\eta}^1 |_{\Gamma_b}=(p_2\p_x+p_4\p_z)v_{\eta}^1 |_{\Gamma_b}=p_4B'\p_xv |_{z=0}=-p_4(d_{\eta} d\cdot h) \p_xv |_{\Gamma_b}.
\eeno
Meanwhile, one has
\beno
(0,1)^t \cdot(d_{\eta}\cP_{S}\cdot h\na v_{ c} )|_{\Gamma_b}=p_4(d_{\eta} d\cdot h) \p_xv_{c} |_{\Gamma_b},
\eeno
which implies that
\beno
\p^{\cP_{S}}_{n_b}v_{\eta}^1 |_{\Gamma_b}=-(0,1)^t \cdot(d_{\eta}\cP_{S}\cdot h\na v_{ c} )|_{\Gamma_b}.
\eeno
In the end, the Dirichlet boundary condition on $\Ga_t$ is easy to verify by the definition of $v_{\eta}^1$, therefore the proof is done.
\end{proof}

From the lemma above, we arrive at the system for $v_{\eta}- v_{\eta}^1$:
\[
\left\{\begin{array}{ll}
\na\cdot\cP_{S}\na(v_{\eta}- v_{\eta}^1)=0\qquad\hbox{on}\quad \cS,\\
(v_{\eta}- v_{\eta}^1)|_{\Gamma_t}=d_\eta f_s\cdot h-B\p_xv|_{\Gamma_t},\qquad \p^{\cP_{S}}_{n_b}(v_{\eta}- v_{\eta}^1) |_{\Gamma_b}=0.
\end{array}\right.
\]
For the moment, we are ready to consider the shape derivative for D-N operator. We will firstly derive the shape derivative for $G(\eta)f_s$, and then we go back to $G(\eta)f$ applying \eqref{relation of two DN op}. 

In fact, from the definition of  D-N operator, we can compute the shape derivative with respect to $\eta$ term by term while noticing that ${\bf n}_t=(1+k^2)^{-\f12}(-k,1)^t$ here:
\beq\label{shape derivative decomp}
\begin{split}
(1+k^2)^{-\f12}d_{\eta}  G({\eta})f_s\cdot h=&d_{\eta}\big(\p^{\cP_S}_{n_t}v\big)\cdot h\\
=&\p_{n_t}^{d_{\eta}\cP_{S}\cdot h}v +\p_{n_t}^{\cP_{S}}v_{\eta}\big|_{\Gamma_t}\\
=&\p_{n_t}^{d_{\eta}\cP_{S}\cdot h}v+\p_{n_t}^{\cP_{S}}v_{\eta}^1+  \p_{n_t}^{\cP_c}(v_{\eta} -v_{\eta}^1)\big|_{\Gamma_t}\\
=&\big(\p_{n_t}^{d_{\eta}\cP_{S}\cdot h}v+\p_{n_t}^{\cP_{S}}v_{\eta }^1\big)\big|_{\Gamma_t}+(1+k^2)^{-\f12}G(\eta)(d_\eta f_s\cdot h-B\p_xv |_{\Gamma_t})\\
\triangleq& I+II,
\end{split}\eeq 
where the last term is derived by the system of $v_{\eta} -v_{\eta}^1$.

Due to the presence of the corner, the last term in the above equation should  have a singular decomposition, which can be done by decomposing $v_{\eta}- v_{\eta}^1$ properly. Since our computations here are only formal ones, we omit the details for singular decompositions. 

For the moment,  we plan to use  $B, \cP_c, v$ to compute the first two terms of the right hand of \eqref{shape derivative decomp}. Remembering the expressions for $d_{\eta}\cP_{S}\cdot h$ and $v_{\eta}^1$, we can prove the following lemma.
\bthm{Lemma}
{\it One has
\[
I= F_1\p_{\tau_t}^2 v+F_2\p_{\tau_t}\p_{n_t}  v+ F_3 \p_{\tau_t}v+ F_4 \p_{n_t}v \big|_{\Gamma_t}
\]
with coefficients $F_i$($i=1,2,3,4$) defined by
\ben
F_1&=& (1+k^2)^{-\f12}B(p_2+kp_4),\label{F_1}\\
  F_2&=&  (1+k^2)^{-\f12}B(kp_2-p_4),\label{F_2}\\
  F_3&=& -(1+k^2)^{-1}k\big(B (p_2+kp_4)\big)',\label{F_3}\\
   F_4&=&-(1+k^2)^{-1} k\big( B(kp_2-p_4)\big)'\label{F_4}
\een where the notation $'$ means taking $\p_z$.
}
\ethm
\begin{proof}
 First of all, let's recall from the previous lemma that the expression of $d_{\eta}  \cP_{S}\cdot h$ is 
\beno
d_{\eta} \cP_{S} \cdot h=-\left(
  \begin{array}{cc}
   2p_2B' & p_4B' \\
    p_4B' & 0 
  \end{array}
\right),
\eeno
which implies that
\beno
(1+k^2)^{\f12}\p_{n_t}^{d_{\eta} \cP_{S} \cdot h}v\big|_{\Gamma_t}=B'\big((2kp_2-p_4)\p_x v+kp_4\p_z v\big)\big|_{\Gamma_t}.
\eeno
A direct computation shows that
\beq\label{v_1 on Ga_t}\begin{split}
(1+k^2)^{\f12}\p_{n_t}^{\cP_{S} }v_{\eta}^1\big|_{\Gamma_t}=&B(p_2-kp_1)\p_x^2v +B(p_4-kp_2)\p_x\p_zv\\
&\quad +B'(p_4-kp_2)\p_xv\big|_{\Gamma_t}.
\end{split}\eeq As a result we find
\beq\label{the sum}
(1+k^2)^{\f12} I=B(p_2-kp_1)\p^2_xv+B(p_4-kp_2)\p_x\p_zv+B'kp_2\p_xv+B'kp_4\p_zv
\eeq

On the other hand, one could express $\p_x,\p_z$ by $\p_{\tau_t},\p_{n_t}$ on $\cS$ as
\beq\label{xz to tb}
\p_x=(1+k^2)^{-\f12}(-\p_{\tau_t}-k\p_{n_t}),\qquad \p_z=(1+k^2)^{-\f12}(-k\p_{\tau_t}+\p_{n_t})
\eeq where ${\bf \tau}_t=(1+k^2)^{-\f12}(-1,-k)^t$ and ${\bf n}_t=(1+k^2)^{-\f12}(-k,1)^t$.
Consequently, rewriting the elliptic equation of $v$ in \eqref{equ:v_eta}, we find that $v$ satisfies 
\[\begin{split}
&(p_1+2kp_2+k^2p_4)\p_{\tau_t}^2 v+2\big(kp_1+(k^2-1)p_2-kp_4\big)\p_{n_t}\p_{\tau_t}v\\
&\qquad+(k^2 p_1-2kp_2+p_4)\p_{n_t}^2 v
-(1+k^2)^{\f12} p'_2(\p_{\tau_t}+k\p_{n_t})v=0
\end{split}\]
with $k^2 p_1-2kp_2+p_4>0$ to be verified later. Restricting the  equation above on the boundary, we have that
\beq\label{eqn for v_nn}\begin{split}
&-(k^2 p_1-2kp_2+p_4)\p_{n_t}^2 v|_{\Gamma_t}\\
&=(p_1+2kp_2+p_4k^2)\p_{\tau_t}^2 v+2\big(kp_1+(k^2-1)p_2-kp_4\big)\p_{\tau_t}\p_{n_t}v\\
& \qquad-(1+k^2)^{\f12} p'_2(\p_{\tau_t}+k\p_{n_t})v\big|_{\Gamma_t}.
\end{split}\eeq

Going back to equation \eqref{the sum} and using \eqref{xz to tb} again to obtain that
\beno
 &&(1+k^2)^{\f32}I\\
&&= -Bk(k^2p_1-2kp_2+p_4)\p^2_{n_t}v-B\big(2k^2p_1+k(k^2-3)p_2-(k^2-1)p_4\big)\p_{\tau_t}\p_{n_t}  v\\
&&\quad-B\big(kp_1+(k^2-1)p_2-kp_4\big)\p_{\tau_t}^2 v-B'k(1+k^2)^{\f12}(p_2+kp_4)\p_{\tau_t}v\\
&&\quad-B'k(1+k^2)^{\f12}(kp_2-p_4)\p_{n_t}v \big|_{\Gamma_t}.
\eeno
As a result, combining \eqref{eqn for v_nn}, some computations lead to the identity
\[
(1+k^2)^{\f32}I
=a_1\p_{\tau_t}^2 v+a_2\p_{\tau_t}\p_{n_t}  v+a_3\p_{\tau_t}v+a_4\p_{n_t}v \big|_{\Gamma_t},
\]
where
\beno
&&a_1=(1+k^2)B(p_2+kp_4),\quad
a_2=(1+k^2)B(kp_2-p_4),\\
 && a_3= -k(1+k^2)^{\f12}\big(B(p_2+kp_4)\big) ',\quad   a_4=-k(1+k^2)^{\f12}\big(B(kp_2-p_4)\big)'.
\eeno
As a result,  we finally arrive at the conclusion of this lemma.
\end{proof}
We introduce the notation 
\[
Z_s=\p_{n_t}v|_{\Ga_t}.
\] Based on the lemma above and noticing that 
\beq\label{partial v to f}
\p_{\tau_t}v|_{\Ga_t}=-(1+k^2)^{-\f12}f'_s,\quad \p_{\tau_t}\p_{n_t}v|_{\Ga_t}=-(1+k^2)^{-\f12}Z'_s,
\eeq one could take a further step to rewrite $I$ as 
\[\begin{split}
I=& F_1\p_{\tau_t}^2 v+F_2\p_{\tau_t}\p_{n_t}  v+ F_3 \p_{\tau_t}v+ F_4 \p_{n_t}v \big|_{\Gamma_t:\,z=kx}\\
=&  \Big((1+k^2)^{-1}F_1(k\cdot)f''_s-(1+k^2)^{-\f12}F_3(k\cdot)f'_s\Big)+\Big(-(1+k^2)^{-\f12}F_2(k\cdot)Z'_s+F_4(k\cdot)Z_s\Big)\\
=& (1+k^2)^{-\f32}\Big(B(k\cdot)\big(p_2(k\cdot)+kp_4\big)f'_s\Big)'-(1+k^2)^{-1}\Big(B(k\cdot)\big(kp_2(k\cdot)-p_4\big)Z_s\Big)'.
\end{split}\]
 Moreover, one can compute directly from $G(\eta)f_s=(1+k^2)^{\f12}\p^{\cP_S}_{n_t} v|_{\Ga_t}$ that 
\beq\label{Z_c}
Z_s=(k^2p_1-2kp_2+p_4)|_{\Ga_t}^{-1}\Big((1+k^2)^{\f12}G(\eta)f_s+(1+k^2)^{-\f12}\big(kp_1+(k^2-1)p_2-kp_4\big)f'_s\Big).
\eeq
On the other hand, remembering that 
\[f_s(x)=f\circ T_S(x)=f(\bar x)\big|_{\bar x=\bar \eta^{-1}(kx)},\] we have 
\[
d_\eta f_s\cdot h(x)=\f{d}{d\eps}f\big(\bar\eta^{-1}_\eps(kx)\big)\big|_{\eps=0}
=\f{f' h}{\bar \eta'}\big|_{\bar x=\bar \eta^{-1}(kx)}
\] with $\bar\eta_\eps(\bar x)=\ga \bar x+\eta(\bar x)+\eps h(\bar x)$.
Combining this together with \eqref{xz to tb} and \eqref{partial v to f}, we can rewrite the last term $II$ of \eqref{shape derivative decomp} as   
\[\begin{split}
II=&(1+k^2)^{-\f12}G(\eta)\Big(\f{f' h}{\bar \eta'}\big(\bar \eta^{-1}(k\cdot)\big)
+(1+k^2)^{-\f12}\big(B\p_{\tau_t}v+kB\p_{n_t}v\big)|_{\Ga_t}\Big)\\
=& (1+k^2)^{-\f12}G(\eta)\Big(\f{f' h}{\bar \eta'}\big(\bar \eta^{-1}(k\cdot)\big)\Big)-(1+k^2)^{-\f32}G(\eta)\big(B(k\cdot)f'_s)\\
&\quad +k(1+k^2)^{-1}G(\eta)\big(B(k\cdot)Z_s\big),
\end{split}\]

Finally, summing up all the computations above, we arrive at the shape derivative of $G(\eta)f_s$ on $\cS$ as follows
\beq\label{shape seri on S}\begin{split}
&d_\eta G(\eta)f_s\cdot h=(1+k^2)^{\f12}(I+II)\\
=&G(\eta)\Big(\f{f' h}{\bar \eta'}\big(\bar \eta^{-1}(k\cdot)\big)\Big)-(1+k^2)^{-1}G(\eta)\big(B(k\cdot)f'_s)\\
&\ +k(1+k^2)^{-\f12}G(\eta)\big(B(k\cdot)Z_s\big)+(1+k^2)^{-1}\Big(B(k\cdot)\big(p_2(k\cdot)+kp_4\big)f'_s\Big)'\\
&\ -(1+k^2)^{-\f12}\Big(B(k\cdot)\big(kp_2(k\cdot)-p_4\big)Z_s\Big)'.
\end{split}\eeq
Now it's the time go back to find the shape derivative of $G(\eta)f$ on $\Om$. Noticing the relationship \eqref{relation of two DN op}  between $G(\eta)f_s$ and $G(\eta)f$, we can have by a direct computation that 
\[\begin{split}
&d_\eta G(\eta)f\cdot h(\bar x)\\
=& \f{d}{d\eps}\Big(\f{\ga+\eta'(\bar x)+\eps h'(\bar x)}{k} G(\eta+\eps h)f_{s,\eps}\big|_{x=\f {\ga \bar x+\eta(\bar x)+\eps h(\bar x)}{k}}\Big)\big|_{\eps=0}\\
=& d_\eta G(\eta)f_s\cdot h\big|_{x=\f{\ga \bar x+\eta(\bar x)}{k}}\f{\ga+\eta'(\bar x)}{k}+\f d{dx} G(\eta)f_s\big|_{x=\f{\ga \bar x+\eta(\bar x)}{k}}\,\f{h'(\bar x)}{k}\ \f{\ga+\eta'(\bar x)}{k}\\
&\ +G(\eta)f_s\big|_{x=\f{\ga \bar x+\eta(\bar x)}{k}}\,\f{h(\bar x)}{k}\\
=&d_\eta G(\eta)f_s\cdot h\big|_{x=\f{\ga \bar x+\eta(\bar x)}{k}}\f{\ga+\eta'(\bar x)}{k}
+\f d{d\bar x}\Big(G(\eta)f_s\big|_{x=\f{\ga \bar x+\eta(\bar x)}{k}}\,\f{h(\bar x)}{k}\Big)
\end{split}\] where $f_{s,\eps}$ means to use $\eta+\eps h$ instead of $\eta$ in $f_s$.
Plugging \eqref{shape seri on S} into the shape derivative above, we find that
\[\begin{split}
&d_\eta G(\eta)f\cdot h(\bar x)\\
=& G(\eta)\Big(\f{f' h}{\bar \eta'}\big(\bar \eta^{-1}(k\cdot)\big)\Big)\big|_{x=\f{\ga \bar x+\eta(\bar x)}{k}}\,\f{\ga+\eta'(\bar x)}{k}\\
&\ -(1+k^2)^{-1}G(\eta)\Big(\big(B(k\cdot)f'_s)-k(1+k^2)^{\f12}\big(B(k\cdot)Z_s\big)\Big)\big|_{x=\f{\ga \bar x+\eta(\bar x)}{k}}\,\f{\ga+\eta'(\bar x)}{k}\\
&\ +(1+k^2)^{-1}\f d{d\bar x}\Big(B(k\cdot)\big(p_2(k\cdot)+kp_4\big)f'_s
 -(1+k^2)^{\f12}B(k\cdot)\big(kp_2(k\cdot)-p_4\big)Z_s\big|_{x=\f{\ga \bar x+\eta(\bar x)}{k}}\Big) \\
&\quad + \f d{d\bar x}\Big(G(\eta)f_s\big|_{x=\f{\ga \bar x+\eta(\bar x)}{k}}\,\f{h(\bar x)}k\Big)
\end{split}\] while recalling that $\bar \eta(\bar x)=\ga\bar x+\eta(\bar x)$.
Since \eqref{relation of two DN op} can be rewritten as 
\[
G(\eta)f(\bar x)=G(\eta)f_s\big|_{x=\f{\ga \bar x+\eta(\bar x)}{k}}\,\f{\ga+\eta'(\bar x)}k,
\] in stead of $f_s$, we could express the shape derivative in terms of $f$: 
\[\begin{split}
&d_\eta G(\eta)f\cdot h(\bar x)\\
=& G(\eta)\Big(\f{h}{\ga+\eta'}f'\Big)-(1+k^2)^{-1} G(\eta)\Big(\f k{\ga+\eta'}\bar B f'-k(1+k^2)^{\f12}\bar B\bar Z_s\Big)\\
&\ +(1+k^2)^{-1}\f d{d\bar x}\Big( \f k{\ga+\eta'}\bar B(\bar p_2+k p_4)f' -(1+k^2)^{\f12}\bar B(k\bar p_2-p_4)\bar Z_s\\
&\quad +(1+k^2)\f{h}{\ga+\eta'}G(\eta)f\Big),
\end{split}\] where we denote
\[\begin{split}
&\bar p_1=p_1\big(\bar\eta(\bar x)\big),\quad \bar p_2=p_2\big(\bar \eta(\bar x)\big),\\
&\bar B=B\big(\bar \eta(\bar x)\big)=-\int^{\bar \eta(\bar x)}_0d_\eta d\cdot h(s)ds\\
&\bar Z_s=Z_s\Big(\f 1k\bar\eta(\bar x)\Big).
\end{split}\]
On the other hand, the shape derivative $d_\eta G(\eta)f\cdot h$ shouldn't depend on the parameter $k$ of course, which means we need to eliminate all the `$k$' in the expression above. 
In fact, direct computations lead to 
\[\begin{split}
&\bar p_2+k p_4=\ga+(1+\ga^2)\Big(k+\f 1k-\f{1}{\ga+\eta'}\Big),\quad k\bar p_2-p_4=k\f{\ga\eta'-1}{\ga+\eta'},\\
&k^2\bar p_1-2k\bar p_2+p_4=k^2\f{1+(\eta')^2}{(\ga+\eta')^2}>0,\quad\hbox{and}\\
&k\bar p_1+(k^2-1)p_2-kp_4=k\f{1+(\eta')^2}{(\ga+\eta')^2}+(1+k^2)\f{\ga\eta'-1}{\ga+\eta'}.
\end{split}\] As a result, one can find from \eqref{Z_c} that  
\[ 
\bar Z_s=\f 1k(1+k^2)^{\f12}\f{\ga+\eta'}{1+(\eta')^2}G(\eta)f+\f 1k(1+k^2)^{\f12}\f{\ga\eta'-1}{1+(\eta')^2}f'+(1+k^2)^{-\f12}\f 1{\ga+\eta'}f'.
\] Moreover, since 
\[
d_\eta d\cdot h(\bar x)=\f{1}{\big(\ga+\eta'(\bar \eta^{-1}(\bar x))\big)^2}
\Big(\f{\eta''(\bar \eta^{-1}(\bar x))}{\ga+\eta'(\bar \eta^{-1}(\bar x))}h(\bar \eta^{-1}(\bar x))+h'(\bar \eta^{-1}(\bar x))\Big),
\] one can express $\bar B$ as 
\beq\label{bar B}
\bar B=\f{h(\bar x)}{\ga+\eta'(\bar x)}-2\int^{\bar x}_0\f{h'(t)}{\ga+\eta'(t)} dt
\eeq while noticing that we take $h(0)=0$ since $\eta(0)=0$ in this paper.

 Summing up the computations above, we finally arrive at the (formal) expression of shape derivative $d_\eta G(\eta)f\cdot h$. Therefore we finish the proof for  Thm \ref{thm:D-N-S}.

\bthm{Remark}
{\it 
Compared with the shape derivative in Thm 3.20 \cite{Lannes}, one can tell that  at least formally the shape derivative  in our paper has a similar form as the shape derivative on a strip domain. This conclusion makes sense since we didn't use any property concerning the corner, and the shape derivative only involves the variation of the upper surface.
}
\ethm

\section{Regularizing Diffeomorphism}
In the end of this paper, we consider about adjusting  transformations $T_S$ and $T_R$ to have a better regularity for $\eta$ (half-order lower regularity indeed). That is, all over our paper, when the norm $|\eta|_{H^{K+2}}$ is used, we show that it can always be replaced by $|\eta|_{H^{K+\f32}}$.

We start with working on $T_S$ first of all. We still denote by $(\bar x,\bar z)$ a point in $\Om$, while by $(x,z)$ a corresponding point in $\cS$ here.  In fact, there are more than one way to find a regularizing transformation. 

For example,  one uses Rmk \ref{one side Dirichlet trace} to define $p(x,z)$  on $\cS$ satisfying the following Dirichlet boundary condition:
\[p(x,z)\big|_{\Ga_t:\,z=kx}=\be(kx)\bar\eta^{-1}(kx) \] where $\be$ is a cut-off function defined on $[0,+\infty)$ and vanish away from $0$. Recall that $\eta(\bar x)=\ga \bar x+\eta(\bar x)$ and $\bar \eta^{-1}(z)$ is the inverse function. 

One can see from Rmk \ref{one side Dirichlet trace} that if $\be\bar \eta^{-1}\in H^{K+\f32}(\R^+)$, then $p(x,z)\in H^{K+2}(\cS)$ with the estimate
\[
\|p\|_{k+2,2}\le C\big(|\eta|_{H^{K+\f32}}\big)|\be \bar \eta^{-1}|_{H^{K+\f32}}\le C\big(|\eta|_{H^{K+\f32}}\big) .
\]  The estimate holds indeed since $|\be \bar \eta^{-1}|_{H^{K+\f32}}$ can be controlled by $|\eta|_{H^{K+\f32}}$.

As a result, we define the regularized transformation $T_S$ as 
\[T_S: \quad (x,\,z)\in S\cap U_{\del S}\mapsto (\bar x,\,\bar z)\in\Om\cap U_\del\]
with
\[\bar x=x+p\big(\eps x+(1-\eps){\f 1k} z,\,z\big)-\f 1k z,\quad \bar z=z-\ga\Big(x+p\big(\eps x+(1-\eps){\f 1k} z,\,z\big)-\f 1k z\Big)\] where   $\eps$ is a small constant to be fixed.
Since a direct computation shows that
\[
\na T_S=\left(\begin{matrix}
\p_x\bar x  &\p_z\bar x\\
\p_x\bar z & \p_z\bar z
\end{matrix}\right)=\left(\begin{matrix}
1+\eps \p_x p  &(1-\eps)\f 1k\p_x p+\p_z p-\f1k\\
-\ga(1+\eps\p_xp) & 1-\ga\big((1-\eps)\f 1k\p_x p+\p_z p-\f1k\big)
\end{matrix}\right)
\] and 
\[Det(\na T_s)=1+\eps\, \p_x p,\]   we know that this new transformation $T_S$ is invertible as long as  the constant $\eps$ is taken  small enough such that 
\[\eps\le \f 1{2\|\p_x p\|_\infty}.\] Besides, one can check that $\Ga_t$, $\Ga_b$ in $\cS$ correspond to $\Ga_t$, $\Ga_b$ in $\Om$ respectively.
 
Some more computations lead to the new matrix 
\[
P_S=(\na T^{-1}_S)\circ T_S=\f 1{1+\eps\,\p_x p}\left(\begin{matrix}
1+\ga\big(\f1k-(1-\eps)\f 1k\p_x p-\p_z p\big) & \ga(1+\eps\,\p_xp)\\
\f1k-(1-\eps)\f 1k\p_x p-\p_z p & 1+\eps\,\p_x p
\end{matrix}\right). \] Therefore, for the  elliptic system of $u$ on $\Om$
\[
\left\{\begin{array}{ll}\Del u= h\qquad \hbox{on}\quad \Om\\
u|_{\Ga_t}=f,\quad \p_{n_b}u|_{\Ga_b}= g
\end{array}\right.
\] the corresponding elliptic system of $v$ on $\cS$ now becomes 
\[
\left\{\begin{array}{ll}\na\cdot \cP_S\na v= h_S\qquad \hbox{on}\quad \cS\\
v|_{\Ga_t}=f\circ T_S,\quad \p^{\cP_S}_{n_b}v|_{\Ga_b}=(1+\ga^2)^\f12  g\circ T_S
\end{array}\right.\] with 
\[\cP_S= P^t_SDet(\na T_S)P_S,\quad h_S=Det(\na T_S)h\circ T_S.\]
When we consider the transformation $T_R$ away from the corner, we can always extend it to a mapping from the infinite flat strip $\{(x,z)|\,0\le z\le 1, x\in \R\}$ to an infinite strip extended  naturally by the part of $\Om$ away from the corner. Therefore, the regularized transformation can be the one in \cite{Lannes} and is omitted here.

Consequently, we could go through all the previous sections using these new transformations. Compared to  old ones, we find that  in all related estimates, the norm $|\eta|_{H^{K+2}}$ can be replaced by the norm $\|p\|_{K+2,2}$, which turns out to retrieve back the norm $|\eta|_{H^{K+\f32}}$.

\bigskip

\section{Appendix}

\subsection{Appendix A: Some known results} To be self-content, we list some results from \cite{PG1,PG2} (sometime slightly adjusted for our paper) here.
We will use the notation $\Ga_j$ as the $j$-th segment (counterclockwise) of the boundary for some polygonal domain $\Om_p$, and $\bf n_j$ as the unit outward vector on $\Ga_j$.
\bthm{Theorem}\label{Thm 1.5.2.3}(Thm 1.5.2.3\cite{PG1}){\it Let $\Om_p$ be a bounded open subset of $\R^2$ whose boundary $\Ga$ is a curvilinear polygon of class $C^1$. Assume that $u|_{\Ga_{j-1}}=u|_{\Ga_{j+1}}=0$ for $u\in H^1(\Om_p)$. Then the mapping $u\mapsto \{f_j\}$, where $f_j=u|_{\Ga_j}$,  is a linear continuous mapping from $H^1(\Om_p)$ onto the subspace of $H^\f12(\Ga_b)$ defined by 
\[\int^{\del_j}_0\f{|f_{j}(\si)|^2}{\si}d\si<\infty\] where $\si$ is the arc length parameter and $\del_j>0$ is some small constant.
}\ethm
\bthm{Theorem}\label{Remark 1.4.4.7}(Rmk 1.4.4.7\cite{PG1}){\it If $u$ belongs to 
$H^\f12(\Om_l)$, then $\na u$ belongs to the dual space $\big(\tilde H^\f12(\Om_l)
\big)^*$ of $\tilde H^\f12(\Om_l)$, where $\Om_l$ is a bounded open subset of $\R^n$ with a Lipschitz boundary.}
\ethm
\bthm{Theorem}\label{Thm 1.5.3.10}(Thm 1.5.3.10\cite{PG1}){\it Let $\Om_p$ be a bounded open set of $\R^2$, whose boundary is curvilinear polygon of class $C^{1,1}$. Then the mapping
\[u\mapsto \p_{n_j}u|_{\Ga_j}\] which is defined on $\mathcal{D}(\bar\Om_p)$ has a  unique continuous extension as an operator from $E(\Del; L^2(\Om_p))=\{u\in H^1(\Om_p)|\,\Del u\in L^2(\Om_p)\}$ into $\big(\tilde H^\f12(\Ga_j)\big)^*$. 
}
\ethm
\bthm{Theorem}\label{Thm 1.5.2}(Thm 1.5.2\cite{PG2}){\it Let $\Om_p$ be a bounded polygonal open subset of $\R^2$. Then the mapping
\[u\mapsto \{u,\,\p_{n_j}u\}|_{\Ga_j}\] which is defined for $u\in H^2(\Om_p)$ has a unique continuous extension as an operator from $D(\Del,L^2(\Om_p))=\{u\in L^2(\Om_p)|\,\Del u\in L^2(\Om_p)\}$ into
\[\big(\tilde H^\f12(\Ga_j)\big)^*\times \big(\tilde H^\f32(\Ga_j)\big)^*.\]
}
\ethm

\subsection{Appendix B: The orthogonal space}
In the beginning of this section, we introduce some more notations for the boundary.  

\noindent -We order $\Ga_t,\Ga_b,\Ga_\del$ counter-clockwisely as $\Ga_1,\Ga_2,\Ga_3$, while the three corner points starting from $(0,0)$ named as $S_1, S_2,S_3$ correspondingly. \\
\noindent -Sometimes we use $\Ga$ to denote the entire boundary $\Ga_1\cup\Ga_2\cup\Ga_3$.  \\
\noindent -We denote by $\al_j$ the angle of the corner $S_j$.\\
\noindent -We denote by $x_j(s)$ the point of $\Ga_j$ which, for $|s|\le \del$ small enough, is at distance $s$ of $S_j$ along $\Ga$. Consequently $x_j(s)\in \Ga_j$ when $s<0$ and $x_{j+1}(s)\in \Ga_{j+1}$ when $s>0$.  \\
\noindent -In this part we say that two functions $\phi_j,\phi_{j+1}$ defined on $\Ga_j, \Ga_{j+1}$ respectively are equivalent at $S_j$ and write 
\[\phi_j\equiv\phi_{j+1}\quad \hbox{at}\quad S_j\] if
\[\int^\del_0\big|\phi_j(x_j(-s))-\phi_{j+1}(x_j(s))\big|ds/s <\infty.\]
We quote a density result Theorem 1.2.8 from \cite{PG2} here, which will be used later.
\bthm{Lemma}\label{dense D}{\it Let $\Om_l$ be a Lipschitz open subset of $\R^n$ then $\cD(\Om_l)$ is dense in $\tilde H^s(\Om_l)$ for all $s\ge0$.}
\ethm
Recalling that 
\[V^2_b=\big\{u\in H^2(\cS_0)\big|\, u|_{\Ga_t}=0,\,\p_{n_b} u+b\p_{\tau_b} u|_{\Ga_b}=0\big\},\] we give a dense result about $V^2_b$ first, which is a modification of Thm 1.6.2\cite{PG2}
\bthm{Lemma}\label{dense}{\it  The space $V^4=H^4(\cS_0)\cap V^2_b$ is dense in $V^2_b$.}\ethm
\begin{proof}This proof is similar as Thm 1.6.2\cite{PG2}. In fact, it's enough to show that  for any linear functional $L$ on $V^2_b$ satisfying $L(H^4\cap V^2_b)=0$,  then $L=0$ on $V^2_b$.  We simply write $V^2$ for $V^2_b$ in this proof.

Since one has the decomposition 
\[H^2(\cS_0)=\dot H^2(\cS_0)\oplus Z^2(\Ga),\quad\hbox{with}\quad Z^2(\Ga)=\big\{\p^l_{n_j}|_{\Ga_j}\big\}_{1\le j\le 3,0\le l\le 1}(H^2)\] we can decompose $L$ acting on any $v\in V^2_b$ as follows
\[\langle L,\,v\rangle=\langle G,\,v- \rho( v|_\Ga)\rangle+\langle g,\,v|_\Ga\rangle\] where $\rho$ is the right inverse of the trace operator, $G\in (\dot H^2)^*$ and $g\in Z^2(\Ga)^*$.

When $L$ vanishes on $H^4\cap V^2$, $L$ also vanishes on $\cD$, so one finds $\langle G,\, v\rangle=0$ for any $v\in \cD$, which implies that $G=0$. Therefore, we have
\[\langle L,\,v\rangle=\langle g,\,v|_\Ga\rangle\]  which implies that $L$ depends on the boundary only. 

Now it remains to show that  $Z^2(\Ga)\cap V^4$ is dense in $Z^2(\Ga)\cap V^2$.
In order to do this, we need descriptions about $Z^2(\Ga)\cap V^4$ and $Z^2(\Ga)\cap V^2$.

 Indeed, $Z^2(\Ga)\cap V^2$ is the image of $V^2$ by the mapping
\[
u\mapsto \{u|_{\Ga_j}=g_j,\,\p_{n_j}u|_{\Ga_j}=h_j\}_{1\le j\le 3}.
\] 
Thanks to Thm \ref{trace thm PG} and Remark \ref{oblique case}, $Z^2(\Ga)\cap V^2$ is the subspace of $\prod_j H^\f32(\Ga_j)\times H^\f12(\Ga_j)$ satisfying the following conditions
\beq\label{common condition}\begin{split}
&g_j=0,\quad j=1,3\\
& h_j+b g'_j=0,\quad j=2\\
&g_j(S_j)=g_{j+1}(S_j),\quad j=1,2,3\end{split}\eeq and 
\beno
&&g'_j\equiv -\cos\al_j g'_{j+1}+\sin\al_jh_{j+1}\quad\hbox{at}\quad S_j,\\
&&h_j\equiv-\cos\al_j h_{j+1}-\sin \al_jg'_{j+1}\quad\hbox{at}\quad S_j,\,\,j=1,2,3.
\eeno 
Similarly, $Z^2(\Ga)\cap V^4$ is the image of $V^4$ by the same mapping
\[
u\mapsto \{u|_{\Ga_j}=g_j,\,\p_{n_j}u|_{\Ga_j}=h_j\}_{1\le j\le 3},
\]
Applying again  Thm \ref{trace thm PG} and Remark \ref{oblique case}, one finds $Z^2(\Ga)\cap V^4$ as the subspace of $\prod_j H^\f72(\Ga_j)\times H^\f52(\Ga_j)$ satisfying \eqref{common condition} and the following conditions for $j=1,2,3$:
\beno
&&g'_j(S_j)= -\cos\al_j g'_{j+1}(S_j)+\sin\al_jh_{j+1}(S_j),\\
&&h_j(S_j)=-\cos\al_j h_{j+1}(S_j)-\sin \al_jg'_{j+1}(S_j),\\
&&-\cos\al_jg''_j(S_j)-\sin\al_jh'_j(S_j)=-\cos\al_jg''_{j+1}(S_j)+\sin\al_j h'_{j+1}(S_j).
\eeno Now we check these conditions for each corner point to conclude the proof. Firstly, for $j=3$, we find $\al_3\in (0,\pi/2)$ and we have Dirichlet conditions on both sides of $\al_3$, which implies that 
\[g_1=g_2=0.\] Consequently, the conditions for $Z^2(\Ga)\cap V^2$ read
\[h_1\equiv 0,\,h_3\equiv-\cos\al_3 h_1\quad\hbox{at}\quad S_3,\] which means $h_1,h_3\equiv 0$ at $S_3$. On the other hand, the conditions for $Z^2(\Ga)\cap V^4$ become
\[h_1(S_3)=0,\,h_3(S_3)=-\cos\al_3h_1(S_3)=0\quad\hbox{and}\quad -\sin\al_3h'_3(S_3)=\sin\al_3h'_1(S_3).\] So the desired density near $S_3$ follows from density of $\cD(\R^+)$ in $\tilde H^\f12(\R^+)$ from Lem \ref{dense D}.

Secondly, when $j=1$, one has $\al_1=\om\in (0,\pi/2)$ and the boundary conditions are mixed type: 
\[
g_1=0,\quad h_2=-bg'_2,\quad\hbox{ with}\quad b=\tan\Phi.
\]  In this case the conditions for $Z^2(\Ga)\cap V^2$ read $g_2(S_1)=0$ and 
\[0\equiv -(\cos\al_1+b\sin\al_1)g'_2,\quad h_1\equiv(b\cos\al_1-\sin\al_1)g'_2\quad\hbox{at}\ \ S_1,\] which implies that 
\[g'_2\equiv0,\quad h_1\equiv 0\quad\hbox{at}\quad S_1.\] On the other hand, the conditions for $Z^2(\Ga)\cap V^4$ rewrite as  $g_2(S_1)=0$ and \[g'_2(S_1)=0,\quad h_1(S_1)=0,\quad-\sin\al_1h'_1(S_1)=-\f1{\cos(k\om)}\cos(k\om+\om) g''_2(S_1),\] therefore the density near $S_1$ can be proved by the density of $\cD(\R^+)$ in $\tilde H^\f12(\R^+)$ and $\tilde H^\f32(\R^+)$. 

For the last case $j=2$, the argument is similar as in the second case.  The proof is finished as a result.
\end{proof}

Now we can give a description of the space $\cN_b$.
\bthm{Lemma}\label{space N}{\it When $\Phi-\om\neq l\pi+\f\pi2$  for some $l\in \Z$, a function $w$ is in the orthogonal space $\cN_b$ of $\Del(V^2_b)$ in $L^2(\cS_0)$ iff $w\in L^2(\cS_0)$ and satisfies
\[\left\{\begin{array}{ll}\Del w=0\qquad\hbox{on}\quad \cS_0\\
w|_{\Ga_t}=0,\quad w|_{\Ga_\del}=0,\quad \p_nw-b\p_\tau w|_{\Ga_b}=0.\end{array}\right.\]  }\ethm
\begin{proof} This proof is a modification of Thm 2.3.3\cite{PG2}. In fact, it suffices to show that when $w$ satisfies the system above, $w$ must be in $\cN_b$, which means that $\int_{\cS_0}w\Del u=0$ for any $u\in V^2_b$. The reverse statement is proved already in Lem 4.4.1.4\cite{PG1}. 

Indeed it's enough to consider $u\in V^4$ from Lem \ref{dense}, and in this case we have $u\in C^2(\bar \cS_0)$.
First of all, since $u|_{\Ga_t}=u|_{\Ga_\del}=0$, we know that $u(S_j)=0$ for $j=1,2,3$. Moreover, taking the tangential derivative along $\Ga_t,\Ga_\del$,  we find 
\[
\na_{\tau_t} u|_{\Ga_t}=\na_{\tau_\del} u|_{\Ga_\del}=0,
\]
 where $\tau_\del$ is the unit tangential vector on $\Ga_\del$. So we have $\na u(S_3)=0$ directly. On the other hand,  remember that $u$ satisfies   
\[
(\p_{n_b}+b\p_{\tau_b})u|_{\Ga_b}=0,\quad\hbox{ with}\quad b=\tan\Phi,
\]
so  when $\Phi-\om\neq l\pi+\f\pi2$,  we find the direction ${\bf \mu}_b={\bf n}_b+b{\bf \tau}_t\nparallel \tau_t$ , which leads to $\na u(S_1)=0$. 
 
 Similarly we also find that $\na u(S_2)=0$, since the boundary $\Ga_\del$ can be changed from the beginning if necessary (which means $\al_2$ can be changed to avoid the parallel case, recalling that when we consider the triangle $\cS_0$, we always focus on some function supported near $S_1$). Combing with the fact that $u\in C^2(\bar\cS_0)$, we arrive at the conclusion
\beq\label{trace space}u|_{\Ga_j}\in\tilde H^\f32(\Ga_j),\quad \p_{n_j}u|_{\Ga_j}\in \tilde H^\f12(\Ga_j)\quad j=1,2,3.\eeq

Secondly, applying the Green Formula for $w\in L^2(\cS_0)$ satisfying the system in this lemma and $u\in V^4$ one has
\[\int_{\cS_0} \Del u \,w-\int_{\cS_0}u\Del w=\sum^3_{j=1}\Big(\int_{\Ga_j}\p_{n_j}u\,w-\int_{\Ga_j} u\p_{n_j} w\Big),\] where the right side holds in the sense  of \eqref{trace space} and 
\[w|_{\Ga_j}\in (\tilde H^\f12(\Ga_j))^*,\quad \p_{n_j} w|_{\Ga_j}\in (\tilde H^\f32(\Ga_j))^*.\] Plugging the conditions of $u,\,w$ into the equality above, one derives that 
\[\int_{\cS_0}\Del u\, w=\int_{\Ga_b}\p_n u\,w-\int_{\Ga_b}u\p_nw=\int_{\Ga_b} u(\p_nw-b\p_\tau w)=0\] where $\int_{\Ga_b}\p_n u\,w=-b\int_{\Ga_b}\p_\tau u \,w=b\int_{\Ga_b}u\p_\tau w$ makes sense since $\p_\tau u\in \tilde H^\f12(\Ga_b)$. In this way the proof is done.\end{proof}
\bigskip

\noindent{\bf Acknowledgement}.   The authors would like to thank for  kind supports from Zhifei Zhang and Chongchun Zeng. The author Mei Ming is supported by NSFC no.11401598.

\end{document}